\documentclass[reqno,1pt]{amsart}
\usepackage{enumerate}
\usepackage{cases}
\usepackage{mathrsfs}
\usepackage[T1]{fontenc}
\usepackage{mathrsfs}
\usepackage{amsmath,latexsym,amssymb,amsfonts,amsbsy, amsthm}
\usepackage[colorlinks=true, linkcolor=blue, citecolor=green]{hyperref} 
\usepackage{bm}
\usepackage[usenames]{color}
\usepackage{xspace,colortbl}
\usepackage{epsfig}
\usepackage{graphicx}
\usepackage{subfigure}
\usepackage{amsmath,amsfonts,amsthm,amssymb,amscd}
\usepackage{color}
\usepackage[title,titletoc,toc]{appendix}
\usepackage{bigints}
\allowdisplaybreaks[4]

\newcommand{\R }{\mathbb{R}}

\newcommand{\al}{\alpha}

\newcommand{\e}{\varepsilon}

\newcommand{\N}{\mathcal{N}}

\newcommand{\la}{\lambda}

\newcommand{\ml}{\mathcal{L}}
\newcommand{\lxi}{\mathcal{L}_\xi}

 \renewcommand{\Re}{\operatorname{Re}}

 \newtheorem{theorem}{Theorem}[section]
 \newtheorem{lemma}[theorem]{Lemma} 
 \newtheorem{prop}[theorem]{Proposition} 
 \newtheorem{remark}[theorem]{Remark}
\newtheorem{definition}[theorem]{Definition}
\newtheorem{cor}[theorem]{Corollary}
\def\theequation{\arabic{section}.\arabic{equation}}

\title[Stability of periodic waves in the FHN equation]
{Nonlinear Stability of Large-Period Traveling Waves Bifurcating from the Heteroclinic Loop in the FitzHugh-Nagumo Equation}

\author[Ji Li, Ke Wang, Qiliang Wu and Qing Yu ]{}

\thanks{Ji Li is supported by the National Natural Science Foundation of China (Grant No. 12171174).
Qing Yu is supported by the National Natural Science Foundation of China (Grant No. 12401202)}
\thanks{e-mail:  liji@hust.edu.cn(Li); ke\_wang@hust.edu.cn(Wang); wuq@ohio.edu(Wu);  yuqing@hust.edu.cn(Yu)}

\begin{document}
\maketitle
\centerline{\scshape Ji Li$^a$, Ke Wang$^a$, Qiliang Wu$^{b,*}$ and Qing Yu$^{a}$
}
\medskip
{\footnotesize
 \centerline{$^a$ School of Mathematics and Statistics,  Huazhong University of Science and Technology,  Wuhan 430074, China}
 \centerline{$^b$ Department of Mathematics, Ohio University, Athens, OH 45701, USA}
}

\medskip
\begin{abstract}
{ A wave front and a wave back that spontaneously connect two hyperbolic equilibria, known as a heteroclinic wave loop, give rise to periodic waves with arbitrarily large spatial periods through the heteroclinic bifurcation. The nonlinear stability of these periodic waves is established in the setting of the FitzHugh-Nagumo equation, which is a well-known reaction-diffusion model with degenerate diffusion. 
First, for general systems, we give the expressions of spectra with small modulus for linearized operators about these periodic waves via the Lyapunov-Schmidt reduction and the Lin-Sandstede method. Second,  applying these spectral results to the FitzHugh-Nagumo equation, we establish their diffusive spectral stability. Finally, we consider the nonlinear stability of these periodic waves against localized perturbations. We introduce a spatiotemporal phase modulation $\varphi$, and couple it with the associated modulated perturbation $\mathbf{V}$ along with the unmodulated perturbation $\mathbf{\widetilde{V}}$ to close a nonlinear iteration argument. 
}

\vspace{0.1in}


\vspace{0.1in}
\end{abstract}

\section{Introduction.}\setcounter{equation}{0}

The FitzHugh-Nagumo (FHN) equation\cite{FHNF, FHNN}
\begin{equation}\label{FHN}
\begin{cases}
u_t = u_{xx}+f(u)-w,\\
w_t = \epsilon (u-\gamma w),
\end{cases}
\end{equation}
where $u(t,x)$ and $w(t,x)$ are scalar functions defined on $\R^+\times\R$, \( f(u) = u(1 - u)(u - a) \) is a  cubic polynomial, $a\in(0,\frac{1}{2})$, $\gamma>0$ and $0<\epsilon\ll1$ are three parameters, is a well-known  phenomenological model for studying excitation and propagation of excitable media far from equilibrium. Originally proposed in the 1960s as a simplification of, the FHN equation \eqref{FHN} has found its application in many fields such as neuroscience \cite{BFOK,Izhi}, cardiac physiology\cite{ABE,LF..}, electronic and optical systems \cite{RAFBJ, SBBSKB} among others; see also \cite{CPRG} and the reference therein.

In addition to its wide applications, the FHN equation is also well-known for its rich dynamic structure as a pattern forming system. Taking traveling waves as an example, the FHN equation \eqref{FHN} accommodates simple pulses, monotone wavefronts, multimodal wavefronts and pulses, pulses with oscillatory tails, spatial periodic waves of small or large amplitudes, and invasion wavefront from spatial periodic to spatial homogeneous states, to name just a few.

Two typical mathematical properties of the FHN motivate this study. On the one hand, the corresponding traveling wave equation \eqref{TWode} admits a heteroclinic loop and a generic heteroclinic bifurcation which generate periodic orbits of large amplitude and period. On the other hand, the PDE \eqref{FHN} is a two-component system of reaction-diffusion equations, which is degenerate in the sense that the diffusion coefficient of $w$ is zero so that semigroup estimates differ from the classical semilinear case.  

We shall focus on spatial periodic waves of large amplitude. For both of excitable and oscillatory regimes, there exist spatial periodic traveling waves. For period-related patterns in the oscillatory case, we refer to the work \cite{CS} and \cite{ACdS}. We study the excitable case, that is, \eqref{FHN} with $0<a<\frac{1}{2}$. In this regime, $(u,w) = (0, 0)$ is a stable spatially constant equilibrium. If $\gamma>0$ is small, $(0,0)$ is the only spatially constant equilibrium. And if $\gamma$ is not too small, there are two more equilibria, one of which is unstable and the other stable.

Introducing the vector notation $\mathbf{U}:=(u,w)^t$, a traveling wave solution with speed $c$ in \eqref{FHN} takes the form
\[
{
\bf U
}(t,x)={\bf U}(x+ct).
\]
For convenience, we rewrite \eqref{FHN} under the co-moving frame \((\Bar{x}, t) = (x + ct, t)\), yielding 
\begin{equation}\label{FHNtr}
\begin{aligned}
u_t &= u_{xx} - cu_x + f(u) - w,\\
w_t &= -cw_x + \epsilon (u - \gamma w),
\end{aligned}
\end{equation}
where we abuse the notation of $x$ by dropping the bar sign. traveling wave solutions of \eqref{FHN} become steady states of \eqref{FHNtr}; that is, traveling wave profiles satisfy the following ODE system 
\begin{equation} \label{TWode}
\begin{aligned}
u' &= v, \\
v' &= c v - f(u)  + w , \\
w' &= \frac{\epsilon}{c} (u - \gamma w),
\end{aligned}
\end{equation}
where \( ' = \frac{d}{dx} \). The slow-fast feature of \eqref{TWode}, i.e. $0<\epsilon\ll 1$, enables construction of homoclinic or heteroclinic solutions. In fact, various non-degeneracy conditions for the onset of homoclinic or heteroclinic bifurcations  were verified; see for instance \cite{CDT, Den91a, Den91b, Has, SA}. Consequently, periodic solutions (of large amplitude and period) exist for nearby parameter values.

To study the Lyapunov stability of a periodic wave $\mathbf{\bar{U}}(\cdot)=(\phi(x),\psi(x))^t$, we consider the evolution of the perturbation $\mathbf{\widetilde{V}}(t,x)=\mathbf{U}(t,x)-\mathbf{\bar{U}}(x)$, where $\mathbf{U}(t,x)$ is a solution of \eqref{FHNtr}, whose initial data is close to $\bar{\bf U}$.
The perturbation $\mathbf{\widetilde{V}}$ satisfies the following semilinear equation
\begin{equation}\label{semi}
(\partial_t-\ml)\mathbf{\widetilde{V}}=\tilde{\mathcal{N}}(\mathbf{\widetilde{V}}),
\end{equation}
where the nonlinear term $\tilde{\mathcal{N}}$ contains no derivative of $\mathbf{\widetilde{V}}$, and the differential operator $\mathcal{L}$ with periodic coefficients is from the linearization of \eqref{FHNtr} at the steady state $\mathbf{\bar{U}}$; that is,
\[
\begin{matrix}
\mathcal{L}:& H^{2}(\mathbb{R})\times H^{1}(\mathbb{R})\subset L^2 (\mathbb{R})\times L^2(\mathbb{R}) &\longrightarrow & L^2 (\mathbb{R})\times L^2(\mathbb{R}) \\
&{\bf V} & \longmapsto &\left(\begin{array}{cc}
\partial_{x}^2-c\partial_x+f'(\phi)&-1\\
\e&-c\partial_x-\e\gamma
\end{array}\right){\bf V}.
\end{matrix}
\]
The spectrum of $\mathcal{L}$ plays a significant role in determining the stability of $\bar{\bf U}$.

\bigskip

\paragraph{\bf Notation} 


 

$\mathbb{N}_{+}$ is the set of all positive integers, while $\mathbb{N}_{0}$ represents the set of all nonnegative integers. The notation $``A\lesssim B "$ means that there exists a constant $C>0$, independent of $A$ and $B$ such that $A\leq CB$. 
Assume $X$ and $Y$ are Banach spaces. $\mathbf{L}(X,Y)$ represents a linear map from $X$ to $Y$. Let $g=\left(\begin{array}{cc}g_1\\g_2\end{array}\right)$ . We say $g\in X\times Y$ if $g_1\in X$ and $g_2\in Y$, and the norm is defined as $\|g\|_{X\times Y}:=(\|g_1\|_X^2+\|g_2\|_Y^2)^{1/2}$. Furthermore, $\langle\cdot, \cdot \rangle_{L^2}$ represents the complex inner product on $L^2$, that is, $\langle f, g\rangle_{L^2}=\int f(x)\cdot \overline{g(x)}dx$, while $\langle\cdot, \cdot \rangle$ represents the scalar product in $\mathbb{R}^n$. 
  
 \bigskip

\subsection{Diffusive spectral stability} 
 
Since $\mathcal{L}$ has periodic coefficients, standard Floquet-Bloch theory implies that its spectrum as an operator acting on $L^2 (\mathbb{R})$ is entirely essential 
and, thanks to the translation invariance of \eqref{FHNtr}, necessarily touches the imaginary axis at the origin. 
This makes the stability analysis for spatially periodic patterns significantly different and more challenging than for wavefronts or pulses, as for the latter ones the translational eigenvalue at the origin can generally be separated from the rest of the spectrum. Nevertheless, in the setting of localized perturbations for spatially periodic patterns, the best one can hope for is that the spectrum is confined to the open left half-plane except for a single critical spectral curve generically touching the origin with a quadratic tangency, which naturally leads to the following definition of spectral diffusive stability\cite{Sch}. 

\begin{definition}
A smooth $T$-periodic steady state solution $\mathbf{\bar{U}}(\cdot)=(\phi(x),\psi(x))^t$ of \eqref{FHNtr} is said to be diffusively spectrally stable provided the following conditions hold:
\begin{enumerate} [(1).]
    \item the $L^2(\mathbb{R})$-spectrum of corresponding linear operator $\mathcal{L}$ satisfies
    $$\sigma(\ml)\subset \{\la\in\mathbb{C}:\Re \la <0\}\cup \{0\}.$$
    \end{enumerate}
    \begin{enumerate}[(2).]
    \item there exists a constant $\theta>0$ such that the $L^2_{per}(0,T)$-spectrum of the Bloch operator $\lxi:=e^{-i\xi x}\ml e^{i\xi x}, \xi \in [-\pi/T,\pi/T)$, acting on $L^2_{per}(0,T)$, satisfies
    $$\Re \sigma(\lxi)\leq -\theta \xi^2,\quad  \xi \in [-\pi/T,\pi/T).$$
    \end{enumerate}
    \begin{enumerate}[(3).]
    \item $\la=0$ is a simple eigenvalue of $\ml_0$ with the eigenfunction given by the derivative of the underlying periodic solution.
\end{enumerate}
\end{definition}

Periodic waves with such a diffusive stability occur in many classical pattern forming models such as the
Ginzburg-Landau equation \cite{van95}, the Gierer-Meinhardt system \cite{VD05}, the Brusslator model \cite{SZJV18}, the Swift-Hohenberg equation \cite{Ue99} and the Lugiato-Lefever equation \cite{DH18,HJPd23}, just to name a few. 

The diffusive spectral stability of periodic waves from homoclinic bifurcation in \eqref{FHNtr} was established in \cite{SA} using the Lin-Sandstede method and the Lyapunov-Schmidt reduction. 
Regarding the bifurcation problems of heteroclinic loops, there have already been many significant and insightful research results. In particular, Deng \cite{ CDT,Den91a} investigated the existence of simple fronts (backs), N-fronts (backs), simple pulses, and periodic solutions bifurcating from heteroclinic loops under certain assumptions.  In addition, Deng \cite{Den91b} studied the existence of infinitely many traveling fronts and backs bifurcating from heteroclinic loops in the FitzHugh-Nagumo  equations, which exhibit a double-twisted structure. 
Other examples include (but are not limited to) the Circadian model \cite{LY}, fluid dynamics \cite{EKR}, coupled arrays of Chua's circuits \cite{GLL}, and the dryland vegetation model \cite{BCD}. For further studies on homoclinic and heteroclinic bifurcations, see \cite{Has}. To the best of our knowledge, the stability of periodic solutions bifurcating from heteroclinic loops has not yet been studied in the literature. The arguments used in \cite{SA} to obtain diffusive spectral stability apply naturally to the case of heteroclinic bifurcation. However, the explicit expression of the critical spectral curve touching the origin is more difficult to obtain. Nevertheless, we are able to generalize those steps in \cite{SA} in the context of heteroclinic bifurcatoin for general reaction-diffusion systems in Section \ref{Sspec}.

\bigskip

\subsection{Nonlinear stability}
There have been mathematically rigorous proofs of nonlinear stability of periodic wave solutions with diffusive spectral stability in various reaction-diffusion equations; see, for instance, \cite{DSSS, JNRZ, JZ11}. However, the degenerate diffusive matrix in the FitzHugh-Nagumo equation causes technical difficulties in showing the spectral mapping properties of the corresponding linear semigroup. In particular, the corresponding semigroup $e^{\ml t}$ is not analytic. As a result, it does not have smooth effects as an analytic semigroup does. This makes the verification of non-linear stability more complicated.

In this work, we establish the following main result on the nonlinear stability of large-period traveling waves bifurcating from the heteroclinic loop in the FHN \eqref{FHN} against localized perturbations. 

\begin{theorem}\label{nonsta}
There exist constants $M,\e_0>0$ such that, whenever $\mathbf{{V}}_0\in (L^1(\mathbb{R})\times L^1(\mathbb{R}))\cap (H^5(\mathbb{R})\times H^4(\mathbb{R}))$ satisfies
$$E_0:=\|\mathbf{{V}}_0\|_{(L^1(\mathbb{R})\times L^1(\mathbb{R}))\cap (H^4(\mathbb{R})\times H^4(\mathbb{R}))}<\e_0,$$
there exist functions
\begin{align*}
\mathbf{\widetilde{V}}&\in C([0,\infty),H^5(\mathbb{R})\times H^4(\mathbb{R}))\cap C^1([0,\infty),H^3(\mathbb{R})\times H^3(\mathbb{R})),\\
\varphi&\in C([0,\infty),H^3(\mathbb{R}))\cap C^1([0,\infty),H^1(\mathbb{R})),
\end{align*}
satisfying $\mathbf{\widetilde{V}}(0)=\mathbf{{V}}_0$ and $\varphi(0)=0$ such that $\mathbf{U}(t)=\mathbf{\bar{U}}+\mathbf{\widetilde{V}}(t)$ is the unique global solution of \eqref{FHNtr} with initial date $\mathbf{\bar{U}}+\mathbf{{V}}_0$, and we have the following estimates
\begin{align*}
\max\{\|\mathbf{U}(t)-\mathbf{\bar{U}}\|_{H^4},\|\varphi(t)\|_{L^2}\}\leq ME_0(1+t)^{-\frac{1}{4}}, \text{ for all }t\geq0,
\end{align*}
and
\begin{align*}
\max\{\|\mathbf{U}(\cdot-\varphi(t,\cdot),t)-\mathbf{\bar{U}}\|_{L^2}, \|\partial_x\varphi(t)\|_{H^3}, \|\partial_t\varphi(t)\|_{H^2}\}\leq ME_0(1+t)^{-\frac{3}{4}}, \text{ for all }t\geq0.
\end{align*}
\end{theorem}
For the norm of $g$ where $g\in H^m\times H^m$,  we consistently use the notation $\|g\|_{H^m}$ and simply for other Sobolev spces throughout this paper.

\bigskip

\subsection{Outline of paper and sketch of proofs}
The outline of the rest of the paper is as follows.
In Section \ref{Sspec} we analyze the spectrum, and in Section \ref{stab} we prove the main theorem on nonlinear stability. We sketch our main steps in the rest of this section. 

\subsubsection{\bf Spectrum analysis}\label{SSspec}
We investigate the existence and spectral stability of periodic solutions bifurcating from heteroclinic loops for general reaction-diffusion systems using the Lyapunov-Schmidt reduction method and the Lin-Sandstede method. 
We consider an abstract reaction-diffusion equation of the form
\begin{equation} \label{102}
    U_t = DU_{xx} +cU_x+ g(U),
\end{equation}
where \(U(t,x) \in \mathbb{R}^N\), \(x \in \mathbb{R}, t \ge 0\), and $D$ is a (non-strictly) positive definite diagonal matrix, with the first $k$ diagocal elements nonzero.  
Setting \( u =(u_1,u_2):= (U, U_{x}) \),  \( U \) is a bounded steady-state, i.e., $U$ is bounded and satisfies 
\[DU_{xx} + c U_{x} + g(U) = 0, \quad U \in \mathbb{R}^N, \]
if and only if, \( u \) is a bounded solution of the system
\begin{equation} \label{104}
u' = F(u):=\begin{pmatrix}
    u_2 \\
    -cD^{-1}u_2-D^{-1}g(u_1)
\end{pmatrix}, \quad u \in \mathbb{R}^{N+k}, 
\end{equation}
where \( ' = \frac{d}{dx} \).

Assumming there is a $T$-periodic steady state \(P(x)\) of \eqref{102}, the linearized operator at \(P(x)\) is denoted as 
\begin{equation} \label{103}
 \mathcal{L} := \partial_{xx} + c \partial_{x} + g_U(P).
\end{equation}
The $L^2(\mathbb{R})$-spectrum of $\mathcal{L}$ consists entirely of its essential spectrum. 
Moreover, the eigenvalue problem
\[\mathcal{L} V = V_{x x} + c V_{x} + g_U(P) V = \lambda V,\quad \lambda\in\mathbb{C},\]
is equivalently written for \(v = (V, V_{x})\) as
\begin{equation} \label{105}
v' = (F_u(p(x)) + \lambda B) v,
\end{equation}
where \(F_u\) is the Jacobian matrix of \(F\), \(p := (P, P_{x})\), and the matrix $B$ is given as
\[
B:=\begin{pmatrix}
    0&0\\
    D^{-1}&0
\end{pmatrix}.
\]

According to Floquet theory, \(\lambda \in \mathbb{C}\) is in the $(L^2(\mathbb{R}))^N$-spectrum of the linearization \(\mathcal{L}\) if and only if there is a \(\xi \in [-\pi/T,\pi/T)\) and a bounded function \(v\) such that
\begin{equation} \label{107}
\begin{cases}
v' = (f_u(p(x)) + \lambda B)v, \quad \text{for} \quad 0<x < T, \\
v(T) = e^{i\xi T}v(0).
\end{cases}
\end{equation}

In line with the case of heteroclinic bifurcation, we consider systems of the form \eqref{104} with at least two parameters; that is, 
\begin{equation}\label{108}
u' = F(u,\mu), \quad (u, \mu) \in \mathbb{R}^n \times \mathbb{R}^2,
\end{equation}
where \( F \) is at least \( C^2 \). We assume that for $\mu=0$ there is a heteroclinic loop consisting of two heteroclinic orbits $h_1$ and $h_2$. We shall prove the existence of periodic orbits in \eqref{108} for $\mu$ near $0$ and  denote as
\(Q(x)\) the $T$-periodic orbit for $\mu=\mu_T$. Also for convenience, assume \(T=2(L_1+L_2)>0\), where we use $2L_1$ and $2L_2$ to represent the time the periodic orbit spends near $h_1$ and $h_2$, respectively, as the period $T$ becomes arbitrarily large.

The corresponding spectrum problem \eqref{107} then equivalently reads
\begin{equation} \label{109}
\begin{aligned}
&v' = (F_u(q(x),\mu_T) + \lambda B) v, \quad  -L_1 < x < L_1 + 2L_2 \\
&v(L_1+ 2L_2) = e^{i \xi T} v(-L_1)
\end{aligned}
\end{equation}
where \(\xi \in [-\pi/T,\pi/T)\),  \(\lambda \in \mathbb{C}\) and $q:=(Q,Q^\prime)$.

The main result of this part is an expansion of the curve of critical eigenvalues of the large spatial periodic pattern in terms of its period. The critical eigenvalue about the \( T = 2(L_1+L_2) \)-periodic wave satisfies, to leading order,
\begin{equation} \label{110}
\begin{aligned}
0=E(\lambda, \xi) = & (1-e^{i \xi T}) \langle \psi_{2}(L_2), h'_1(-L_1) \rangle\langle \psi_{1}(L_1), h_2'(-L_2) \rangle \\
& + (1-e^{-i \xi T})\langle \psi_{1}(-L_1), h_2'(L_2) \rangle\langle \psi_{2}(-L_2), h'_{1}(L_1) \rangle \\
& - \left( \langle \psi_{2}(-L_2), h'_{1}(L_1) \rangle - \langle \psi_{2}(L_2), h'_1(-L_1) \rangle \right) \lambda M_1 \\
& - \left(\langle \psi_{1}(-L_1), h_2'(L_2) \rangle - \langle \psi_{1}(L_1), h_2'(-L_2) \rangle \right) \lambda M_2.  \\
\end{aligned}
\end{equation}
In the above, the nonzero constant $M_i$ is given by the Melnikov-type integral
$$M_i:= \int_{-\infty}^{\infty} \langle \psi_i(x), Bh_i^{'}(x) \rangle dx , \quad i = 1, 2,$$
where $\psi_i(x)$ is the unique bounded solution (up to constant multiples) of the adjoint variational equation
\begin{equation*}
w'= -F_u(h_i(x), 0)^* w, \quad x \in \mathbb{R}.
\end{equation*}
From the expression \eqref{110}, we can see that the spectrum of the periodic solution $q(x)$
 near $0$ is closely related to the decay properties of the heteroclinic orbit $h_i(x)$
as $x$ tends to infinity. As the period $T=2(L_1+ L_2)$ increases to infinity, the critical spectral curve shrinks to $0$. Physically, the inner product terms in the expression describe the interactions of wave solutions near the system boundaries and the characteristics of energy transfer, where $\psi_{i}(x)$ represents the system's response to perturbations, and $h'_i(x)$ describes the rate of change of the heteroclinic orbit $h_i(x)$. For example, 
$\langle \psi_{1}(-L_1), h_2'(L_2) \rangle$ describes how a perturbation at position $-L_1$
affects changes in the system at position $L_2$ , while $\langle \psi_{2}(L_2), h_1'(-L_1) \rangle$  describes the reverse process. This symmetry ensures that as the wave travels through one period, the energy transfer within the system remains balanced.

\bigskip

\subsubsection{\bf Nonlinear stability} 
To study the FHN \eqref{FHN}, we first apply the results of Section \ref{Sspec} to show that the critical eigenvalues near 0 admits the following expansion
$$\la_c(\xi)=ib\xi-d\xi^2+O(|\xi|^3),\quad  |\xi|\ll 1,$$
for some $b\in\mathbb{R}$ and $d>0$, which establishes the spectral stability.

Linear decay estimates are then established via the Gearhart-Pr\"{u}ss theorem.
It turns out that the semigroup $e^{\ml t}$  decays at the rate $(1+t)^{-1/4}$. 
The unmodulated perturbation $\mathbf{\widetilde{V}}=\mathbf{U}(t,x)-\mathbf{\bar{ U}}(x)$, which satisfies the semilinear equation \eqref{semi}, decays in $H^m(\mathbb{R}),m\in\mathbb{N}_0$ at most with rate $(1+t)^{-1/4}$. We are unable to close a nonlinear iteration scheme using this slow decay; see Remark \ref{unv}. Importantly, however, it is worth noting that we can establish the nonlinear damping estimate for the unmodulated perturbation $\mathbf{\widetilde{V}}$ and, as a result, the norm of higher derivatives of $\mathbf{\widetilde{V}}(t)$ can be controlled by $\|\mathbf{\widetilde{V}}(t)\|_{L^2}$ using nonlinear damping estimates. Take $\|\mathbf{\widetilde{V}}(t)\|_{H^4}$ for example. As remarked in Section \ref{Secnon}, one can deduce that
\begin{equation}\label{nonlinearD}
\frac{d}{dt}\|\mathbf{\widetilde{V}}(t)\|^2_{H^4}\leq -\e\gamma\|\mathbf{\widetilde{V}}(t)\|^2_{H^4}+C\|\mathbf{\widetilde{V}}(t)\|_{L^2}^2,
\end{equation}which yields the inequality
\begin{equation*}
\|\mathbf{\widetilde{V}}(t)\|^2_{H^4}\leq e^{-\e\gamma t}\|\mathbf{{V}}_0\|^2_{H^4}+C\int_0^te^{-\e\gamma(t-s)}\|\mathbf{\widetilde{V}}(s)\|^2_{L^2}ds,
\end{equation*}
by the Gronwall's inequality. 

To overcome difficulties induced by the slow $(1+t)^{-1/4}$ decay, we exploit the Bloch wave decomposition to decompose the semigroup $e^{\ml t}$ into an exponentially decaying part $S_e(t)$ and the critical part $S_c(t)$. The critical part $S_c(t)$ consists of a critical component $\mathbf{\bar{U}}(\cdot)s_p(t)$ with decaying rate $(1+t)^{-1/4}$ and another part $\tilde{S}_c(t)$ with faster decaying rate $(1+t)^{-3/4}$.  Moreover, $S_c(t)$ is infinitely smooth as it defines a bounded linear map from $L^1(\mathbb{R})\cap L^2(\mathbb{R})$ to $H^m(\mathbb{R})$ for every $m\in\mathbb{N}_0$ while the exponentially decaying part $S_e(t)$ is not infinitely smoothing. 
Taking advantage of the above decomposition, 
we introduce the spatio-temporal phase modulation $\varphi(t,x)$ and consider the modulated perturbation $\mathbf{V}(t,x)=\mathbf{U}(t,x-\varphi(t,x))-\mathbf{\bar{ U}}(x)$. The phase modulation $\varphi(t,x)$ is chosen to compensate for the slowest decay $(1+t)^{-1/4}$ such that the linear part of $\mathbf{V}$ decays at rate $(1+t)^{-3/4}$. 

The modulated perturbation $\mathbf{V}$ satisfies the following quasilinear equation
\begin{align}\label{quasi}
(\partial_t-\ml)(\mathbf{V}+\varphi\bar{\mathbf{U}}')=\mathcal{N}(\mathbf{V},\varphi,\varphi_t)+(\partial_t-\ml)(\varphi_x\mathbf{V}).
\end{align}
Then, by virtue of the Duhamel formulation and a suitable choice of $\varphi$, we obtain
$$\mathbf{V}(t)=\tilde{S}(t)\mathbf{V}_0+\int_0^t \tilde{S}(t-s) \mathcal{N}(\mathbf{V}(s),\varphi(s),\varphi_t(s))ds+\varphi_x(t)\mathbf{V}(t),$$
where $\tilde{S}(t)$ decays at rate $(1+t)^{-3/4}$. 
This integral equation and the expression of $\varphi$ can be used to control the $L^2$-norm of $\mathbf{V}(t)$, which decays with rate $(1+t)^{-3/4}$. As a result, estimates of $\mathbf{\widetilde{V}}$ can be obtained using the relation $\mathbf{\widetilde{V}}(t,x)-\mathbf{V}(t,x)=\mathbf{U}(t,x)-\mathbf{{U}}(t,x-\varphi(t,x))$ and the mean value theorem.  

However, it is difficult to directly establish a nonlinear damping estimate for the modulated perturbation $\mathbf{V}$ that satisfies the quasilinear equation \eqref{quasi}. More specifically, it is difficult to directly control $\mathbf{V}_x$ and $\mathbf{V}_{xx}$ using the integral equation \eqref{quasi}, since the derivatives, $\mathbf{V}_x$ and $\mathbf{V}_{xx}$, occur in the nonlinear term $\mathcal{N}$ and the exponentially decaying part $S_e(t)$ is not sufficiently smoothing.
To overcome the above difficulties, inspired by the work \cite{HJPd23}, we try to close a nonlinear iteration scheme by coupling $\mathbf{\widetilde{V}},\varphi,\mathbf{V}$ and their derivatives as a system, which, together with the mean value theorem, shows that $\mathbf{V}_x$ and $\mathbf{V}_{xx}$ can be controlled by the unmodulated perturbation $\|\mathbf{\widetilde{V}}\|^2_{H^4}$, the phase modulation $\varphi$ and their derivatives. In addition, $\|\mathbf{\widetilde{V}}\|^2_{H^4}$ can be controlled using the nonlinear damping estimate \eqref{nonlinearD}.
\bigskip

\section{general case: existence and spectral analysis of periodic waves bifurcating from heteroclinic loop}\label{Sspec}   \setcounter{equation}{0}

\subsection{Existence of periodic solutions}

As explained in Section \ref{SSspec}, we consider the system of ordinary differential equations with two bifurcation parameters
\begin{equation}\label{201}
u' = F(u,\mu), \quad (u, \mu) \in \mathbb{R}^n \times \mathbb{R}^2,
\end{equation}
where $\mu = (\mu_1, \mu_2)$, and \( F \) is at least \( C^2 \).
First, we make the following hypotheses.
\vskip 0.1in
\textbf{Hypothesis (H1).}{\it
The system (\ref{201}) has two equilibrium points $e_1(\mu)$ and $e_2(\mu)$, both of which are hyperbolic when $\mu= {\bf 0}$. In other words, denoting $e_i( {\bf 0}):=e_1^*$ for $i=1,2$, we have \( F(e_1^*,  {\bf 0}) = F(e_2^*, {\bf 0})=  {\bf 0} \) and there exist constants \(\alpha_i^s, \alpha_i^u > 0\) such that for any $\nu\in\sigma (F_u(e_i^*, 0))$, \( \Re \nu < -\alpha_i^s \) or \( \Re \nu > \alpha_i^u \).
Furthermore, we set $\dim W^s(e_1^*, {\bf 0})=\dim W^s(e_2^*,  {\bf 0})$, where
$W^s(e_i(\mu),\mu)$ and $W^u(e_i(\mu),\mu)$ represent the stable and unstable manifolds of the equilibrium $e_i(\mu) (i=1,2)$, respectively.}

\vskip 0.1in

\textbf{Hypothesis (H2).}{\it
The system (\ref{201}) possesses a heteroclinic loop when \(\mu =  {\bf 0}\), composed of two heteroclinic solutions \(h_1(x)\) and \(h_2(x)\), satisfying
\[\lim_{x\to -\infty} h_1(x)=\lim_{x\to+\infty}h_2(x)=e_1^*, \quad \lim_{x\to+\infty}h_1(x)=\lim_{x\to-\infty} h_2(x)=e_2^*.\]
}
\vskip 0.1in

\textbf{Hypothesis (H3).} {\it The unique nontrivial bounded solution of the variational equation
\begin{equation}\label{202}
v' = F_u(h_i(x),  {\bf 0}) v, \quad x \in \mathbb{R},
\end{equation}
is given by \(h_i'(x)\) for $i=1,2$, up to constant scalar multiples.}
\vskip 0.1in

\begin{figure}
    \centering
    \includegraphics[width=0.5\linewidth]{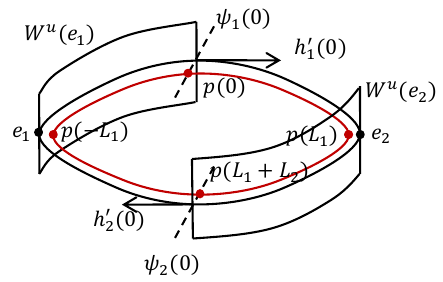}
    \caption{$e_1$ and $e_2$ are equilibrium points of system \eqref{201}. The black solid line from $e_1$ to $e_2$ and the black solid line from $e_2$ to $e_1$ represent the heteroclinic orbits $h_1(x)$ and $h_2(x)$, respectively. The red solid line represents the periodic orbit $p(x)$ bifurcating from the heteroclinic loop. In the schematic diagram, we divide the periodic orbit into four segments at $x=-L_1$( or equivalently $L_1+2L_2$), $0$, $L_1$, and $L_1+L_2$.}
    \label{fig}
\end{figure}

On the one hand, this last hypothesis implies that the adjoint variational equation
\begin{equation}\label{203}
w'= -F_u(h_i(x),  {\bf 0})^* w, \quad x \in \mathbb{R},
\end{equation}
also has a unique nontrivial bounded solution, which we denote by $\psi_i(x)$; see Figure \ref{fig}. Moreover, $\psi_i(x)$ satisfies
\begin{equation}\label{240}
  |\psi_i(x)| \leq Ce^{-\alpha |x|}, \quad x \in \mathbb{R},\quad i=1,2,  
\end{equation}
where $\al$ is defined as
$$\al :=\min\{\al_1^s,\al_1^u,\al_2^s,\al_2^u\}.$$

On the other hand, hypothesis (H3) implies that \(\lambda = 0\) has geometric multiplicity one as an eigenvalue of the PDE linearization of the underlying traveling front (or back). In other words, system \eqref{202} admits exponential dichotomies on both $\mathbb{R}_{\pm}$ with projections $P_{i,\pm}^{u,s}(x)$. Denote the state transition matrix of \eqref{202} by $\Phi_{i}(x,y)$. Then the range of the projection $P_{i,+}^{s}(0)$ is determined by
$$\mathrm{R}(P_{i,+}^{s}(0))=\{v_0\in\mathbb{R}^n:\sup_{x\geq0}|\Phi_{i}(x,0)v_0|<\infty\},$$
while the kernel may be any complementary subspace. Similarly, the kernel of the projection $P_{i,-}^{s}(0)$ is determined by
$${\rm Ker}(P_{i,-}^{s}(0))=\{v_0\in\mathbb{R}^n:\sup_{x\leq0}|\Phi_{i}(x,0)v_0|<\infty\},$$
while the range may be any complementary subspace, see \cite{Cop}. For convenience, we denote $\Phi_{i,\pm}^{u,s}(x,y):=\Phi_{i}(x,y)P_{i,\pm}^{u,s}(y)$, which is called unstable or stable evolution of \eqref{202}.
Thus, we have the following lemma about exponential dichotomies.

\begin{lemma} \label{lem2.1}
 The operators \(\Phi^s_{i,+}( x, y)\) and \(\Phi^u_{i,+}(x, y)\) satisfy
\[|\Phi^s_{i,+}(x, y)| \leq Ce^{-\alpha_j^s|x-y|},~~~|\Phi^u_{i,+}(y, x)| \leq Ce^{-\alpha_j^u|x-y|}, \text{ for } x \geq y \geq 0,\]
and the projection \(P^s_{i,+}(x)\) satisfies
\begin{equation}\label{204}
|P^s_{i,+}(x) - P^s_j| \leq Ce^{-\alpha_j^s x}, \text{ for } x \geq 0,
\end{equation}
where $P^s_j$ is the stable spectral projection of $F_u(e_j,0)$.  Analogously, the operators \(\Phi^s_{i,-}(y, x)\) and \(\Phi^u_{i,-}(y, x)\) satisfy
\[|\Phi^s_{i,-}(y, x)| \leq Ce^{-\alpha_i^s|y-x|},~~~|\Phi^u_{i,-}(x, y)| \leq Ce^{-\alpha_i^u|y-x|}, \text{ for } x \leq y \leq 0,\]
and the projection \(P^u_{i,-}(x)\) satisfies
\begin{equation}\label{205}
|P^u_{i,-}(x) - P^u_i| \leq Ce^{-\alpha_i^u |x|},\text{ for }x \leq 0,
\end{equation}
where $P^u_j$ is the unstable spectral projection of $F_u(e_j,0)$.
In addition, there exist spaces \(Y_i^u\)  and \(Y_i^s\) such that \(Y_i^c \oplus Y_i^s \oplus Y_i^u \oplus Y_i^\bot = \mathbb{R}^n\) and
\begin{equation} \label{206}
\begin{aligned}
\mathrm{R}(P^s_{i,+}(0)) &= Y_i^c \oplus Y_i^s, & \mathrm{R}(P^u_{i,+}(0)) &= Y_i^u \oplus Y_i^\bot, \\
\mathrm{R}(P^u_{i,-}(0)) &= Y_i^c \oplus Y_i^u, & \mathrm{R}(P^s_{i,-}(0)) &= Y_i^s \oplus Y_i^\bot,\quad i,j=1,2, i\neq j,
\end{aligned}
\end{equation}
where the subspaces $Y_i^c$ and $Y_i^\bot$ are defined as
$$Y_i^c := {\rm{span}}\{h_i'(0)\},\quad Y_i^\bot :=  {\rm{span}}\{\psi_i(0)\}.$$
\end{lemma}

\smallskip
We remark that the adjoint system also has exponential dichotomies on $\mathbb{R}_{\pm}$ with the stable and unstable evolution denoted by $\widetilde{\Phi}_{i,\pm}^{s,u}(x, y)$, see \cite{sand02}. Moreover, we have 
\begin{equation*}
    \begin{aligned}
(\Phi_{1,+}^u(x, y))^* = \widetilde{\Phi}_{1,+}^s(y, x),
\quad (\Phi_{1,-}^s(x, y))^* = \widetilde{\Phi}_{1,-}^u(y, x),
    \end{aligned}
\end{equation*}
which we used in the proof of Lemma \ref{lem2.4}.

\textbf{Hypothesis (H4).} The Melnikov integrals
\[N_i := \int_{-\infty}^{\infty} \langle \psi_i(x), D_{\mu} F(h_i(x), {\bf 0}) \rangle dx \in \mathbb{R}^2, \quad i = 1, 2,\]
are linearly independent (and in particular nontrivial).

Then the existence of periodic solutions is stated as follows.
\begin{theorem}\label{thm3.1}
Assume that the hypotheses (H1)-(H4) are satisfied. There exist positive constants \( C \), \( \delta \), and \( L_* \), such that the following result holds.\\
 (1) If \( L_1 + L_2 > L_* \), then equation (\ref{201}) has a periodic solution \( p(x) \) with period \(T:= 2(L_1 + L_2) \) at \( \mu = \mu_T \), which satisfies 
\[\sup_{|x| \leq L_1} |p(x) - h_1(x)| + \sup_{|x-L_1-L_2| \leq L_2} |p(x) - h_2(x-L_1-L_2)| < \delta, \quad |\mu_T| < \delta,\]
if and only if
\begin{equation} \label{301}
    \begin{aligned}
 \langle \psi_1(L_1), h_2(-L_2) \rangle - \langle \psi_1(-L_1), h_2(L_2) \rangle - \mu \int_{-\infty}^{\infty} \langle \psi_1(x), F_{\mu}(h_1(x), {\bf 0}) \rangle dx  + R_1(\mu) =  0, \\
\langle \psi_2(L_2), h_1(-L_1) \rangle - \langle \psi_2(-L_2), h_1(L_1) \rangle - \mu \int_{-\infty}^{\infty} \langle \psi_2(x), F_{\mu}(h_2(x), {\bf 0}) \rangle dx + R_2(\mu) =  0, \\
    \end{aligned}
\end{equation}
at \(\mu = \mu_T\) , where \(R(\mu)\) is differentiable in \(\mu\) and
\begin{equation} \label{302}
    \begin{aligned}
 R(\mu) \leq C(e^{-\alpha L} + |\mu|)(|\mu| + e^{-2\alpha L}), \quad \partial_{\mu} R(\mu) \leq C(e^{-\alpha L} + |\mu|),
    \end{aligned}
\end{equation}
with 
$$L=\min\{L_1,L_2\}.$$
(2) Furthermore, for any such periodic solution \(p(x)\), the following estimates hold:
\begin{equation}\label{303}\begin{aligned} 
\sup_{|x| \leq L_1} |p(x) - h_1(x)| + \sup_{|x-L_1-L_2| \leq L_2} |p(x) - h_2(x-L_1-L_2)| \leq & C(|\mu_T| + e^{-\alpha L}), \\
|P_1^s (p(L_1 + 2L_2)-h_2(L_2))| + |P_2^u( p(L_1) - h_2(-L_2))|\leq & C(|\mu_T| + e^{-2\alpha L}), \\
|P_1^u( p(L_1 + 2 L_2)-h_1(-L_1))| + |P_2^s (p(L_1 ) - h_1(L_1))|\leq & C(|\mu_T| + e^{-2\alpha L}), \\
 |\mu_T| \leq & C e^{-2\alpha L}.
\end{aligned}
\end{equation}
\end{theorem}
\begin{remark}
Although Deng has proved the existence of periodic solutions of \eqref{201}, here we give another proof and more detailed estimates about these periodic solutions, which are prepared for locating the spectrum. 
    Equality \eqref{301}, estimates \eqref{302} and the first inequality in \eqref{303} can be found in \cite{Lin}. Then solving \eqref{301}, we conclude that $|{\mu}_T|\leq Ce^{-\al L}$ by using \eqref{302}. These inequalities in \eqref{303} will be used to locate the spectrum of related periodic waves. 
\end{remark}

\begin{proof}
Introduce the decomposition
\begin{equation*}
    \begin{aligned}
   p(x) =& h_1(x) + w_1(x), \quad x \in (-L_1, L_1),\\
   p(x) =& h_2(x-L_1-L_2) + w_2(x), \quad x \in (L_1, L_1+2L_2).
    \end{aligned}
\end{equation*}
Then $p(x)$ is a periodic solution of \eqref{201} with period \(2(L_1+L_2)\) if and only if \(w\) satisfies
\begin{equation} \label{304}
\begin{aligned}
(\romannumeral 1) \quad & w_{1,\pm}'(x) = F(h_1(x)+ w_{1,\pm}, \mu)- F(h_1(x),  {\bf 0}),  \quad x \in (-L_1, L_1),\\
(\romannumeral 2) \quad & w_{2,\pm}'(x) = F(h_2(x-L_1-L_2) + w_{2,\pm}, \mu) - F(h_2(x-L_1-L_2) ,  {\bf 0}),  \quad  x \in (L_1, L_1+2L_2),\\
(\romannumeral 3) \quad & P(Y_1^c, Y_1^s \oplus Y_1^u \oplus Y_1^\bot)w_{1,\pm}(0) = {\bf 0}, \\
(\romannumeral 4) \quad & P(Y_2^c, Y_2^s \oplus Y_2^u \oplus Y_2^\bot)w_{2,\pm}(L_1+L_2) = {\bf 0}, \\
(\romannumeral 5) \quad & w_{1,+}(0) - w_{1,-}(0) \in Y_1^\bot, \\
(\romannumeral 6) \quad & w_{2,+}(L_1+L_2) - w_{2,-}(L_1+L_2) \in Y_2^\bot, \\
(\romannumeral 7) \quad & w_{1,+}(L_1) - w_{2,-}(L_1)= h_2(-L_2) - h_1(L_1), \\
(\romannumeral 8) \quad & w_{2,+}(L_1+2L_2) - w_{1,-}(-L_1) = h_1(-L_1) - h_2(L_2),
\end{aligned}
\end{equation}
and
\begin{equation}\label{305}
\begin{aligned}
&\langle \psi_1(0), w_{1,+}(0) - w_{1,-}(0) \rangle = 0,  \\
&\langle \psi_2(0), w_{2,+}(L_1+L_2) - w_{2,-}(L_1+L_2) \rangle = 0,
\end{aligned}
\end{equation}
where $P(X,Y)$ represents a projection on $X$ along $Y$, and the subscripts $\{1,-\},\{1,+\},\{2,-\}$ and $\{2,+\}$ always correspond to $x\in(-L_1,0),(0,L_1),(L_1,L_1+L_2)$ and $(L_1+L_2,L_1+2L_2)$.
Furthermore, Lin (Theorem 3.1 in \cite{Lin}) proved that for every $L_1+L_2> L_*$ and $|\mu|<\delta$, there exists a unique piecewise continuous solution of \eqref{304}, which is differentiable in $\mu$ for fixed $x$ and satisfies
$$|w|\leq C(e^{-\al L}+|\mu|),\quad |\partial_\mu^l w|\leq C,\quad \text{ for }l\geq1.$$

Note that
\begin{equation}\label{F}
\begin{aligned}
&F(h(x) + w, \mu) - F(h(x),  {\bf 0})  \\
 = & F_u(h(x), {\bf 0}) w + F_\mu(h(x),  {\bf 0}) \mu + \int_0^1 [F_{u}(h(x) + \tau w,  {\bf 0}) - F_{u}(h(x),  {\bf 0})] d\tau w \\
&+ \int_0^1 [F_{\mu}(h(x) + w, \tau\mu) - F_{\mu}(h(x), {\bf 0})] d\tau \mu.
\end{aligned}
\end{equation}
Thus, this theorem can be proved by applying a similar method used in the next subsection and Appendix, so we omit the proof.
\end{proof}

\subsection{Spectral analysis of periodic waves with large periods}
In this subsection, we will analyze the spectrum of periodic waves with a large period. We focus on the spectrum with small modulus. Equivalently, for \( \xi \in [-\pi/T,\pi/T) \), we locate $\lambda\in\mathbb{C}$ such that the following boundary value problem admits nontrivial bounded solution:
\begin{equation} \label{207}
v' = (F_u(p(x), \mu_T) + \lambda B) v,
\end{equation}
 \begin{equation} \label{208}
v(L_1+2L_2) = e^{i\xi T}v(-L_1).
 \end{equation}
We state the main result in this subsection.
\begin{theorem}\label{thm2.1}
Assuming that hypotheses (H1)-(H4) are satisfied, there exist positive constants
\(C\) and \(\delta\) such that the following property holds.
The boundary-value problem (\ref{207})-(\ref{208}) has a solution \((\lambda, \xi, v)\) for \(\lambda \in \mathbb{C}\) with \(|\lambda| < \delta\) and \(\xi \in [-\pi/T,\pi/T)\) if and only if
\begin{equation*}
\begin{aligned}
0=E(\lambda, \xi) = & (1-e^{i \xi T}) \langle \psi_{2}(L_2), h'_1(-L_1) \rangle\langle \psi_{1}(L_1), h_2'(-L_2) \rangle \\
& + (1-e^{-i \xi T})\langle \psi_{1}(-L_1), h_2'(L_2) \rangle\langle \psi_{2}(-L_2), h'_{1}(L_1) \rangle \\
& - \left( \langle \psi_{2}(-L_2), h'_{1}(L_1) \rangle - \langle \psi_{2}(L_2), h'_1(-L_1) \rangle \right) \lambda \int_{-\infty}^{\infty} \left\langle \psi_{1}(x), Bh'_1(x) \right\rangle dx \\
& - \left(\langle \psi_{1}(-L_1), h_2'(L_2) \rangle - \langle \psi_{1}(L_1), h_2'(-L_2) \rangle \right) \lambda \int_{-\infty}^{\infty} \left\langle \psi_{2}(x), Bh'_2(x) \right\rangle dx  \\
& +(1- e^{i \xi T})R(\lambda, \xi) + \lambda\tilde{R}(\lambda, \xi),
\end{aligned}
\end{equation*}
where \(R(\lambda, \xi)\) and \(\tilde{R}(\lambda, \xi)\) are analytic in \((\lambda, \xi)\) and satisfy, for fixed $j,l\geq0$, the following conditions:
\begin{equation*}
    \begin{aligned}
\left| \partial^{j}_{\lambda}\partial^{\ell}_{\xi} R(\lambda, \xi) \right| \leq &Ce^{-2 \alpha L}\big((|p'(L_1)|+ |p'(L_1+2L_2)|) ( |\lambda| + e^{-\alpha L} + |\mu_{T}| \\
&+ \sup_{|x| \leq L_1} |p(x) - h_1(x)|  + \sup_{|x-L_1-L_2| \leq L_2} |p(x) - h_2(x-L_1-L_2)|)^{2} \\
&+ e^{- \alpha L}(|P_1^s (p(L_1 + 2L_2)-h_2(L_2))| + |P_2^u( p(L_1) - h_2(-L_2))|  \\
&+ |P_1^u( p(L_1 + 2 L_2)-h_1(-L_1))| + |P_2^s (p(L_1 ) - h_1(L_1))| + |\mu_T| + e^{- 2\alpha L} )\big), \\
\left| \partial^{\ell}_{\xi} \tilde{R}(\lambda, \xi) \right| \leq & C e^{-2\alpha L}\big(e^{-\alpha L} + |\lambda| + |\mu_{T}| + \sup_{|x| \leq L_1} |p(x) - h_1(x)| \\
&+ \sup_{|x-L_1-L_2| \leq L_2} |p(x) - h_2(x-L_1-L_2)|\big),  \\
\left| \partial^{j+1}_{\lambda} \partial^{\ell}_{\xi} \tilde{R}(\lambda, \xi) \right| \leq & C.
 \end{aligned}
\end{equation*}

\end{theorem}

\begin{proof}
 We  prove this theorem by two steps. First, we give a reformation of the original problem by using the Lyapunov-Schmidt reduction method and solve a general boundary value problem. Second, we substitute the reformation into boundary value problem solved in step 1 and give the estimates about the spectrum.
 
{\bf Step 1: reformation of the original boundary value problem.} Rearrange (\ref{207})-(\ref{208}) as follows:
\begin{equation} \label{209}
\begin{aligned}
&v'_{1, \pm}  = (F_u(p(x), \mu_T) + \lambda B) v_{1, \pm}, \quad x \in (-L_1, L_1), \\
&v'_{2, \pm} = (F_u(p(x), \mu_T) + \lambda B) v_{2, \pm}, \quad x \in (L_1, L_1+2L_2), \\
&v_{1,-}(0) = v_{1,+}(0), \\
&v_{1,+}(L_1) = v_{2,-}(L_1), \\
&v_{2,-}(L_1+L_2) = v_{2,+}(L_1+L_2), \\
&v_{2,+}(L_1+2L_2) = e^{i\xi T}v_{1,-}(-L_1).
\end{aligned}
\end{equation}
Note that the function $p'(x)$ satisfies (\ref{209}) for \(\lambda = \xi = 0\).  Then, we write
\begin{equation} \label{210}
\begin{aligned}
v_{1,\pm}(x) &= p'(x) d_1 + w_{1,\pm}(x), \quad  x \in (-L_1, L_1), \\
v_{2,\pm}(x) &= p'(x) d_2 + w_{2,\pm}(x), \quad x \in (L_1, L_1 + 2L_2), \\
\end{aligned}
\end{equation}
where $d_1, d_2 \in \mathbb{C}$ are constants to be determined. Define projection \(Q_1\) and \(Q_2: \mathbb{R}^n \to \mathbb{R}^n\) by
\[ {\rm{R}}(Q_1) = {{\rm{span}}\{p'(0)\}, \quad {\rm Ker}}(Q_1) = Y_1^s \oplus Y_1^u \oplus Y_1^\bot,\]
and,
\[{\rm{R}}(Q_2) =  {\rm{span}}\{p'(L_1+L_2)\}, \quad {\rm Ker}(Q_2) = Y_2^s \oplus Y_2^u \oplus Y_2^\bot.\]
It should be noted that if the periodic solution $p(x)$ is $\delta$-close to the heteroclinic loop, then ${\rm R}(Q_i)$ is close to the space $Y_{i}^c$, which is a complement of ${\rm Ker}(Q_i)$. Thus, $Q_i ( i = 1,2 )$ is well-defined, and its norm only relies on $\delta$, not on $p(x)$ itself. 

By employing the ansatz (\ref{210}), we have the equivalent form of the boundary-value problem (\ref{209}):
\begin{equation} \label{211}
\begin{aligned}
(\romannumeral 1) \quad & w_{1,\pm}' = F_u(h_1(x),  {\bf 0}) w_{1,\pm} + (F_u(p(x), \mu_T) - F_u(h_1(x),  {\bf 0}) + \lambda B) w_{1,\pm} \\
&\qquad \quad + \lambda Bp'(x) d_1,   \quad x \in (-L_1, L_1),\\
(\romannumeral 2) \quad & w_{2,\pm}' = F_u(h_2(x-L_1-L_2),  {\bf 0}) w_{2,\pm} + (F_u(p(x), \mu_T) - F_u(h_2(x-L_1-L_2),  {\bf 0}) + \lambda B) w_{2,\pm} \\
&\qquad \quad+ \lambda Bp'(x) d_2,   \quad  x \in (L_1, L_1+ 2L_2),\\
(\romannumeral 3) \quad & Q_1w_{1,\pm}(0) =  {\bf 0},\\
(\romannumeral 4) \quad & Q_2w_{2,\pm}(L_1+L_2) =  {\bf 0},\\
(\romannumeral 5) \quad & w_{1,+}(0) - w_{1,-}(0) \in Y_1^\bot,\\
(\romannumeral 6) \quad & w_{2,+}(L_1+L_2)-w_{2,-}(L_1+L_2)\in Y_2^\bot,\\
(\romannumeral 7) \quad & w_{1,+}(L_1)-w_{2,-}(L_1)= p'(L_1) (d_2-d_1), \\
(\romannumeral 8) \quad & w_{2,+}(L_1+2L_2)-e^{i\xi T}w_{1,-}(-L_1) = p'(L_1+2L_2)(e^{i\xi T}d_1-d_2),\\
\end{aligned}
\end{equation}
and
\begin{equation}\label{212}
\begin{aligned}
 &\langle \psi_1(0), w_{1,+}(0) - w_{1,-}(0) \rangle = 0,  \\
 &\langle \psi_2(0), w_{2,+}(L_1+L_2) - w_{2,-}(L_1+L_2) \rangle = 0.
\end{aligned}
\end{equation}
Note that (\ref{211})(\romannumeral 5)-(\romannumeral 6) and (\ref{212}) are met if and only if $w_{1,+}(0) = w_{1,-}(0)$ and $w_{2,+}(L_1+L_2)= w_{2,-}(L_1+L_2)$.

Thus, solving the boundary value problem (\ref{207})-(\ref{208}) is reduced to solve (\ref{211})-(\ref{212}). Note that \eqref{211} and the boundary value problem of existence \eqref{304} are similar in form by substituting \eqref{F} into \eqref{304}. Thus, we study the following boundary value problem.
\begin{equation}\label{243}
\begin{aligned}
(\romannumeral 1) \quad & w_{1,\pm}' = F_u(h_1(x),  {\bf 0}) w_{1,\pm} + H_1(x) w_{1,\pm} + g_{1,\pm}(x),  \quad x \in (-L_1, L_1)\\
(\romannumeral 2) \quad & w_{2,\pm}' = F_u(h_2(x-L_1-L_2),  {\bf 0}) w_{2,\pm} + H_2(x) w_{2,\pm} + g_{2,\pm}(x), \quad  x \in (L_1, L_1+2L_2)\\
(\romannumeral 3) \quad & Q_1w_{1,\pm}(0) =  {\bf 0}, \\
(\romannumeral 4) \quad & Q_2w_{2,\pm}(L_1+L_2) =  {\bf 0}, \\
(\romannumeral 5) \quad & w_{1,+}(0) - w_{1,-}(0) \in Y_1^\bot, \\
(\romannumeral 6) \quad & w_{2,+}(L_1+L_2) - w_{2,-}(L_1+L_2) \in Y_2^\bot, \\
(\romannumeral 7) \quad & w_{1,+}(L_1) - w_{2,-}(L_1)= D_1, \\
(\romannumeral 8) \quad & w_{2,+}(L_1+2L_2) - e^{i\xi T}w_{1,-}(-L_1) = D_2, \\
\end{aligned}
\end{equation}
for  $g = (g_{1,-}, g_{1,+},g_{2,-},g_{2,+}) \in V_{w}:= C^0([-L_1, 0], \mathbb{C}^n) \times C^0([0, L_1], \mathbb{C}^n) \times C^0([L_1, L_1+L_2], \mathbb{C}^n) \times C^0([L_1+L_2, L_1+2L_2], \mathbb{C}^n) $,   $D \in \mathbb{C}^{2n}$  and  $H = (H_1, H_2) \in V_{H} := C^{0} ( [-L_1, L_1], \mathbf{L}(\mathbb{C}^n) ) \times C^{0} ( [L_1, L_1+2L_2], \mathbf{L}(\mathbb{C}^n) )$. Let $U_{\delta}$ denote the $\delta$-neighborhood of $H = 0$ in $V_{H}$.

Next, we shall seek  solutions to (\ref{243}) in $V_{w}$ for $\xi \in [-\pi/T,\pi/T)$.

\begin{lemma}\label{lem2.5}
For fixed $l\geq0$, there exist positive constants $C$, $L_{*}$, and $\delta$, such that for $L_1+L_2 > L_{*}$, the following statement holds. A solution operator exists, denoted by 
$$W: [-\pi/T,\pi/T) \times U_{\delta} \rightarrow \mathbf{L}( \mathbb{C}^{2n} \times V_{w}, V_{w} ), \quad( \xi, H ) \mapsto W( \xi, H ),$$
which depends analytically on the parameter $\xi$ and the perturbation $H$, such that $w$ satisfies the boundary-value problem (\ref{243}) if and only if $w = W( \xi, H )(D, g)$.
Moreover, the operator satisfies the following estimate:
\begin{equation}\label{244}
| \partial^\ell_{(\xi, H)} W | \leq C.
\end{equation}
 Furthermore, we have
\begin{equation}\label{245}
\begin{aligned}
&(w_{1,-}(-L_1), w_{1,+}(L_1), w_{2,-}(L_1), w_{2,+}(L_1 + 2L_2)) \\
= &(-e^{-i\xi T} P_1^s D_2, P_2^u D_1, -P_2^s D_1, P_1^u D_2) + \widetilde{\mathcal{W}}_2(\xi)(D, g) + \widetilde{\mathcal{W}}_2(\xi)(0, H W(\xi, H)(D, g)),
\end{aligned}
\end{equation}
where the definition of \(\widetilde{\mathcal{W}}_2\) can be found in Proposition \ref{prop2.1}. 

The jump condition for \(w\) at \(x = 0\) is given as follow:
\begin{equation}\label{246}
\begin{aligned}
0=\Xi_1 &= \langle \psi_1(L_1), P_2^u D_1 \rangle + e^{-i \xi T} \langle \psi_1(-L_1), P_1^s D_2 \rangle \\
&\quad - \int_{-L_1}^{L_1} \langle \psi_1(x), g_1(x) \rangle \, dx + \hat{R}_{2,1}(\xi, H)(D, g),
\end{aligned}
\end{equation}
and at \(x = L_1 + L_2\), the jump conditions are given by:
\begin{equation}\label{247}
\begin{aligned}
0=\Xi_2 &= \langle \psi_2(L_2), P_1^u D_2 \rangle + \langle \psi_2(-L_2), P_2^s D_1 \rangle \\
&\quad - \int_{L_1}^{L_1 + 2L_2} \langle \psi_2(x - L_1 - L_2), g_2(x) \rangle \, dx + \hat{R}_{2,2}(\xi, H)(D, g),
\end{aligned}
\end{equation}
where \(g_1 = (g_{1,-}, g_{1,+})\) and \(g_2 = (g_{2,-}, g_{2,+})\).

Let \(\hat{R}_2 = (\hat{R}_{2,1}, \hat{R}_{2,2})\), where \(\hat{R}_2: [-\pi/T, \pi/T) \times U_{\delta} \to \mathbf{L}(\mathbb{C}^{2n} \times V_w, \mathbb{C}^2)\) is analytic with respect to \((\xi, H)\), and its Taylor expansion in \(H\) is given by:
$$\hat{R}_2(\xi, H)(D, g) = T_0(\xi)(D, g) + T_1(\xi)(D, g)[H] + T_2(\xi, H)(D, g)[H, H],$$
where the following bounds hold:
\begin{equation}\label{248}
\begin{aligned}
| \partial^\ell_\xi T_0(\xi)(D, g) | &\leq C e^{-\alpha L}(e^{-\alpha L} |D| + |g|), \\
| \partial^\ell_\xi T_1(\xi)(D, g) | &\leq C(e^{-\alpha L} |D| + |g|), \\
| \partial^\ell_{(\xi, H)} T_2(\xi, H) | &\leq C.
\end{aligned}
\end{equation} 
 \end{lemma}
 
\begin{proof}
The proof of this lemma can be found in Appendix \ref{A}.
\end{proof}

{\bf Step 2: substituting and estimating.} Returning to the boundary-value problem (\ref{211}), a comparison between the original boundary-value problem (\ref{213}) and the boundary-value problem (\ref{243}) studied in the previous section shows that they become identical when we set
\begin{equation}\label{256}
    \begin{aligned}
D_1 &= p'(L_1)(d_2-d_1), \\
D_2 &= p'(L_1+2L_2)(e^{i \xi T}d_1-d_2), \\
g_1(x) &= \lambda Bp'(x)d_1, \quad x \in \left(-L_1, L_1 \right),\\
g_2(x) &= \lambda Bp'(x)d_2, \quad x \in \left(L_1, L_1+2L_2 \right),\\
H_1(x) &= F_{u}(p(x), \mu_{T}) - F_{u}(h_1(x),  {\bf 0}) + \lambda B, \quad x \in \left(-L_1, L_1 \right),\\
H_2(x) &= F_{u}(p(x), \mu_{T}) - F_{u}(h_2(x-L_1-L_2),  {\bf 0}) + \lambda B, \quad x \in \left(L_1, L_1+2L_2 \right).
  \end{aligned}
\end{equation}
Then, according to Theorem \ref{thm3.1}, we conclude
\begin{equation}\label{257}
    \begin{aligned}
|H| \leq C(|\lambda| + |\mu_{T}| + \sup_{|x| \leq L_1} |p(x) - h_1(x)|  + \sup_{|x-L_1-L_2| \leq L_2} |p(x) - h_2(x-L_1-L_2)|) \leq C \delta.
  \end{aligned}
\end{equation}
Therefore, Lemma \ref{lem2.5} can be applied. We obtain a solution \( w \) that is analytic in \( (\lambda, \xi) \), with the jumps given by
\begin{equation}\label{258}
    \begin{aligned}
\Xi_1 = &\langle \psi_{1}(L_1), P^u_{2}p'(L_1) \rangle (d_{2}-d_1)
+ \langle \psi_{1}(-L_1), P^s_{1}p'(L_1+2L_2) \rangle (d_{1}-e^{-i \xi T}d_2) \\
&- \lambda  \int_{-L_1}^{L_1} \langle \psi_{1}(x), Bp'(x) \rangle d_{1} dx + \hat{R}_{2, 1}(\xi, H)(D, g),
\end{aligned}
\end{equation}
and,
\begin{equation}\label{259}
    \begin{aligned}
\Xi_2 = &\langle \psi_{2}(L_2), P^u_{1}p'(L_1+2L_2) \rangle (e^{i\xi T}d_1-d_2)
+ \langle \psi_{2}(-L_2), P^s_{2}p'(L_1) \rangle (d_2-d_1) \\
&- \lambda  \int_{L_1}^{L_1+2L_2} \langle \psi_{2}(x-L_1-L_2), Bp'(x) \rangle d_{2} dx + \hat{R}_{2, 2}(\xi, H)(D, g).
\end{aligned}
\end{equation}

First, due to \(\hat{R}_{2,i}(\xi,H)\) $(i=1,2)$ is linear in \((D, g)\), we obtain
\begin{align*}
\hat{R}_{2,i}(\xi, H)(D, g) =& \hat{R}_{2,i}(\xi, H)(D, 0) + \hat{R}_{2,i}(\xi, H)(0, g),
\end{align*}
where \(D\), \(g\) and \(H\) are replaced by the expressions in (\ref{256}). Moreover,
\begin{align*}
 \hat{R}_{2,i}(\xi, H)(D, 0)&= \acute{R}^{(1)}_{2,i}(\lambda, \xi)(d_2 - d_1) + \acute{R}^{(2)}_{2,i}(\lambda, \xi) (e^{i \xi T}d_1 - d_2), \\
 \hat{R}_{2,i}(\xi, H)(0, g) &= \lambda\Grave{R}^{(1)}_{2,i}(\lambda, \xi) d_1 + \lambda\Grave{R}^{(2)}_{2,i}(\lambda, \xi) d_2.
\end{align*}

Substituting the estimates (\ref{257}) for \(H\) into the estimates (\ref{248}) for \(\hat{R}_{2,i}\) in Lemma \ref{lem2.5}, it is obvious that
\begin{equation} \label{260}
    \begin{aligned}
|\acute{R}^{(j)}_{2,i}(\lambda, \xi)| \leq & C\left|\hat{R}_{2,i}(\xi, H)(D, 0)\right| \\
\leq & C|D|\left(e^{-\alpha L} + |H|\right)^2 \\
\leq & C (|p'(L_1)|+|p'(L_1+2L_2)|) (|\lambda| + |\mu_{T}| + \sup_{|x| \leq L_1} |p(x) - h_1(x)|  \\
&+ \sup_{|x-L_1-L_2| \leq L_2} |p(x) - h_2(x-L_1-L_2)| ) ^2,\\
|\Grave{R}^{(j)}_{2,i}(\lambda, \xi)| \leq &  C\left|\hat{R}_{2,i}(\xi, H)(0,g/\la)\right| \\
\leq & C\left(e^{-\alpha L} + |H|\right) \\
\leq & C (e^{-\alpha L} + |\lambda| + |\mu_{T}| + \sup_{|x| \leq L_1} |p(x) - h_1(x)|  \\
&+ \sup_{|x-L_1-L_2| \leq L_2} |p(x) - h_2(x-L_1-L_2)|),
   \end{aligned}
\end{equation}
for $i,j = 1,2 $. Similar estimates can be derived for the derivatives with respect to \((\lambda, \xi)\).

In order to find $(d_1, d_2)\neq0$ to solve $\Xi_1=0$ and $\Xi_2=0$, we consider the following algebraic equation :
\begin{equation}\label{261}
\begin{pmatrix}
\Xi_1 \\
\Xi_2
\end{pmatrix}
=
\begin{pmatrix}
A_{11} & A_{12}  \\
A_{21} & A_{22}
\end{pmatrix}
\begin{pmatrix}
d_1  \\
d_2
\end{pmatrix}
=
\begin{pmatrix}
0  \\
0
\end{pmatrix},
\end{equation}
where
\begin{equation}\label{262}
\begin{aligned}
A_{11} =& \langle \psi_{1}(-L_1), P^s_1 p'(L_1+2L_2) \rangle - \langle \psi_{1}(L_1), P^u_2 p'(L_1) \rangle - \lambda \int_{-L_1}^{L_1} \langle \psi_{1}(x), Bp'(x) \rangle 
dx \\ 
&- \acute{R}^{(1)}_{2,1}(\lambda, \xi) + e^{i \xi T}\acute{R}^{(2)}_{2,1}(\lambda, \xi) +  \lambda \Grave{R}^{(1)}_{2,1}(\lambda, \xi),\\
A_{12} =& \langle \psi_{1}(L_1), P^u_2 p'(L_1) \rangle - e^{-i \xi T}\langle \psi_{1}(-L_1), P^s_1 p'(L_1+2L_2) \rangle 
+ \acute{R}^{(1)}_{2,1}(\lambda, \xi) -  \acute{R}^{(2)}_{2,1}(\lambda, \xi) \\
&+ \lambda \Grave{R}^{(2)}_{2,1}(\lambda, \xi),\\
A_{21} =& e^{i \xi T}\langle \psi_{2}(L_2), P^u_1 p'(L_1+2L_2) \rangle - \langle \psi_{2}(-L_2), P^s_2p'(L_1) \rangle -\acute{R}^{(1)}_{2,2}(\lambda, \xi)  + e^{i \xi T} \acute{R}^{(2)}_{2,2}(\lambda, \xi) \\
&+ \lambda \Grave{R}^{(1)}_{2,2}(\lambda, \xi),\\
A_{22} =& \langle \psi_{2}(-L_2), P^s_2p'(L_1) \rangle - \langle \psi_{2}(L_2), P^u_1 p'(L_1+2L_2) \rangle - \lambda \int_{L_1}^{L_1+2L_2} \langle \psi_{2}(x-L_1-L_2), Bp'(x) \rangle dx \\
&+ \acute{R}^{(2)}_{2,2}(\lambda, \xi) - \acute{R}^{(1)}_{2,2}(\lambda, \xi) + \lambda \Grave{R}^{(2)}_{2,2}(\lambda, \xi).
\end{aligned}
\end{equation}

Next, we replace the functions \( p'(x) \) that appear in the scalar products by the functions $ h'_{1}(x)$  or \(h'_{2}(x) \), and change the interval of integration from \([-L_1, L_1] \) or \( [L_1, L_1+2L_2] \) to $\mathbb{R}$.

Exploiting the Taylor expansion of \( F \), and using Lemma \ref{lem2.1} and (\ref{240}), we get
\begin{equation*} 
    \begin{aligned}
&\left| \left\langle \psi_{1}(-L_1), P^s_{1}p'(L_1+2L_2)- h_2'(L_2) \right\rangle \right| \\
= &\left| \left\langle \psi_{1}(-L_1), P^s_{1}p'(L_1+2L_2)- P^s_{2,+}(L_2)h_2'(L_2) \right\rangle \right| \\
\leq & Ce^{- \alpha L}\left(\left| P^s_{1}p'(L_1+2L_2)- P^s_{1}h_2'(L_2) \right| + \left|  P^s_{1}- P^s_{2,+}(L_2) \right| \left| h_2'(L_2)  \right|\right) \\
\leq & Ce^{- \alpha L}\left(\left| P^s_{1}(F(p(L_1+2L_2),\mu)-F(h_2(L_2), {\bf 0})) \right| + \left|  P^s_{1}- P^s_{2,+}(L_2) \right| \left| h_2'(L_2)  \right|\right)\\
\leq & Ce^{- \alpha L}\left(\left| P^s_{1}(F_u(h_2(L_2), {\bf 0})((p(L_1+2L_2)-h_2(L_2))\right|+\left|F_\mu(h_2(L_2), {\bf 0})\mu\right|\right.\\
&+\left.O((\left|\mu\right|+\left|p(L_1+2L_2)-h_2(L_2)\right|)^2 )+ \left|  P^s_{1}- P^s_{2,+}(L_2) \right| \left| h_2'(L_2)  \right|\right)\\
\leq & Ce^{- \alpha L}\left(\left| P^s_{1}(F_u(h_2(L_2), {\bf 0})-F_u(e_1, {\bf 0}))((p(L_1+2L_2)-h_2(L_2)))+  P^s_{1}F_u(e_1, {\bf 0})(p(L_1+2L_2)-h_2(L_2))\right|\right) \\
&+ O(e^{-3\al L})\\
\leq & Ce^{- \alpha L}\left( \left|P_1^s(p(L_1+2L_2)-h_2(L_2)) \right| + \left|\mu_{T} \right|  + e^{- 2\alpha L} \right)\\
\leq& Ce^{-3\al L}, 
\end{aligned}
\end{equation*} where we use $P_1^sF_u(e_1, {\bf 0})=F_u(e_1, {\bf 0})P_1^s$.
Similarly, we have 
\begin{equation*}\begin{split}
 &\left| \left\langle \psi_{1}(L_1), P^u_{2}p'(L_1)- h_2'(-L_2) \right\rangle \right| \\
\leq & Ce^{- \alpha L}\left( \left|P_2^u( p(L_1)-h_2(-L_2) )\right| + \left|\mu_{T} \right|  + e^{- 2\alpha L} \right), \\
&\left| \left\langle \psi_{2}(-L_2), P^s_{2}p'(L_1)- h'_1(L_1) \right\rangle \right| \\
\leq & Ce^{- \alpha L}\left( \left| P_2^s( p(L_1)-h_1(L_1)) \right| + \left|\mu_{T} \right|  + e^{- 2\alpha L} \right), \\
&\left| \left\langle \psi_{2}(L_2), P^u_{1}p'(L_1+2L_2)- h'_1(-L_1) \right\rangle \right| \\
\leq & Ce^{- \alpha L}\left( \left|P_1^u( p(L_1+2L_2)-h_1(-L_1) ) \right| + \left|\mu_{T} \right|  + e^{- 2\alpha L} \right). \\
  \end{split}
\end{equation*}
Analogously, we see that
\begin{equation} \label{264}
    \begin{aligned}
&\left| \int_{-L_1}^{L_1} \langle \psi_{1}(x), Bp'(x) \rangle dx - \int_{-\infty}^{\infty} \left\langle \psi_{1}(x), Bh'_1(x) \right\rangle dx \right| \\
\leq &C(e^{-2 \alpha L} + \sup_{|x| \leq L_1} |p(x) - h_1(x)| +|\mu_T|), \\
&\left| \int_{L_1}^{L_1+2L_2} \langle \psi_{2}(x-L_1-L_2), Bp'(x) \rangle dx - \int_{-\infty}^{\infty} \left\langle \psi_{2}(x), Bh'_2(x) \right\rangle dx \right| \\
\leq &C(e^{-2 \alpha L} + \sup_{|x-L_1-L_2| \leq L_2} |p(x) - h_2(x-L_1-L_2)| +|\mu_T|).
 \end{aligned}
\end{equation}

Denote
\begin{equation*}
\begin{aligned}
M_1:= \int_{-\infty}^{\infty} \left\langle \psi_{1}(x), Bh'_1(x) \right\rangle dx , \quad M_2:= \int_{-\infty}^{\infty} \left\langle \psi_{2}(x), Bh'_2(x) \right\rangle dx.
\end{aligned}
\end{equation*}
In summary, the jumps (\ref{258}) and (\ref{259}) are given by
\begin{equation}\label{265}
    \begin{aligned}
\Xi_1 = &\langle \psi_{1}(L_1), h_2'(-L_2) \rangle (d_{2}-d_1)
+ \langle \psi_{1}(-L_1), h_2'(L_2) \rangle (d_1-e^{-i \xi T}d_2) -\lambda M_1 d_1 \\
&+ \acute{R}^{(1)}_{2, 1}(\lambda, \xi)(d_2 - d_1) + \acute{R}^{(2)}_{2, 1}(\lambda, \xi) (e^{i \xi T}d_1 - d_2) + \lambda \Grave{R}^{(1)}_{2, 1}(\lambda, \xi)d_1 + \lambda \Grave{R}^{(2)}_{2, 1}(\lambda, \xi)d_2\\
&+ O \big(  e^{- \alpha L} ( \left|P_1^s( p(L_1+2L_2)-h_2(L_2)) \right|+ \left|P_2^u( p(L_1)-h_2(-L_2)) \right| + \left|\mu_{T} \right| + e^{-2\alpha  L}) \\
&+ |\lambda|(e^{-2 \alpha L} +|\mu_T|  + \sup_{|x| \leq L_1} |p(x) - h_1(x)| ) \big),
\end{aligned}
\end{equation}
and
\begin{equation}\label{266}
    \begin{aligned}
\Xi_2 = &\langle \psi_{2}(L_2), h'_1(-L_1) \rangle (e^{i\xi T}d_1-d_2)
+ \langle \psi_{2}(-L_2), h'_{1}(L_1) \rangle (d_2-d_1) - \lambda  M_2 d_2 \\
&+ \acute{R}^{(1)}_{2, 2}(\lambda, \xi)(d_2 - d_1) + \acute{R}^{(2)}_{2, 2}(\lambda, \xi)(e^{i \xi T}d_1 - d_2) + \lambda \Grave{R}^{(1)}_{2, 2}(\lambda, \xi)d_1 + \lambda \Grave{R}^{(2)}_{2, 2}(\lambda, \xi)d_2 \\
&+ O \big(  e^{- \alpha L} ( \left|P_2^s( p(L_1)-h_1(L_1) )\right|+ \left|P_1^u( p(L_1+2L_2)-h_1(-L_1) )\right| + \left|\mu_{T} \right| + e^{-2\alpha  L}) \\
&+ |\lambda|(e^{-2 \alpha L}  +|\mu_T|  + \sup_{|x-L_1-L_2| \leq L_2} |p(x) - h_2(x-L_1-L_2)|) \big).
\end{aligned}
\end{equation}
Based on the above analysis, we rewrite expressions (\ref{262})
\begin{equation*}
\begin{aligned}
A_{11} =& \langle \psi_{1}(-L_1), h'_2(L_2) \rangle - \langle \psi_{1}(L_1),  h'_2(-L_2) \rangle - \lambda M_1 + R^{(1)}_{2,1}(\lambda, \xi) + \lambda \tilde{R}^{(1)}_{2,1}(\lambda, \xi),\\
A_{12} =& \langle \psi_{1}(L_1),  h'_2(-L_2) \rangle - e^{-i \xi T}\langle \psi_{1}(-L_1),  h'_2(L_2) \rangle +  R^{(2)}_{2,1}(\lambda, \xi) + \lambda \tilde{R}^{(2)}_{2,1}(\lambda, \xi),\\
A_{21} =& e^{i \xi T}\langle \psi_{2}(L_2),  h'_1(-L_1) \rangle - \langle \psi_{2}(-L_2), h'_1(L_1) \rangle +  R^{(1)}_{2,2}(\lambda, \xi) + \lambda \tilde{R}^{(1)}_{2,2}(\lambda, \xi),\\
A_{22} =& \langle \psi_{2}(-L_2), h'_1(L_1) \rangle - \langle \psi_{2}(L_2), h'_1(-L_1) \rangle - \lambda M_2 + R^{(2)}_{2,2}(\lambda, \xi) + \lambda \tilde{R}^{(2)}_{2,2}(\lambda, \xi).
\end{aligned}
\end{equation*}
Equation (\ref{261}) has a non-zero solution if and only if,
\begin{equation} \label{267}
E(\lambda, \xi) := A_{11}A_{22}-A_{12}A_{21} = 0.
\end{equation}
Precisely,
\begin{equation} \label{268}
\begin{aligned}
E(\lambda, \xi) = & (1-e^{i \xi T}) \langle \psi_{2}(L_2), h'_1(-L_1) \rangle\langle \psi_{1}(L_1), h_2'(-L_2) \rangle \\
& + (1-e^{-i \xi T})\langle \psi_{1}(-L_1), h_2'(L_2) \rangle\langle \psi_{2}(-L_2), h'_{1}(L_1) \rangle \\
& - \left( \langle \psi_{2}(-L_2), h'_{1}(L_1) \rangle - \langle \psi_{2}(L_2), h'_1(-L_1) \rangle \right) \lambda M_1 \\
& - \left(\langle \psi_{1}(-L_1), h_2'(L_2) \rangle - \langle \psi_{1}(L_1), h_2'(-L_2) \rangle \right) \lambda M_2 \\
& +  (1- e^{i \xi T})R(\lambda, \xi) + \lambda \tilde{R} (\lambda, \xi).
\end{aligned}
\end{equation}
The remainder terms \(R(\lambda, \xi)\) and \(\tilde{R}(\lambda, \xi)\) in the expansion are analytic in \((\lambda, \xi)\).  Combining (\ref{260}), (\ref{265}) and (\ref{266}), then we have
\begin{equation} \label{269}
    \begin{aligned}
\left| \partial^{j}_{\lambda}\partial^{\ell}_{\xi} R(\lambda, \xi) \right| \leq &Ce^{-2 \alpha L}\big((|p'(L_1)|+ |p'(L_1+2L_2)|) ( |\lambda| + e^{-\alpha L} + |\mu_{T}| \\
&+ \sup_{|x| \leq L_1} |p(x) - h_1(x)|  + \sup_{|x-L_1-L_2| \leq L_2} |p(x) - h_2(x-L_1-L_2)|)^{2} \\
&+ e^{- \alpha L}(|P_1^s (p(L_1 + 2L_2)-h_2(L_2))| + |P_2^u( p(L_1) - h_2(-L_2))|  \\
&+ |P_1^u( p(L_1 + 2 L_2)-h_1(-L_1))| + |P_2^s (p(L_1 ) - h_1(L_1))| + |\mu_T| + e^{- 2\alpha L} )\big), \\
\left| \partial^{\ell}_{\xi} \tilde{R}(\lambda, \xi) \right| \leq & C\big( |\lambda| + (|p'(L_1)|+|p'(L_1+2L_2)|) (|\lambda| + |\mu_{T}| + \sup_{|x| \leq L_1} |p(x) - h_1(x)| \\
&+ \sup_{|x-L_1-L_2| \leq L_2} |p(x) - h_2(x-L_1-L_2)| ) ^2 + e^{-2\alpha L}\big(e^{-\alpha L} + |\lambda| + |\mu_{T}| \\
&+ \sup_{|x| \leq L_1} |p(x) - h_1(x)| 
+ \sup_{|x-L_1-L_2| \leq L_2} |p(x) - h_2(x-L_1-L_2)|\big) \big),  \\
\left| \partial^{j+1}_{\lambda} \partial^{\ell}_{\xi} \tilde{R}(\lambda, \xi) \right| \leq & C.
 \end{aligned}
\end{equation}

\end{proof}

\begin{remark} \label{rem2.3}
Using the estimate (\ref{303}) from Theorem \ref{thm3.1} in the following section, we can obtain  the remainder terms \(R(\lambda, \xi)\) and \(\tilde{R}(\lambda, \xi)\) satisfy
\begin{equation} \label{271}
    \begin{aligned}
\left| \partial^{j}_{\lambda}\partial^{\ell}_{\xi} R(\lambda, \xi) \right| \leq & C e^{-3 \alpha L}(|\la|+e^{-\al L})^2,\\
\left| \partial^{\ell}_{\xi} \tilde{R}(\lambda, \xi) \right| \leq & C(|\la|+e^{-3\alpha L}),\\
\left| \partial^{j+1}_{\lambda} \partial^{\ell}_{\xi} \tilde{R}(\lambda, \xi) \right| \leq & C.
 \end{aligned}
\end{equation}
 Therefore, for large enough $T$, the equation $E(\lambda, \xi) = 0 $, approximately reads
\begin{equation*}
\begin{aligned}
0 = & (1-e^{i \xi T}) \langle \psi_{2}(L_2), h'_1(-L_1) \rangle\langle \psi_{1}(L_1), h_2'(-L_2) \rangle \\
& + (1-e^{-i \xi T})\langle \psi_{1}(-L_1), h_2'(L_2) \rangle\langle \psi_{2}(-L_2), h'_{1}(L_1) \rangle \\
& - \left( \langle \psi_{2}(-L_2), h'_{1}(L_1) \rangle - \langle \psi_{2}(L_2), h'_1(-L_1) \rangle \right) \lambda M_1 \\
& - \left(\langle \psi_{1}(-L_1), h_2'(L_2) \rangle - \langle \psi_{1}(L_1), h_2'(-L_2) \rangle \right) \lambda M_2.
\end{aligned}
\end{equation*}
\end{remark}

\subsection{Solving for $\lambda(\xi)$}
In the last subsection, we have obtained the expression of $E(\lambda,\xi)$. In this subsection, we make some further assumptions to solve $E(\lambda,\xi)=0$.
To obtain further results, we need the following assumptions (A1)-(A3). First, we assume that the equilibria $e_1$ and $e_2$ both have simple, real leading eigenvalues; see (A1).

\textbf{(A1).} Assume that the eigenvalue of $F_u(e_i, {\bf 0})$ that is closest to the imaginary axis is real, negative and simple. 

Next, we assume that the solutions have the strong inclination property. In other words, \(h_i(x)\) and \(\psi_i(x)\) converge along the leading directions to the equilibria and zero, respectively; see (A2).

\textbf{(A2).} Assume that the limits
\begin{equation} \label{307}
\begin{aligned}
\lim_{t \to -\infty} e^{-\alpha^u_i x} h'_i(x) &= v_i^-, &
\lim_{t \to \infty} e^{\alpha^s_{i+1} x} h'_i(x) &= v_{i+1}^+, \\
\lim_{t \to -\infty} e^{-\alpha^s_i x} \psi_i(x) &= w_i^+, &
\lim_{t \to \infty} e^{\alpha^u_{i+1} x} \psi_i(x) &= w_{i+1}^-,
\end{aligned}
\end{equation}
are nonzero for \(i = 1, 2\). The index \(i\) is taken modulo two.

Note that \(v_i^{\pm}\) and \(w_i^{\pm}\) are right and left eigenvectors of \( F_u(h_i, 0)\) for the eigenvalues \(\alpha_i^{s,u}\) and $i=1,2$. For situations without strong inclination property, the  analysis in Theorem \ref{thm3.2} below still can be used. 

\textbf{(A3).} Assume that the Melnikov integrals
\[M_i := \int_{-\infty}^{\infty} \langle \psi_i(x), Bh_i^{'}(x) \rangle dx , \quad i = 1, 2,\]
are nonzero.

\begin{theorem} \label{thm3.2}
Assume that the hypothesis (H1)-(H4) and assumptions (A1)-(A3) are satisfied. Fixed $l\geq 0$, there exist positive constants \(C\), \(\delta\), and a function \(\lambda(\xi)\) that is analytic in \(\xi \in [-\pi/T,\pi/T)\),  such that the boundary-value problem (\ref{207})-(\ref{208}) has a solution \((\lambda, \xi , v)\) for \(|\lambda| < \delta\), \(\xi \in [-\pi/T,\pi/T)\) and  \(L_1+L_2 > 1/\delta\) if and only if \(\lambda = \lambda(\xi)\), which has the following expansion:
\begin{equation}\label{308}
\begin{aligned}
\lambda(\xi) =&(1-e^{i\xi T}) \left( \frac{U_1 U_2 - e^{-i\xi T}S_1 S_2}{(S_2 - U_1)M_1 + (S_1 - U_2)M_2 } + R(\xi) \right), \quad \xi \in [-\pi/T,\pi/T),
\end{aligned}
\end{equation}
where
\begin{equation}\label{309}
    \begin{aligned}
S_1 :=  e^{-\alpha_1^s(L_1+L_2)}\langle w_{1}^{+}, v_1^+ \rangle, \qquad
U_1 :=  e^{-\alpha_1^u(L_1+L_2)}\langle w_{1}^{-}, v_1^- \rangle,\\
S_2 :=  e^{-\alpha_2^s(L_1+L_2)}\langle w_{2}^{+}, v_2^+ \rangle, \qquad
U_2 :=  e^{-\alpha_2^u(L_1+L_2)}\langle w_{2}^{-}, v_2^- \rangle.
    \end{aligned}
\end{equation}
The remainder term is analytic in \(\xi\) with
\begin{equation}\label{310}
\begin{aligned}
   |\partial_{\xi}^{\ell} R(\xi)| \leq C e^{-(\alpha + \delta) ( L_1 + L_2) }.
\end{aligned}
\end{equation}
Moreover, the critical eigenvalues near $\lambda = 0$ have a quadratic tangency at zero with the imaginary axis.
\end{theorem}

\begin{proof}
Recall that (\ref{268}) is the reduced equation for solutions of the eigenvalue problem that we derived in from Theorem \ref{thm2.1}. Then we solve for the zeros of
\begin{equation} \label{311}
E(\lambda, \xi)=0,
\end{equation}
where
\begin{equation*}
\begin{aligned}
E(\lambda, \xi) = & (1-e^{i \xi T}) \langle \psi_{2}(L_2), h'_1(-L_1) \rangle\langle \psi_{1}(L_1), h_2'(-L_2) \rangle \\
& + (1-e^{-i \xi T})\langle \psi_{1}(-L_1), h_2'(L_2) \rangle\langle \psi_{2}(-L_2), h'_{1}(L_1) \rangle \\
& - \left( \langle \psi_{2}(-L_2), h'_{1}(L_1) \rangle - \langle \psi_{2}(L_2), h'_1(-L_1) \rangle \right) \lambda M_1 \\
& - \left(\langle \psi_{1}(-L_1), h_2'(L_2) \rangle - \langle \psi_{1}(L_1), h_2'(-L_2) \rangle \right) \lambda M_2  \\
& + (1 - e^{i \xi T} )R(\lambda, \xi) + \lambda\tilde{R}(\lambda, \xi) .
\end{aligned}
\end{equation*}
According to the assumption (A2), we have
\begin{equation} \label{312}
    \begin{aligned}
 &\langle \psi_{1}(-L_1), h_2'(L_2) \rangle = e^{-\alpha_1^s(L_1+L_2)}\langle w_{1}^{+}, v_1^+ \rangle + O(e^{-(\alpha_1^s+ \delta)(L_1+L_2)} ),\\
&\langle \psi_{1}(L_1), h_2'(-L_2) \rangle = e^{-\alpha_2^u(L_1+L_2)}\langle w_{2}^{-}, v_2^- \rangle + O(e^{-(\alpha_2^u+ \delta)(L_1+L_2)} ),\\
 &\langle \psi_{2}(-L_2), h_1'(L_1) \rangle = e^{-\alpha_2^s(L_1+L_2)}\langle w_{2}^{+}, v_2^+ \rangle + O(e^{-(\alpha_2^s+ \delta)(L_1+L_2)} ),\\
&\langle \psi_{2}(L_2), h_1'(-L_1) \rangle = e^{-\alpha_1^u(L_1+L_2)}\langle w_{1}^{-}, v_1^- \rangle + O(e^{-(\alpha_1^u + \delta)(L_1+L_2)} ).
    \end{aligned}
\end{equation}
Then, substituting  (\ref{312}) into (\ref{311}) and using the estimates (\ref{271}) of the remainder terms, we obtain the  equation
\begin{equation} \label{313}
\begin{aligned}
0 = &\lambda \left((S_2 - U_1)M_1 +(S_1 - U_2)M_2 + O( e^{- (\alpha +\delta) (L_1+L_2)}+|\la|) \right) \\
&- (1 - e^{ i\xi T}) \left(U_1 U_2 -  e^{ -i\xi T}S_1 S_2 + O(e^{- (2\alpha + \delta) (L_1 + L_2)}) \right),
 \end{aligned}
\end{equation}
where $U_i$ and $S_i$  are defined in (\ref{309}).
It is easy to verify that $E(\lambda, \xi)$ satisfies
\begin{equation*}
\begin{aligned}
    E(0,0)=0, \qquad \frac{\partial E}{\partial \lambda}|_{ \lambda = \xi =0 } \neq 0.
\end{aligned}
\end{equation*}
Thus, we obtain the expression (\ref{308}) by using the Implicit Function Theorem. The remainder terms $R(\xi)$ are analytic in \((\lambda, \xi)\) and satisfies the estimate (\ref{310}). Moreover, the remainder term $R(\xi)$ is real whenever $(\la,e^{i\xi T})$ is real.

We expand the real part of (\ref{308}) in a Taylor series. The derivatives of this expression, evaluated at \(\xi = 0\) are
 \begin{equation*}
     \begin{aligned}
\lambda'(0) &= -iT \left( \frac{U_1 U_2 - S_1 S_2}{(S_2 - U_1)M_1 + (S_1 - U_2)M_2 } + R(0) \right), \\
\lambda''(0) &= T^2 \left( \frac{U_1 U_2 + S_1 S_2}{(S_2 - U_1)M_1 + (S_1 - U_2)M_2 } + R(0) \right) - 2iT R'(0).
     \end{aligned}
 \end{equation*}

Note that $R(0)=0$. Then we obtain
\begin{equation*}
    \begin{aligned}
\text{Re} \, \lambda(\xi) &= \text{Re} \, \lambda(0) + \xi \, \text{Re} \, \lambda'(0) + \frac{\xi^2}{2} \, \text{Re} \, \lambda''(0) + O(\xi^3) \\
&= \frac{1}{2} T^2 \left( \frac{U_1 U_2 + S_1 S_2}{(S_2 - U_1)M_1 + (S_1 - U_2)M_2 } + R(0) \right) \xi^2 - T \text{Re} (i R'(0)) \xi^2 + O(\xi^3).
    \end{aligned}
\end{equation*}
Using the estimate (\ref{310}), we obtain
\[\text{Re} \, \lambda(\xi) =  \frac{1}{2} T^2 \left( \frac{U_1 U_2 + S_1 S_2}{(S_2 - U_1)M_1 + (S_1 - U_2)M_2 }  + O(e^{-(\alpha + \delta) (L_1+L_2)}) + O(\xi) \right) \xi^2.\]
The proof of the theorem is completed.
\end{proof}

The existence of a heteroclinic loop may arise from certain symmetrical properties of the system. Therefore, we further consider the case where the leading eigenvalues of the system's two equilibrium points are equal.

\begin{cor} \label{cor3.1}
Assume, in addition to the assumptions of Theorem \ref{thm3.2}, that $\alpha_1^s < \alpha_1^u$, $\alpha_2^s<\alpha_2^u$ and $\alpha_1^s = \alpha_2^s$. Fixed $l\geq 0$, there are positive constants \(C\) and \(\delta\) and a function \(\lambda(\xi)\) that is analytic in \(\xi \in [-\pi/T,\pi/T)\) such that the boundary-value problem (\ref{207})-(\ref{208}) has a solution \((\lambda, \xi , v)\) for \(|\lambda| < \delta\), \(\xi \in [-\pi/T,\pi/T)\) and  \(L_1+L_2 > 1/\delta\) if and only if, \(\lambda = \lambda(\xi)\). Then we have the expansion
\begin{equation}\label{314}
\begin{aligned}
\lambda(\xi) = (1 - e^{-i\xi T}) e^{-\alpha_1^s (L_1+L_2)} \left( \frac{\langle w_1^+, v_1^+ \rangle \langle w_2^+, v_2^+ \rangle}{\langle w_2^+, v_2^+ \rangle M_1 + \langle w_1^+, v_1^+ \rangle M_2} + R(\xi) \right), \quad \xi \in [-\pi/T,\pi/T)
\end{aligned}
\end{equation}
where the remainder term $R(\xi)$ satisfies
\begin{equation*}
\begin{aligned}
|\partial_{\xi}^{\ell} R(\xi)| \leq C e^{- \delta( L_1 + L_2) }.
\end{aligned}
\end{equation*}
Moreover, the critical eigenvalues near $\lambda = 0$ has a quadratic tangency at zero with the imaginary axis.
\end{cor}

\begin{cor} \label{cor3.2}
Assume, in addition to the assumptions of Theorem \ref{thm3.2}, that $\alpha_1^s < \alpha_1^u$, $\alpha_2^s < \alpha_2^u$  and $\alpha_1^s < \alpha_2^s$. Fixed $l\geq0$, there are positive constants \(C\) and \(\delta\) and a function \(\lambda(\xi)\) that is analytic in \(\xi \in [-\pi/T,\pi/T)\) such that the boundary-value problem (\ref{207})-(\ref{208}) has a solution \((\lambda, \xi , v)\) for \(|\lambda| < \delta\), \(\xi \in [-\pi/T,\pi/T)\) and  \(L_1+L_2 > 1/\delta\) if and only if, \(\lambda = \lambda(\xi)\). Then we have the expansion
\begin{equation*}
\begin{aligned}
   \lambda(\xi) = ( 1 - e^{-i \xi T}  ) e^{-\alpha_2^s (L_1+L_2)} \left( \frac{\langle w_2^+, v_2^+ \rangle}{M_2} + R(\xi) \right), \quad  \xi \in [-\pi/T,\pi/T)
\end{aligned}
\end{equation*}
where the remainder term $R(\xi)$ satisfies
\begin{equation*}
\begin{aligned}
|\partial_{\xi}^{\ell} R(\xi)| \leq C e^{- \delta( L_1 + L_2) }.
\end{aligned}
\end{equation*}
 Moreover, the critical eigenvalues near $\lambda = 0$ has a quadratic tangency at zero with the imaginary axis.
\end{cor}

\begin{remark}
If the leading eigenvalues of two equilibrium points $e_1$ and $e_2$ are such that one is stable and the other unstable, it is often insufficient to consider only which eigenvalue is smaller. For example, if
$\alpha_1^u < \alpha_1^s $ but $\alpha_2^s < \alpha_2^u $, then it is necessary not only to compare $\alpha_1^u$ and $\alpha_2^s$ but also to compare $\alpha_1^u + \alpha_2^u$ and $\alpha_1^s + \alpha_2^s$. Additional conclusions for other specific cases can be analyzed and derived in a similar manner.
\end{remark}

\begin{remark}
If the leading eigenvalues of the equilibrium point of the heteroclinic loop are a pair of complex conjugate roots, a similar approach can be used to derive the corresponding expression for locating the eigenvalues near zero; see also \cite{SA}. For example, in the reference \cite{GLL}, chaotic traveling wave solutions in coupled Chua's circuits are studied. It is found that at the equilibrium points of the heteroclinic loop, the leading eigenvalues of the system are a pair of complex numbers with negative real parts. In this model, the expected result is that the periodic orbit bifurcating from the heteroclinic loop will continuously change its stability with respect to the period $T$.
\end{remark}

\section{an example: nonlinear stability of periodic waves in FitzHugh-Nagumo equations}\label{stab}
\subsection{Spectral properties}    \setcounter{equation}{0}
Consider the FitzHugh-Nagumo equation \eqref{FHN} in the traveling coordinates
\begin{equation}\label{402}
\begin{aligned}
u_t &= u_{xx} - cu_x + f(u) - w,\\
w_t &= -cw_x + \epsilon (u - \gamma w).
\end{aligned}
\end{equation}
Waves that travel with speed \( c \) are solutions of the ODE
\begin{equation} \label{403}
\begin{aligned}
u' &= v, \\
v' &= c v - f(u)  + w , \\
w' &= \frac{\epsilon}{c} (u - \gamma w) .
\end{aligned}
\end{equation}

\begin{figure}
    \centering
    \subfigure[]{\includegraphics[width=0.4\linewidth]{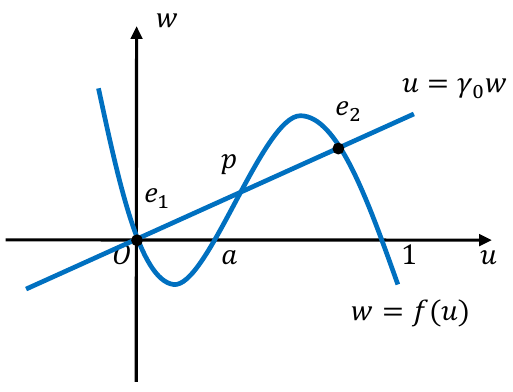}}
\subfigure[]{\includegraphics[width=0.4\linewidth]{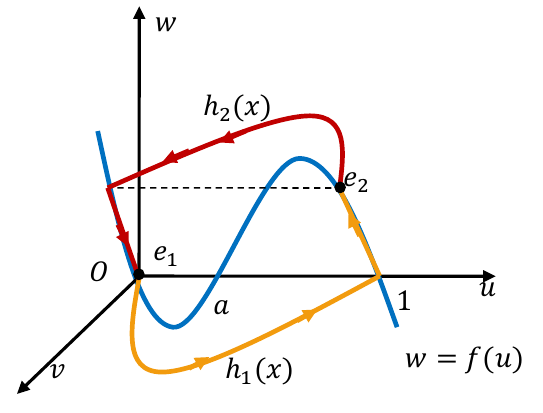}}
    \caption{(a) System (\ref{403}) has three equilibrium points. $\gamma =\gamma_0$ is chosen so that $e_1$ and $e_2$ are symmetric with respect to the inflection point $p(\bar{u},0,\bar{w})$ of the cubic function $w=u(1-u)(u-a)$. (b)  Heteroclinic solutions $h_1(x)$ and $h_2(x)$ exist at the same time for a certain $(\gamma,c) = (\gamma_0, c(\epsilon))$, $h_1(x)$  from $e_1$ to $e_2$, and $h_2(x)$  from $e_2$ to $e_1$ at the same time, forming a heteroclinic loop. }
    \label{fig1}
\end{figure}

 Choose $\gamma = \gamma_0$ such that system (\ref{403}) has three equilibrium points $e_1 (0,0,0)$, $e_2(u_2,0,w_2)$ and $p(\Bar{u},0,\Bar{w})$, where $e_1 (0,0,0)$ and $e_2(u_2,0,w_2)$ are symmetric with respect to point $p(\Bar{u},0,\Bar{w})$; see Figure \ref{fig1}(a). It is well known that (\ref{403}) has  heteroclinic solutions \( h_1(x) \) and $h_2(x)$ for certain parameter values $( \gamma_0,  c(\epsilon))$ for each \( \epsilon > 0 \) sufficiently small:
\begin{equation*}
\lim_{x \to -\infty} h_1(x) = e_1, \quad \lim_{x \to +\infty} h_1(x) = e_2,
\end{equation*}
\begin{equation*}
\left( \lim_{x \to -\infty} h_2(x) = e_2, \quad \lim_{x \to +\infty} h_2(x) = e_1, \text{ respectively} \right).
\end{equation*}
See Figure \ref{fig1}(b) and also see \cite{Den91b, RT}. Yanagida proved in \cite{Yan} that the simple fronts \( h_1(x) \) and backs \( h_2(x) \) building the heteroclinic loop are spectrally stable with respect to the partial differential equation \eqref{402}; that is, the spectrum of the linearized operator is contained in the left half-plane except for a simple eigenvalue at zero. 

\begin{figure}
    \centering

\includegraphics[width=0.5\linewidth]{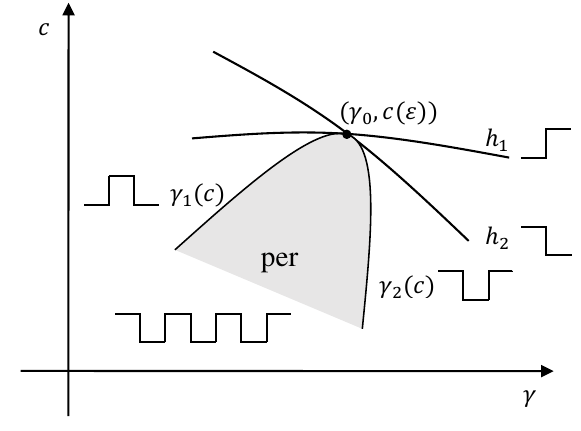}
\caption{A bifurcation diagram of the  equation (\ref{403}) with respect to the 
parameters $(\gamma, c)$. The parameter curves $\gamma_1(c)$ and 
$\gamma_2(c)$ represent the existence of homoclinic orbits to the equilibrium points $e_1$ and $e_2$, respectively. The shaded area between the curves $\gamma_1(c)$ and 
$\gamma_2(c)$ corresponds to the parameter region where periodic orbits exist.}
    \label{bifur}
\end{figure}

The existence of period solutions bifurcating from the heteroclinic loop has been established in suitable \( (\gamma, c) \)-parameter region, see \cite{CDT, Den91a}. They demonstrated that when the parameters $(\gamma, c)$ satisfy
$\gamma_1(c)<\gamma<\gamma_2(c)$, there exist periodic solutions, where $\gamma_1(c)$ and $\gamma_2(c)$ are smooth curves with respect to $c$, see Figure \ref{bifur}.
In the previous section, we also used the Lyapunov-Schmidt reduction method to prove the existence of the large spatial period traveling wave solutions under certain assumptions, along with some corresponding estimates, which applies to \eqref{403}. 

We shall analyze $L^2({\mathbb{R}})$-spectrum of  the periodic  wave denoted by $\mathbf{\bar{U}}(x):=(\phi(x),\psi(x))^t$ of (\ref{402}) by applying Theorem \ref{thm3.2}. Linearized system \eqref{402} about periodic traveling waves $\mathbf{\bar{U}}(x)$, we obtain the following linear differential operator $\mathcal{L}:H^{2}(\mathbb{R})\times H^{1}(\mathbb{R})\subset L^2 (\mathbb{R})\times L^2(\mathbb{R})\to L^2 (\mathbb{R})\times L^2(\mathbb{R}) $ described as
$$\mathcal{L}=\left(\begin{array}{cc}
\partial_{x}^2-c\partial_x+f'(\phi)&-1\\
\e&-c\partial_x-\e\gamma
\end{array}\right),$$
which has periodic coefficients.
\begin{theorem}\label{FHNSpe}
  Fix \( a \in \left(0, \frac{1}{2}\right) \), there exists a constant \( \epsilon_* = \epsilon_*(a) \) such that the following result holds. For every \( \epsilon \) with \( 0 < \epsilon < \epsilon_* \), there exists an \( L_* = L_*(\epsilon) \) such that the heteroclinic loop in the FitzHugh-Nagumo system is accompanied by periodic wave trains with period \( 2(L_1 + L_2) \) for any \( (L_1 + L_2) > L_* \), and all such traveling wave trains are spectrally stable.
\end{theorem}

\begin{proof}
We can divide the right half of the complex plane \(\mathbb{C}\) into three subsets, namely \(\Omega_1\), \(\Omega_2\), and \(\Omega_3\). The definitions are as follows:
\begin{equation*}
\begin{aligned}
\Omega_1 =& \{\lambda \in \mathbb{C} : |\lambda| < r_1\}, \\
\Omega_2 =& \{\lambda \in \mathbb{C} : \operatorname{Re} \lambda \ge -r_1, \, r_1 \le |\lambda| \le r_2\}, \\
\Omega_3 =& \{\lambda \in \mathbb{C} : |\operatorname{arg} \lambda| \le \frac{2\pi}{3}, \, |\lambda| > r_2 \}.
\end{aligned}
\end{equation*}
where $r_1, r_2 $ are positive constants to be determined and $r_1$  is sufficiently small.

Hypotheses (H1)-(H4) and assumptions (A1)-(A3) are verified in the work of Deng \cite{Den91b} and Sandstede \cite{San}. Then by applying Corollary \ref{cor3.1}, there exists a constant $r_1>0$ such that for $|\lambda|<r_1$, $\la$ is an spectrum of $\ml$ if and only if $\la$ satisfies the expression \eqref{314}.
We choose \(\psi_1(x)\),  \(\psi_2(x)\) such that \(\langle w_1^+, v_1^+ \rangle > 0\) and \(\langle w_2^+, v_2^+ \rangle > 0\). Since the aforementioned heteroclinic loop is double-twisted, i.e.,  \(\langle w_1^+, v_1^+ \rangle \) and
\(\langle w_2^-, v_2^- \rangle \) have the same sign, and
\(\langle w_2^+, v_2^+ \rangle \) and \(\langle w_1^-, v_1^- \rangle \) have the same sign. Thus, we can assume that
\(\langle w_i^{+}, v_i^{+} \rangle \) and \(\langle w_i^{-}, v_i^{-} \rangle \) for $i=1,2$ are all positive values. With this choice, it follows from pp.206-207 in \cite{San}  that
$$M_i = \int_{-\infty}^{\infty} \langle \psi_i(x), Bh_i'(x) \rangle dx = \int_{-\infty}^{\infty} \langle \psi_i(x), F_c(h_i(x), c) \rangle dx < 0,$$
for $i= 1,2$.
Thus, we have $$\langle w_2^+, v_2^+ \rangle M_1 + \langle w_1^+, v_1^+ \rangle M_2< 0.$$
Using expressions (\ref{314}), we can obtain $\lambda(\xi)< 0$ for $\xi \in [-\pi/T,\pi/T)$.
That is, the spectrum of the periodic solutions $p(x)$ near the origin lies in the left half of the complex plane.
 It is standard that there exists a large constant $r_2>0$ such that the unbounded subset $\Omega_3$ and the compact subset $\Omega_2$ contain no spectrum of $\ml$, see \cite{CdS16,Mag}. Then the proof is completed.
 \end{proof}


\begin{remark}
In the above discussion, we applied Theorem \ref{thm3.2} to study the
spectral stability of a FitzHugh-Nagumo model. In fact, this theorem can also be applied to many other mathematical models that have been shown to possess heteroclinic loops. For instance, M. Engel, C. Kuehn, and B. de Rijk\cite{EKR} explored a turbulent pipe flow model, which is a mixed system combining a bistable reaction-diffusion system with advective terms of Burgers' type. Their research demonstrates that, in different parameter regions, the model admits heteroclinic loops with single twist or double twist. Another interesting example is the work by Li and Yu\cite{LY}, considering a reaction-diffusion mechanical model derived from coupling a modified FitzHugh-Nagumo model with a mechanical system. This model allows for the existence of heteroclinic loops with double twist.
\end{remark}

%
%

Recall that we denote the periodicity by $T$. As mentioned in foregoing sections, a complex number $\la$ lies in $L^2(\mathbb{R})$-spectrum of $\ml$ if and only if
\begin{align}\label{eigen}
\ml V=\la V, \quad
V(T)=e^{i\xi T}V(0),
\end{align}
has a nontrivial solution for some $\xi\in[-\pi/T,\pi/T)$.
Substituting $V(x)=e^{i\xi x}W(x)$ into \eqref{eigen}, we deduce that \eqref{eigen} has a nontrivial solution if and only if
\begin{align}\label{eigen2}
\la W=e^{-i\xi x}\ml e^{i\xi x}W=: \lxi W,
\end{align}
has a nontrivial solution $W\in L^2_{per}(0,T)$. Thus, we have the following spectral decomposition
$$\sigma(\ml)=\bigcup_{\xi\in[-\pi/T,\pi/T)}\sigma(\lxi),$$
where
$$\lxi=\left(\begin{array}{cc}
(\partial_{x}+i\xi)^2-c(\partial_x+i\xi)+f'(\phi)&-1\\
\e&-c(\partial_x+i\xi)-\e\gamma
\end{array}\right),$$
 acts on $L^2_{per}(0,T)$. By virtue of the transform $V(x)=e^{i\xi x}W(x)$, $\la$ is an eigenvalue of $\lxi$ with $L^2_{per}(0,T)$-eigenfunction $W$ if and only if $V$ is a nontrivial solution of \eqref{207}-\eqref{208}.

Inspired by the characterization of spectrum of $\ml$, one intends to decompose a function in $L^2(\mathbb{R})$ into the superpositions of functions of the form $e^{i\xi x}W(x)$. For a function $g\in L^2(\mathbb{R})$, we introduce the inverse Bloch-Fourier representation
$$g(x)=\frac{1}{2\pi}\int_{-\pi/T}^{\pi/T} e^{i\xi x} \check{g}(\xi,x) d\xi, $$
where $\check{g}(\xi,x)=\sum_{k\in\mathbb{Z}} e^{2\pi ilx/T}\hat{g}(\xi+2\pi lx/T)$ is a T-periodic function of $x$ and $\hat{g}(\cdot)$ denoting the usual Fourier transformation of $g$, i.e.,
$$\hat{g}(\xi)=\int_{-\infty}^{\infty}e^{-i\xi x}g(x)dx.$$
Then, for fixed $m\in\mathbb{N}_0$, the Bloch transform
$$\mathcal{B}: H^m(\mathbb{R})\to L^2([-\pi/T,\pi/T);H^m_{per}(0,T)),\quad \mathcal{B}g=\check{g},$$
is a linear bounded operator, since
\begin{equation}\label{Par}\begin{split}
\|g\|_{H^m(\mathbb{R})}^2&\simeq \int_{-\infty}^\infty (1+\xi^2)^m|\hat{g}(\xi)|^2d\xi\\
&\simeq \sum_{l\in\mathbb{Z}}\int_{-\pi/T}^{\pi/T}(1+(2\pi l/T)^2)^m|\hat{g}(\xi+2\pi l/T)|^2d\xi\\
&\simeq \int_{-\pi/T}^{\pi/T} \|\check{g}(\xi,\cdot)\|^2_{H^m_{per}(0,T)}d\xi\\
&=\|\check{g}\|^2_{L^2([-\pi/T,\pi/T);H^m_{per}(0,T))},
\end{split}\end{equation}
where the notation $``A \simeq B"$ means that $A\lesssim B$ and $B\lesssim A$. Specifically, for $m=0$, we have \begin{align}\label{Parseval}\|g\|^2_{L^2(\mathbb{R})}=\frac{1}{2\pi T}\|\check{g}\|^2_{L^2([-\pi/T,\pi/T);L^2_{per}(0,T))}.\end{align}
Then for $U\in D(\ml)$, we have the representation formula
$$\ml U(x)=\frac{1}{2\pi}\int_{-\pi/T}^{\pi/T} e^{i\xi x} \lxi\check{U}(\xi,x) d\xi.$$
Moreover, the operator $\ml$ can also be regarded as
$$\ml:H^{m+2}({\mathbb{R}})\times H^{m+1}({\mathbb{R}})\subset H^{m}({\mathbb{R}})\times H^{m}({\mathbb{R}})\to H^{m}({\mathbb{R}})\times H^{m}({\mathbb{R}}), \text{ for } m\in\mathbb{N}_0. $$
Correspondingly, the Bloch operators $\lxi$ defined as
$$\lxi:H^{m+2}_{per}(0,T)\times H^{m+1}_{per}(0,T)\subset H^{m}_{per}(0,T)\times H^{m}_{per}(0,T) \rightarrow H^{m}_{per}(0,T)\times H^{m}_{per}(0,T).$$
Note that the periodic solution $\mathbf{\bar{U}}$ is smooth, and the Bloch operators $\lxi$ are closed in $L^2_{per}(0,T)\times L^2_{per}(0,T)$ and $H^m_{per}(0,T)\times H^m_{per}(0,T)$ for $m\in\mathbb{N}_+$ with  compactly embedding domains $H^{m+2}_{per}(0,T)\times H^{m+1}_{per}(0,T)$. Then, by a standard bootstrapping argument, the Bloch operators $\lxi$ have the same spectral properties when acting on $H^m_{per}(0,T)\times H^m_{per}(0,T)$  for all $m\in\mathbb{N}_0$. After the above preliminaries, we state more details about the spectrum of $\lxi$.

\begin{lemma}\label{spec}
For fixed $m\in\mathbb{N}_0$, the Bloch operators $\lxi:H^{m+2}_{per}(0,T)\times H^{m+1}_{per}(0,T)\subset H^{m}_{per}(0,T)\times H^{m}_{per}(0,T)\to  H^{m}_{per}(0,T)\times H^{m}_{per}(0,T)$ have the following spectral properties:
\begin{enumerate}[(i).]
  \item For any fixed $\xi_0\in(0,\pi/T)$, there exists a $\delta_0>0$ such that
      $$\Re \sigma(\lxi)<-\delta_0,\text{ for }\xi\in[-\pi/T,\pi/T) \text{ with } |\xi|>\xi_0.$$
  \item There exist constants $\xi_1\in[0,\pi/T)$ and $\delta_1>0$ such that the spectrum of $\lxi$ can be decomposed into two disjoint subsets
      $$\sigma(\lxi)=\sigma_s(\lxi)\cup \sigma_c(\lxi), \text{ for any } |\xi|<\xi_1,$$
      with the following properties:
      \begin{enumerate}[(a).]
        \item $\Re \sigma_s(\lxi)<-\delta_1$ and $\Re \sigma_c(\lxi)>-\delta_1$.
        \item $\sigma_c(\lxi)$ consists of a simple eigenvalue $\la_c(\xi)$, which is analytic in $\xi$ with the following expansion:
            $$\la_c(\xi)=ib\xi-d\xi^2+O(|\xi|^3), \text{ for some } b\in\mathbb{R} \text{ and } d>0.$$
        \item the eigenfunction $\Phi_\xi(x)$ associated with $\la_c(\xi)$ is smooth in $x$ and analytic in $\xi$, satisfying
            $$\|\Phi_\xi-\mathbf{\bar{U}}'\|_{H^m_{per}(0,T)}\leq C|\xi|, \text{ for some constant }C>0,$$
            where $\mathbf{\bar{ U}}'=(\phi',\psi')^t$ is the derivative of periodic traveling wave.
      \end{enumerate}
\end{enumerate}
\end{lemma}
\begin{proof}
It is sufficient to prove this lemma for $m=0$. The statements (i), (ii)(a) and (ii)(b) (except for simplicity of $\la_c(\xi)$) can be obtained by \eqref{308}, which gives the expansion of eigenvalue $\la_c(\xi)$ with respect to $\xi$. 
The algebraic multiplicity of $\la_c(\xi)$ is equal to the multiplicity of $\la$ as a zero of the function $E(\la,\xi)$, which called Evans function, see \cite{KP13,LY,San}. Thus, $\la_c(\xi)$ is a simple eigenvalue of $\lxi$ by virtue of $M_1\cdot M_2\neq0$.
Note that $\lxi$ depends continuously on $\xi$. Then the statement (ii)(c) follows from the conclusions (ii)(a)-(b) and standard spectral perturbation theory \cite{kato}.
\end{proof}

\subsection{Semigroup estimates}\label{Seclin}
In this section, we prove that $\ml$ generates a $\mathcal{C}^0$-semigroup $e^{\ml t}$ by using a bounded perturbation theory and the Bloch operators generate $\mathcal{C}^0$-semigroups by applying Gearhart-Pr\"{u}ss theorem. Furthermore, decomposing the semigroup $e^{\ml t}$ into exponentially decaying part $S_e(t)$ and the critical part $S_c(t)$, we obtain the following linear stability (Theorem \ref{linearsta}) of the periodic traveling wave $\mathbf{\bar{U}}=(\phi,\psi)^t$.
\begin{theorem}\label{linearsta}
There exists a constant $C>0$ such that
$$\|e^{\ml t}g\|_{L^2}\leq C(1+t)^{-\frac{1}{4}}\|g\|_{L^1\cap L^2}, \quad \text{ for all }g\in L^1(\mathbb{R})\cap L^2(\mathbb{R}) \text{ and } t\geq0.$$
Moreover, we have
$$\|e^{\ml t}g-\mathbf{\bar{U}}'s_p(t)g\|_{L^2}\leq C(1+t)^{-\frac{3}{4}}\|g\|_{L^1\cap L^2}, \quad \text{ for all }g\in L^1(\mathbb{R})\cap L^2(\mathbb{R}) \text{ and } t\geq0,$$
where $s_p(t)g$ is defined by \eqref{sp}.
\end{theorem}

\subsubsection{Semigroup properties}
\begin{theorem}[Gearhart-Pr\"{u}ss \cite{Pru84}]
Assume the linear operator $\mathcal{A}:\mathcal{D(A)}\subset X\to X$ generates a $\mathcal{C}^0$-semigroup on Hilbert space X and the spectra of $\mathcal{A}$ satisfy $\sigma(\mathcal{A})\subset \{\la\in\mathbb{C}:\Re \la<0\}$ . If
$$\sup_{\mu\in\mathbb{R}}\|(\mathcal{A}-i\mu)^{-1}\|_{\mathbf{L}(X)}<\infty,$$
then there exists $\delta>0$ such that
$$\|e^{\mathcal{A}t}\|_{\mathbf{L}(X)}\leq C_\delta e^{-\delta t}.$$
\end{theorem}

\begin{prop}\label{group}
The linear operator $\mathcal{L}$ and the Bloch operators $\mathcal{L}_\xi,\xi\in[-{\pi}/{T},{\pi}/{T})$ have the following semigroup properties.
\begin{enumerate}[(i).]
\item For $m\in\mathbb{N}_0$, $\mathcal{L}$ generates a $\mathcal{C}^0$-semigroup on $H^{m}(\mathbb{R})$.
\item For any fixed $\xi_0\in(0,\pi/T)$ and $m\in\mathbb{N}_0$, there exist positive constants $C_0$ and $\mu_0$ such that
    $$\|e^{\mathcal{L}_\xi t}\|_{\mathbf{L}(H^m_{per}(0,T))}\leq C_0 e^{-\mu_0 t}, \quad \text{for all }t\geq0 \text{ and } |\xi|>\xi_0.$$
\item Suppose $\Pi(\xi)$ is the spectral projection onto the eigenspace associated to the eigenvalue $\la_c(\xi)$. Then there exist positive constants $C_1$ and $\mu_1$ such that
    $$\|e^{\mathcal{L}_\xi t}(I-\Pi(\xi))\|_{\mathbf{L}(H^m_{per}(0,T))}\leq C_1 e^{-\mu_1 t}, \quad \text{for all }t\geq0 \text{ and } |\xi|<\xi_1,$$
    where $\xi_1$ is defined in Lemma \ref{spec}.
\end{enumerate}
\end{prop}
\begin{proof}
(i). We rewrite the linear operator $\mathcal{L}$ as
      $$\ml=\tilde{\ml}+\left(\begin{array}{cc}f'(\phi)&-1\\\e&-\e\gamma
      \end{array}\right),$$
      where $\tilde{\ml}=\left(\begin{array}{cc}\partial_{x}^2-c\partial_x&0\\0&-c\partial_x
      \end{array}\right)$. $\tilde{\ml}$ generates a $\mathcal{C}^0$-semigroup on $H^m(\mathbb{R})\times H^m(\mathbb{R})$ for $m\in\mathbb{N}_0$ since both $\partial_{xx}+c\partial_x$ and $c\partial_x$ generate $\mathcal{C}^0$-semigroup on $H^m(\mathbb{R})\times H^m(\mathbb{R})$; see \cite{KP13}. Note that $\ml-\tilde{\ml}$ is a bounded operator defined on the space $H^m(\mathbb{R})\times H^m(\mathbb{R})$, then the conclusion follows from Theorem 3.1.1 of \cite{kato}.

(ii). For convenience, we use the notation
$$\hat{\partial}_x^m=(\partial_x+i\xi)^m,\quad m\in\mathbb{N}_+.$$
Clearly, we have the following identity 
$$\langle h,\hat{\partial}_xg\rangle_{L^2}=-\langle\hat{\partial}_x h,g\rangle_{L^2}, \quad \text{ for }h,g\in H^1_{per}(0,T),$$
which is derived by by parts. Specially, we have
\begin{equation}\label{hh}
\langle h,\hat{\partial}_xh\rangle_{L^2}\in i\mathbb{R}, \quad \text{ for }h\in H^1_{per}(0,T).
\end{equation}
Consider $$(\lxi-\la)\mathbf{U}=\mathbf{F},$$
where $\mathbf{U}=(u,w)^t\in H^2_{per}(0,T)\times H^1_{per}(0,T)$ and $\mathbf{F}=(f_1,f_2)^t\in L^2_{per}(0,T)\times L^2_{per}(0,T)$.
Recall that 
$$\|\mathbf{U}\|_{L^2}^2=\langle u,u\rangle_{L^2}+\langle w,w\rangle_{L^2}.$$
 Let $\la\in i\mathbb{R}\setminus\{0\}$. Using the identity \((\lxi-\la)\mathbf{U}=\mathbf{F}\), we obtain
\begin{equation}\label{E0}
\begin{split}0=\Re \la\cdot \|\mathbf{U}\|^2_{L^2}=&\Re( \langle u, \la u\rangle_{L^2} +\langle w, \la w\rangle_{L^2})\\
=&\Re( \langle u, \hat{\partial}_{x}^2u-c\hat{\partial}_xu+f'(\phi)u-w-f_1\rangle_{L^2} +\langle w, -c\hat{\partial}_xw+\e(u-\gamma w)-f_2\rangle_{L^2})\\
\overset{\eqref{hh}}{\leq}&-\|\hat{\partial}_x u\|_{L^2}^2-\e\gamma\|w\|_{L^2}^2+C(\|u\|_{L^2}^2+(\|u\|_{L^2}+\|f_2\|_{L^2})\|w\|_{L^2}+\|\mathbf{F}\|_{L^2}^2)\\
\leq&-\|\hat{\partial}_x u\|_{L^2}^2-\frac{\e\gamma}{2}\|\mathbf{U}\|_{L^2}^2+C(\|u\|_{L^2}^2+\|\mathbf{F}\|_{L^2}^2),
\end{split}\end{equation}where we use the Young's inequality $\|w\|_{L^2}\|f_2\|_{L^2}\leq \frac{\e\gamma}{2}\|w\|_{L^2}^2+\frac{2}{\e\gamma}\|f_2\|_{L^2}^2$ for sufficiently large $\tilde{C}$. On the other hand, a direct calculation shows that
\begin{equation}\label{u}
\begin{split}
\|u\|_{L^2}^2=\frac{1}{\la}\langle u,\la u\rangle_{L^2}\leq \frac{C}{|\la|}(\|\hat{\partial}_x u\|_{L^2}^2+\|\mathbf{U}\|_{L^2}^2+\|\mathbf{F}\|_{L^2}^2),\quad \text{for }\la\neq0.
\end{split}
\end{equation}
Combining \eqref{E0} and \eqref{u}, we conclude that there exists a constant $\la_0>0$ such that for all $\la$ with $|\la|>\la_0$, we have $\|\mathbf{U}\|_{L^2}^2\leq C(\la_0)\|\mathbf{F}\|_{L^2}^2$. That is, 
$$\|\mathbf{U}\|_{L^2}^2\leq C(\la_0)\|(L_\xi-\la)\mathbf{U}\|_{L^2}^2, \text{ for } |\la|>\la_0.$$
 By using Gearhart-Pr\"{u}ss theorem, we obtain the semigroup decay estimates on $L_{per}^2(0,T)$.
The proof of $H^m,m\in\mathbb{N}_+$ is almost the same, except for a slight difference.

We define the norm
$$\|g\|_{H^m_\xi(0,T)}:=\left(\|g\|_{L^2(0,T)}^2+\sum_{k=1}^m\|\hat{\partial}_x^mg\|_{L^2(0,T)}^2\right)^{1/2},$$
which is equivalent to the usual normal $\|g\|_{H^m(0,T)}$ for $\xi\in[-\pi/T,\pi/T).$ Consider $$(\lxi-\la)\mathbf{U}=F,$$
where $\mathbf{U}=(u,w)^t\in H^{m+2}_{per}(0,T)\times H^{m+1}_{per}(0,T)$ and $\mathbf{F}=(f_1,f_2)^t\in H^m_{per}(0,T)\times H^m_{per}(0,T)$. 
Define 
$$\|\mathbf{U}\|^2_{H^m_\xi}:=\sum_{k=0}^m(\langle \hat{\partial}^k_x u,\hat{\partial}^k_x u\rangle_{L^2}+\langle \hat{\partial}^k_xw,\hat{\partial}^k_xw\rangle_{L^2}),$$ which are real-valued.
 Then\begin{equation}\label{Em}
\begin{split}0=\Re \la\cdot \|\mathbf{U}\|^2_{H^m_\xi}=&\Re\sum_ {k=0}^m( \langle \hat{\partial}^k_x u, \la\hat{\partial}^k_x u\rangle_{L^2} +\langle \hat{\partial}^k_x w, \la \hat{\partial}^k_x w\rangle_{L^2})\\
=&\Re\sum_ {k=0}^m( \langle \hat{\partial}^k_x u, \hat{\partial}^k_x(\hat{\partial}_{x}^2u-c\hat{\partial}_xu+f'(\phi)u-w-f_1)\rangle_{L^2} \\ &+\Re\sum_ {k=0}^m(\langle \hat{\partial}^k_xw, \hat{\partial}^k_x(-c\hat{\partial}_xw+\e(u-\gamma w)-f_2)\rangle_{L^2})\\
\leq&-\|\hat{\partial}_x^{m+1} u\|_{L^2}^2-\sum_{k=1}^m
\|\hat{\partial}_x^{k} u\|_{L^2}^2-\e\gamma\|w\|_{H^m_\xi}^2\\
&+C(\|u\|_{H^m_\xi}^2+(\|u\|_{H^m_\xi}+\|f_2\|_{H^m_\xi})(\|w\|_{H^m_\xi}+\|f_1\|_{H^m_\xi}))\\
\leq&-\|\hat{\partial}_x^{m+1} u\|_{L^2}^2-\sum_{k=1}^m
\|\hat{\partial}_x^{k} u\|_{L^2}^2-\frac{\e\gamma}{2}\|w\|_{H^m_\xi}^2\\
&+C(\|u\|_{H^m_\xi}^2+(\|u\|_{H^m_\xi}+\|f_2\|_{H^m_\xi})\|f_1\|_{H^m_\xi}+\|f_2\|_{H^m_\xi}^2)\\
\leq&-\frac{\e\gamma}{2}\|\mathbf{U}\|_{H^m_\xi}^2+C(\|u\|_{L^2}^2+\|\mathbf{F}\|_{H^m_\xi}^2),
\end{split}\end{equation}
where we use the Young's inequality $ab\leq \tilde{C}a^2+\frac{1}{\tilde{C}}b^2$ and the Sobolev inequality $\|g\|_{H^m_\xi}^2\leq \frac{1}{\tilde{C}}\|\hat{\partial}_x^{m+1}g\|_{L^2}^2+\tilde{C}\|g\|_{L^2}^2$ for sufficiently large $\tilde{C}$. By using \eqref{u}, we also have
\begin{equation}\label{u1}
\begin{split}
\|u\|_{L^2}^2\leq \frac{C}{|\la|}(\| \mathbf{U}\|_{H^m_\xi}^2+\|\mathbf{F}\|_{H^m_\xi}^2),\quad \text{for }\la\neq0 \text{ and } m\in\mathbb{N}^+.
\end{split}
\end{equation}
Then, combining \eqref{Em} and \eqref{u1}, we conclude that
$$\|\mathbf{U}\|_{H^m_\xi}^2\leq C \|\mathbf{F}\|_{H^m_\xi}^2, \text{ for } |\la| \text{ sufficiently large.}$$

(iii). Note that the estimates \eqref{E0} and \eqref{Em} also apply to $\lxi$ for $|\xi|<\xi_1$. In addition, $\lxi(I-\Pi(\xi))$ satisfies the spectral conditions in Gearhart-Pr\"{u}ss Theorem. Then we have the semigroup estimates as we wanted. 
\end{proof}

\subsubsection{Semigroup decomposition}
By virtue of Proposition \ref{group}, it is straightforward to verify that
\begin{align}\label{eLt}e^{\ml t}U(x)=\frac{1}{2\pi}\int_{-\pi/T}^{\pi/T}e^{i\xi x}e^{\lxi t}\check{U}(\xi,x)d\xi, \text{ for } U\in \mathcal{D}(\ml).\end{align}
Then, by using the exponentially decaying estimates proved in Proposition \ref{group}, we decompose the $\mathcal{C}^0$-semigroup $e^{\ml t}$ into an exponentially decaying part and a critical part, which turns out to exhibit an algebraic decaying rate. Introduce a smooth nonnegative cutoff function $\rho(\xi)$ as 
$$\rho(\xi):=\left\{\begin{array}{ll}1,\quad |\xi|<{\xi_1}/2,\\0,\quad |\xi|>\xi_1.\end{array}\right.$$
Then, we decompose $e^{\ml t}$ into
\begin{align}
e^{\ml t}U(x)=S_e(t)U(x)+S_c(t)U(x),
\end{align}
where \begin{equation*}\begin{split}
S_e(t)U(x):=&\frac{1}{2\pi}\int_{-\pi/T}^{\pi/T}(1-\rho(\xi))e^{i\xi x}e^{\lxi t}\check{U}(\xi,x)d\xi\\
&+\frac{1}{2\pi}\int_{-\pi/T}^{\pi/T}\rho(\xi)e^{i\xi x}e^{\lxi t}(1-\Pi(\xi))\check{U}(\xi,x)d\xi,
\end{split}
\end{equation*}
and
\begin{equation*}\begin{split}
S_c(t)U(x):=\frac{1}{2\pi}\int_{-\pi/T}^{\pi/T}\rho(\xi)e^{i\xi x}e^{\lxi t}\Pi(\xi)\check{U}(\xi,x)d\xi.
\end{split}
\end{equation*}
The following lemma follows from Proposition \ref{group} and the equality \eqref{Par}.

\begin{lemma}\label{exp}
For any $m\in\mathbb{N}_0$, there exist constants $C,\mu>0$ such that
$$\|S_e(t)\|_{\mathbf{L}(H^{m}(\mathbb{R}))}\leq Ce^{-\mu t}, \text{ for all }t\geq0.$$
\end{lemma}

Next, we further discuss the critical part $S_c(t)$. Recall that $\Pi(\xi)$ is the spectral projection onto the one-dimensional eigenspace of $\lxi$ associated to $\la_c(\xi)$. Suppose $\widetilde{\Phi}_\xi$ is the smooth eigenfunction of adjoint operator $\lxi^*$ associated to eigenvalue $\overline{\la_c(\xi)}$ which satisfies $\langle {\Phi}_\xi, \widetilde{\Phi}_\xi\rangle_{L^2}=1$.  Then for any $g\in L_{per}^2(0,T)$, we have
$$\Pi(\xi)g=\langle g, \widetilde{\Phi}_\xi\rangle_{L^2}\Phi_\xi.$$

In order to further analyze $S_c(t)$, we introduce a smooth cutoff function $\chi:[0,\infty)\to \mathbb{R}$ as
$$\chi(t):=\left\{\begin{array}{ll}0,\quad t\in[0,1],\\1,\quad t\in[2,\infty).\end{array}\right.$$
Then we rewrite $S_c(t)$ as
\begin{equation}\begin{split}
S_c(t)U(x)=&\ \chi(t)S_c(t)U(x)+(1-\chi(t))S_c(t)U(x),\\
=&\ \frac{\chi(t)}{2\pi}\int_{-\pi/T}^{\pi/T}\rho(\xi)e^{i\xi x+\la_c(\xi)t}\langle \check{U}(\xi,\cdot), \widetilde{\Phi}_\xi\rangle_{L^2} \Phi_\xi(x)d\xi+(1-\chi(t))S_c(t)U(x)  \\
=&\ \mathbf{\bar{U}}'(x)\left(\frac{\chi(t)}{2\pi}\int_{-\pi/T}^{\pi/T}\rho(\xi)e^{i\xi x+\la_c(\xi)t}\langle \check{U}(\xi,\cdot), \widetilde{\Phi}_\xi\rangle_{L^2} d\xi\right)\\
&\ + \frac{\chi(t)}{2\pi}\int_{-\pi/T}^{\pi/T}\rho(\xi)e^{i\xi x+\la_c(\xi)t}i\xi \left(\frac{\Phi_\xi(x)-\mathbf{\bar{U}}'(x)}{i\xi}\right)\langle \check{U}(\xi,\cdot), \widetilde{\Phi}_\xi\rangle_{L^2} d\xi\\
&\ +(1-\chi(t))S_c(t)U(x)\\
=&\mathbf{\bar{U}}'(x)s_p(t)U(x)+\tilde{S}_c(t)U(x),
\end{split}
\end{equation}
where
\begin{equation}\label{sp}s_p(t)U(x):=\frac{\chi(t)}{2\pi}\int_{-\pi/T}^{\pi/T}\rho(\xi)e^{i\xi x+\la_c(\xi)t}\langle \check{U}(\xi,\cdot), \widetilde{\Phi}_\xi\rangle_{L^2} d\xi.\end{equation}
The following lemma shows that $s_p(t)$ exhibits the slowest decaying rate $(1+t)^{-1/4}$ and $\tilde{S}_c(t)$ exhibits the faster decaying rate $(1+t)^{-3/4}$.

\begin{lemma}\label{crit}
For any fixed $l,j,m\in\mathbb{N}_0$, there exists 
 a constant $C>0$ such that for all $t\geq0$, the following inequalities hold:
\begin{align}
\|\partial_x^l\partial_t^js_p(t)U\|_{L^2}&\leq C(1+t)^{-\frac{l+j}{2}}\|U\|_{L^2},\quad U\in L^2(\mathbb{R}),\\
\|\partial_x^l\partial_t^js_p(t)U\|_{L^2}&\leq C(1+t)^{-\frac{1}{4}-\frac{l+j}{2}}\|U\|_{L^1},\quad U\in L^1(\mathbb{R})\cap L^2(\mathbb{R}),\\
\|\partial_x^m\tilde{S}_c(t)U\|_{L^2}&\leq C(1+t)^{-\frac{3}{4}}\|U\|_{L^1},\quad U\in L^1(\mathbb{R})\cap L^2(\mathbb{R}).
\end{align}
\end{lemma}
\begin{proof}
Consider $t\geq2$. By definition of $\chi(t)$, we have $\chi(t)=1$ and $\chi^{(j)}(t)=0$ for $j\in\mathbb{N}$ and $t\geq2$. As a consequence, we have
$$\partial_x^l\partial_t^js_p(t)U(x)=\frac{1}{2\pi}\int_{-\pi/T}^{\pi/T}\rho(\xi)(i\xi)^l(\la_c(\xi))^je^{i\xi x+\la_c(\xi)t}\langle \check{U}(\xi,\cdot), \widetilde{\Phi}_\xi\rangle_{L^2} d\xi.$$
Then, by virtue of \eqref{Parseval}, it holds that
$$\|\partial_x^l\partial_t^js_p(t)U\|_{L^2}^2=\frac{1}{2\pi T}\int_{-\pi/T}^{\pi/T}\int_0^T|\rho(\xi)(i\xi)^l(\la_c(\xi))^je^{\la_c(\xi)t}\langle \check{U}(\xi,\cdot), \widetilde{\Phi}_\xi\rangle_{L^2} |^2dxd\xi,$$
for $U\in L^2(\mathbb{R})$. Recall that $\la_c(\xi)=ib\xi-d\xi^2+O(|\xi|^3)$. Thus, we conclude
\begin{align*}
\|\partial_x^l\partial_t^js_p(t)U\|_{L^2}&\lesssim \max_{\xi\in[-\frac{\pi}{T},\frac{\pi}{T}]}(\xi^{l+j}e^{-d\xi^2t})\cdot\|U\|_{L^2}\\
&\lesssim (1+t)^{-\frac{l+j}{2}}\|U\|_{L^2},\quad U\in L^2(\mathbb{R}),\end{align*}
where we use 
\begin{align}\label{l2}
|\langle \check{U}(\xi,\cdot), \widetilde{\Phi}_\xi\rangle_{L^2}|\lesssim \|\check{U}(\xi,\cdot)\|_{L^2(0,T)}, \quad \text { for } U\in L^2(\mathbb{R}).
\end{align}
Moreover, we have
\begin{align*}
\|\partial_x^l\partial_t^js_p(t)U\|_{L^2}&\lesssim \left(\int_{-\pi/T}^{\pi/T}|\xi^{l+j}e^{-d\xi^2t}|^2d\xi\right)^{1/2}\cdot\|U\|_{L^1}\\
&\lesssim (1+t)^{-\frac{1}{4}-\frac{l+j}{2}}\|U\|_{L^1},\quad U\in L^1(\mathbb{R})\cap L^2(\mathbb{R}),
\end{align*}
where we use \begin{align}\label{l1}
|\langle \check{U}(\xi,\cdot), \widetilde{\Phi}_\xi\rangle_{L^2}|^2&\lesssim \|\hat{U}\|_{L^{\infty}(\mathbb{R})}\lesssim\|U\|_{L^1(\mathbb{R})}^2,\quad\quad \text{ for }U\in L^1(\mathbb{R})\cap L^2(\mathbb{R}).
\end{align}
Thus, we have proved the first two inequalities for $t\geq 2$. As for the third inequality, we have
\begin{align*}
\|\partial_x^m\tilde{S}_c(t)U\|_{L^2}&\lesssim \left(\sum_{k=1}^{m+1} \int_{-\pi/T}^{\pi/T}|\xi^{k}e^{-d\xi^2t}|^2d\xi\right)\cdot \|U\|_{L^1}\\
&\lesssim (1+t)^{-\frac{3}{4}}\|U\|_{L^1}, \quad \text{ for }t\geq2,
\end{align*}
by virtue of the fact that
$$\sup_{\xi\in[-\xi_1,\xi_1)}\left\|\partial_x^k\left(\frac{\Phi_\xi-\mathbf{\bar{U}}}{i\xi}\right)\right\|_{L^\infty}<\infty, \text{ for }k=0,1,...,m.$$
Then the third inequality holds for $t\geq2$.

 Note that $\chi(t)=0$ for $t\in[0,1]$, and $\chi^{(j)}(t), j\in \mathbf{N}_0$ are bounded for $t\in[1,2]$. Then a similar calculation shows the first these inequalities hold for $t\in[0,2]$. 

\end{proof}

\begin{remark}
    Lemma \ref{crit} shows that $S_c(t)$ is infinitely smoothing in the sense that it defines a bounded map from $L^1(\mathbb{R})\cap L^2(\mathbb{R})$ into $H^m(\mathbb{R})$, even for $t=0$. However, the exponentially decaying part $S_e(t)$ is not infinitely smoothing, see Lemma \ref{exp}.
\end{remark}

\subsubsection{Integration by parts}
Prepared for the nonlinear scheme, we give the following integration-by-parts-type estimates.

\begin{prop}\label{part}
For fixed $l,j\in\mathbb{N}_0$, the following inequalities hold.
\begin{enumerate}[(i).]
  \item For $g=\left(\begin{array}{cc}g_1\\g_2\end{array}\right),h=\left(\begin{array}{cc}h_1\\h_2\end{array}\right)\in H^1(\mathbb{R})$ and $t\geq0$, we have
      \begin{align*}
      \left\|\partial_x^l\partial_t^j s_p(t)\left(\begin{array}{cc}g_1\cdot\partial_xh_1\\g_2\cdot\partial_xh_2\end{array}\right)\right\|_{L^2}&\lesssim (1+t)^{-\frac{1}{4}-\frac{l+j}{2}}\sum_{i=1}^{2}(\|g_ih_i\|_{L^1}+\|\partial _xg_i\cdot h_i\|_{L^1})\\
      &\lesssim(1+t)^{-\frac{1}{4}-\frac{l+j}{2}}\|g\|_{H^1}\|h\|_{L^2},\\
      \left\|\tilde{S}_c(t)\left(\begin{array}{cc}g_1\cdot\partial_xh_1\\g_2\cdot\partial_xh_2\end{array}\right)\right\|_{L^2}&\lesssim (1+t)^{-\frac{3}{4}}\sum_{i=1}^2(\|g_ih_i\|_{L^1}+\|\partial_x g_i\cdot h_i\|_{L^1})\\
      &\lesssim (1+t)^{-\frac{3}{4}}\|g\|_{H^1}\|h\|_{L^2}.
      \end{align*}
  \item For $g=\left(\begin{array}{cc}g_1\\g_2\end{array}\right),h=\left(\begin{array}{cc}h_1\\h_2\end{array}\right) \in H^2(\mathbb{R})$ and $t\geq0$, we have
      \begin{align*}
      \left\|\partial_x^l\partial_t^j s_p(t)\left(\begin{array}{cc}g_1\cdot\partial_x^2h_1\\g_2\cdot\partial_x^2h_2\end{array}\right)\right\|_{L^2}&\lesssim (1+t)^{-\frac{1}{4}-\frac{l+j}{2}}\sum_{i=1}^2(\|g_ih_i\|_{L^1}+\|\partial_xg_i\cdot h_i\|_{L^1}+\|\partial_x^2 g_i\cdot h_i\|_{L^1})\\
      &\lesssim (1+t)^{-\frac{1}{4}-\frac{l+j}{2}}\|g\|_{H^2}\|h\|_{L^2},\\
      \left\|\tilde{S}_c(t)\left(\begin{array}{cc}g_1\cdot\partial_x^2h_1\\g_2\cdot\partial_x^2h_2\end{array}\right)\right\|_{L^2}&\lesssim (1+t)^{-\frac{3}{4}}\sum_{i=1}^2(\|g_ih_i\|_{L^1}+\|\partial_xg_i\cdot h_i\|_{L^1}+\|\partial_x^2 g_i\cdot h_i\|_{L^1})\\
      &\lesssim  (1+t)^{-\frac{3}{4}}\|g\|_{H^2}\|h\|_{L^2}.
      \end{align*}
\end{enumerate}
\end{prop}
\begin{proof}
By the definition of Bloch transform, if $g\in H^1(\mathbb{R})$, then $\mathcal{B}(g)(\xi,\cdot)$ is differentiable and
$$\partial_x\mathcal{B}(g)(\xi,x)=\mathcal{B}(\partial_xg)(\xi,x)-i\xi\check{g}(\xi,x).$$
Then integration by parts shows that
\begin{align*}\langle \mathcal{B}(\left(\begin{array}{cc}g_1h_1\\g_2h_2\end{array}\right))(\xi,\cdot), \partial_x\widetilde{\Phi}_\xi\rangle_{L^2}=&i\xi\langle \mathcal{B}(\left(\begin{array}{cc}g_1h_1\\g_2h_2\end{array}\right))(\xi,\cdot), \widetilde{\Phi}_\xi\rangle_{L^2}-\langle \mathcal{B}(\left(\begin{array}{cc}\partial_xg_1\cdot h_1+g_1\partial_x\cdot h_1\\\partial_xg_2\cdot h_2+g_2\partial_x\cdot h_2\end{array}\right))(\xi,\cdot), \widetilde{\Phi}_\xi\rangle_{L^2}.
\end{align*}
As a consequence, we have
\begin{align*}
s_p(t)(\left(\begin{array}{cc}g_1\cdot\partial_xh_1\\g_2\cdot\partial_xh_2\end{array}\right))(x)=&-s_p(t)(\left(\begin{array}{cc}\partial_xg_1\cdot h_1\\\partial_xg_2\cdot h_2\end{array}\right))(x)+\partial_xs_p(t)(\left(\begin{array}{cc}g_1h_1\\g_2h_2\end{array}\right))(x)\\
&-\frac{\chi(t)}{2\pi}\int_{-\pi/T}^{\pi/T}\rho(\xi)e^{i\xi x+\la_c(\xi)t}\langle \mathcal{B}(\left(\begin{array}{cc}g_1h_1\\g_2h_2\end{array}\right))(\xi,\cdot), \partial_x\widetilde{\Phi}_\xi\rangle_{L^2} d\xi,
\end{align*}
and
\begin{align*}
\tilde{S}_c(t)(\left(\begin{array}{cc}g_1\cdot\partial_xh_1\\g_2\cdot\partial_xh_2\end{array}\right))(x)=&-\tilde{S}_c(t)(\left(\begin{array}{cc}\partial_xg_1\cdot h_1\\\partial_xg_2\cdot h_2\end{array}\right))(x)\\
&+\frac{1-\chi(t)}{2\pi}\int_{-\pi/T}^{\pi/T}\rho(\xi)e^{i\xi x+\la_c(\xi)t}i\xi\Phi_\xi(x)\langle \mathcal{B}(\left(\begin{array}{cc}g_1h_1\\g_2h_2\end{array}\right))(\xi,\cdot), \widetilde{\Phi}_\xi\rangle_{L^2} d\xi\\
&-\frac{1-\chi(t)}{2\pi}\int_{-\pi/T}^{\pi/T}\rho(\xi)e^{i\xi x+\la_c(\xi)t}\Phi_\xi(x)\langle \mathcal{B}(\left(\begin{array}{cc}g_1h_1\\g_2h_2\end{array}\right))(\xi,\cdot), \partial_x\widetilde{\Phi}_\xi\rangle_{L^2} d\xi\\
&+\frac{\chi(t)}{2\pi}\int_{-\pi/T}^{\pi/T}\rho(\xi)e^{i\xi x+\la_c(\xi)t}(i\xi)^2\left(\frac{\Phi_\xi(x)-\mathbf{\bar{U}}}{i\xi}\right)\langle \mathcal{B}(\left(\begin{array}{cc}g_1h_1\\g_2h_2\end{array}\right))(\xi,\cdot), \widetilde{\Phi}_\xi\rangle_{L^2} d\xi\\
&-\frac{\chi(t)}{2\pi}\int_{-\pi/T}^{\pi/T}\rho(\xi)e^{i\xi x+\la_c(\xi)t}i\xi\left(\frac{\Phi_\xi(x)-\mathbf{\bar{U}}}{i\xi}\right)\langle \mathcal{B}(\left(\begin{array}{cc}g_1h_1\\g_2h_2\end{array}\right))(\xi,\cdot), \partial_x\widetilde{\Phi}_\xi\rangle_{L^2} d\xi.
\end{align*}
Then, similar to the proof of Lemma \ref{crit}, the inequalities in (i) can be proved by the above two equalities and estimates in Lemma \ref{crit}.

Note that for $g\in H^2(\mathbb{R})$, we have
\begin{align}\label{sec}
\partial_x^2\mathcal{B}(g)(\xi,x)=\mathcal{B}(\partial_x^2 g)(\xi,x)-2i\xi\mathcal{B}(\partial_x g)(\xi,x)-\xi^2\check{g}(\xi,x).
\end{align}
Then (ii) can be proved similarly by combining \eqref{sec} and (i).
\end{proof}

\subsection{Nonlinear stability}\label{Secnon}
In this subsection, we devote to analyze nonlinear stability of periodic traveling wave $\mathbf{\bar{U}}=(\phi,\psi)^t$. First, we introduce the unmodulated perturbation $\mathbf{\widetilde{V}}=\mathbf{U}(t,x)-\mathbf{\bar{ U}}(x)$, which satisfies a semilinear equation and decays in $H^m(\mathbb{R}),m\in\mathbb{N}_0$ at most with rate $(1+t)^{-1/4}$. This decaying rate is not sufficient to close a nonlinear iteration scheme; see Remark \ref{unv}. Next, we introduce the spatio-temporal phase modulation $\varphi(t,x)$ and the modulated perturbation $\mathbf{V}(t,x)=\mathbf{U}(t,x-\varphi(t,x))-\mathbf{\bar{ U}}(x)$. The phase modulation $\varphi(t,x)$ is chosen to compensate the slowest decaying $(1+t)^{-1/4}$ such that $\mathbf{V}$ decays at rate $(1+t)^{-3/4}$. Then we can close a nonlinear iteration scheme by using $\mathbf{\widetilde{V}},\varphi,\mathbf{V}$ and their derivatives. Theorem \ref{nonsta} is the main result of this subsection.

\subsubsection{Unmodulated perturbation}
Let
$$\mathbf{\widetilde{V}}(t,x)=\mathbf{U}(t,x)-\mathbf{\bar{ U}}(x),
\quad \mathbf{\widetilde{V}}=(\tilde{u},\tilde{w})^t,$$ where $\mathbf{U}(t,x)$ is a solution of \eqref{FHN} and $\bar{\mathbf{U}}(x)$ is the periodic traveling wave. Then the unmodulated perturbation $\mathbf{\widetilde{V}}$ satisfies
\begin{equation}\label{unmod}
(\partial_t-\ml)\mathbf{\widetilde{V}}=\widetilde{\N}(\mathbf{\widetilde{V}}),
\end{equation}
where the nonlinear term $\widetilde{\N}(\mathbf{\widetilde{V}})=(-\tilde{u}^3+(a+1-3\phi)\tilde{u}^2,0)^t$.
Using the embedding inequality $\|g\|_{L^\infty(\mathbb{R})}\lesssim \|g\|_{H^1(\mathbb{R})}$, we can obtain the following estimates on $\widetilde{\N}(\mathbf{\widetilde{V}})$.

\begin{lemma}\label{tildeNV}
Suppose that $\tilde{u}\in H^4(\mathbb{R})$ with its $H^2$-norm having a fixed upper bound, that is there is a fixed constant $M>0$ such that $\|\tilde{u}\|_{H^2}\leq M$. Then
\begin{equation}\begin{split}
\|\widetilde{\N}(\mathbf{\widetilde{V}})\|_{L^1}\lesssim \|\tilde{u}\|_{L^2}^2,\quad \|\widetilde{\N}(\mathbf{\widetilde{V}})\|_{H^4}\lesssim \|\tilde{u}\|_{H^4}^2.
\end{split}\end{equation}
\end{lemma}

Note that the operator $\ml$ generates a $\mathcal{C}^0$-semigroup on $H^3(\mathbb{R})$, see Proposition \ref{group} (i). Then standard arguments imply the local existence and uniqueness of $\mathbf{\widetilde{V}}$ as a solution of \eqref{unmod}. See Proposition 4.3.9 of \cite{CH98} or Theorem 6.1.3 of \cite{pazy}.

\begin{prop}\label{tildeV}
Given $\mathbf{{V}}_0\in H^5(\mathbb{R})\times H^4(\mathbb{R})$, there exists a maximal time $T_{\rm{\rm{max}}}\in(0,\infty]$ such that \eqref{unmod} has a unique solution
$$\mathbf{\widetilde{V}}\in C([0,T_{\rm{max}}),H^5(\mathbb{R})\times H^4(\mathbb{R}))\cap {C}^1([0,T_{\rm{max}}),H^3(\mathbb{R})\times H^3(\mathbb{R})), $$
satisfying the initial value condition $\mathbf{\widetilde{V}}(0)=\mathbf{{V}}_0$. Furthermore, if $T_{\rm{max}}<\infty$, then
$$\lim_{t\uparrow T_{\rm{max}}}\|\mathbf{\widetilde{V}}(t)\|_{H^3}=\infty.$$
\end{prop}

The unmodulated perturbation $\mathbf{\tilde{V}}(t)$ has the following nonlinear damping estimates. That is, the higher derivatives of  $\mathbf{\widetilde{V}}(t)$ can be controlled by its $L^2$-norm.

\begin{prop}\label{damp}(nonlinear damping)
The unmodulated perturbation $\mathbf{\widetilde{V}}(t)$ satisfies the following inequality:
\begin{equation}\label{E}
\|\mathbf{\widetilde{V}}(t)\|^2_{H^4}\leq e^{-\e\gamma t}\|\mathbf{{V}}_0\|^2_{H^4}+C\int_0^te^{-\e\gamma(t-s)}\|\mathbf{\widetilde{V}}(s)\|^2_{L^2}ds,
\end{equation}
if $\|\tilde{u}(t)\|_{H^4}$ is sufficiently small.

\end{prop}
\begin{proof}
   Similar to the proof of Proposition \ref{group} (ii), we have
\begin{equation*}\begin{split}
\frac{1}{2}\frac{d}{dt}\|\mathbf{\widetilde{V}}(t)\|^2_{H^4}&=\Re \sum_{j=0}^4 (\langle {\partial}^j_x \tilde{u},{\partial}^j_x \tilde{u}_t\rangle_{L^2}+\langle {\partial}^j_x \tilde{w},{\partial}^j_x\tilde{w}_t\rangle_{L^2})\\
&\leq -\e\gamma\| \mathbf{\widetilde{V}}(t)\|_{H^4}^2+C_1\|\tilde{u}(t)\|_{L^2}^2+C_2\|\tilde{u}(t)\|_{H^4}\cdot\|\tilde{u}(t)\|_{H^4}^2\\
&\leq -\e\gamma\| \mathbf{\widetilde{V}}(t)\|_{H^4}^2+C_1\|\mathbf{\tilde{V}}(t)\|_{L^2}^2+C_2\|\tilde{u}(t)\|_{H^4}\cdot\|\mathbf{\tilde{V}}(t)\|_{H^4}^2,
\end{split}\end{equation*}by using Lemma \ref{tildeNV}. Then
provided by $\|\tilde{u}(t)\|_{H^4}\leq \frac{\e\gamma}{2C_2}$, it implies that
\begin{equation}
\frac{d}{dt}\|\mathbf{\widetilde{V}}(t)\|^2_{H^4}\leq -\e\gamma\|\mathbf{\widetilde{V}}(t)\|^2_{H^4}+C\|\mathbf{\widetilde{V}}(t)\|_{L^2}^2.
\end{equation}
By integrating both sides of the above inequality, the proposition is thus proved.
\end{proof}

\begin{remark}\label{unv}
By virtue of linear stability analysis in Section \ref{Seclin}, the unmodulated perturbation $\mathbf{\widetilde{V}}(t)$ decays in $H^1(\mathbb{R})$ at most with rate $(1+t)^{-\frac{1}{4}}$, which indicates that the decaying rate $(1+t)^{-1/4}$ is too weak to close a nonlinear iteration scheme. Intuitively, applying Duhamel's principle, we obtain the integral representation
\begin{equation}\label{tildeVr}
\mathbf{\widetilde{V}}(t)=e^{\ml t}\mathbf{{V}}_0+\int_0^te^{\ml (t-s)}\widetilde{\mathcal{N}}(\mathbf{\widetilde{V}}(s))ds.
\end{equation}
 Let $t>0$ and assume that $\|\mathbf{\widetilde{V}}(s)\|_{H^1}\lesssim (1+s)^{-1/4}$ for $s\in[0,t]$. Similar to Lemma \ref{tildeNV}, we can prove that
$\|\widetilde{\mathcal{N}}(\mathbf{\widetilde{V}}(s))\|_{H^1}\lesssim \|\mathbf{\widetilde{V}}\|^2_{H^1},$
where we use the fact that $\|\mathbf{\widetilde{V}}(s)\|_{H^1}$ is uniformly bounded for $s\in[0,t]$. Then, using Lemma \ref{exp}, Lemma \ref{crit}, Lemma \ref{tildeNV},
and expression \eqref{tildeVr}, we have
\begin{equation}\label{1/4}\begin{split}\|\mathbf{\widetilde{V}}(t)\|_{H^1}&\leq C(1+t)^{-\frac{1}{4}} \|\mathbf{V}_0\|_{H^1\cap L^1}+C\int_0^t (1+t-s)^{-\frac{1}{4}}\|\mathbf{\widetilde{V}}(s)\|^2_{H^1}ds\\
&\leq C(1+t)^{-\frac{1}{4}} \|\mathbf{V}_0\|_{H^1\cap L^1}+C\int_0^t (1+t-s)^{-\frac{1}{4}}(1+s)^{-\frac{1}{2}}ds\\
&\leq C(1+t)^{\frac{1}{4}}.
\end{split}\end{equation}
Therefore, we further consider the modulated perturbation in the next subsection.
\end{remark}

\subsubsection{Modulated perturbation}
Let
\begin{equation}
\label{moduv}\mathbf{V}(t,x)=\mathbf{U}(t,x-\varphi(t,x))-\mathbf{\bar{ U}}(x),
\quad \mathbf{V}=(u,w)^t,\end{equation}
where the spatio-temporal phase $\varphi(t)$ satisfies $\varphi(0)=0$. Substituting \eqref{moduv} into \eqref{402} yields the following quasilinear equation
\begin{align}\label{hatv}
(\partial_t-\ml)(\mathbf{V}+\varphi\bar{\mathbf{U}}')=\mathcal{N}(\mathbf{V},\varphi,\varphi_t)+(\partial_t-\ml)(\varphi_x\mathbf{V}),
\end{align}
where the nonlinear term is given by
$$\mathcal{N}(\mathbf{V},\varphi,\varphi_t)=\mathcal{Q}(\mathbf{V},\varphi)+\partial_x\mathcal{R}(\mathbf{V},\varphi,\varphi_t),$$
with
$$\mathcal{Q}(\mathbf{V},\varphi)=(1-\varphi_x)\left(\begin{array}{cc}-u^3+(a+1-3\phi)u^2\\0\end{array}\right),$$
and
$$\mathcal{R}(\mathbf{V},\varphi,\varphi_t)=
\left(\begin{array}{cc}-\varphi_t u+\varphi_{xx}u+2\varphi_x u_x+\frac{\varphi_x^2}{1-\varphi_x}(\phi'+u_x+c\varphi_x u)\\-\varphi_tw+c\varphi_xw\end{array}\right).$$
By using the embedding inequality $\|g\|_{L^\infty(\mathbb{R})}\lesssim \|g\|_{H^1(\mathbb{R})}$, a simple calculation shows that the nonlinear term $\mathcal{N}$ has the following estimates.

\begin{lemma}\label{modun}
Suppose that $\mathbf{V}=(u, w)^t\in H^2(\mathbb{R})\times H^1(\mathbb{R})$ and $(\varphi,\varphi_t)\in H^3(\mathbb{R})\times H^1(\mathbb{R})$ with the $H^1$-norm of $u$ having a upper bound and $\|\varphi\|_{H^2}\leq 1/2$. Then we have
\begin{equation}\begin{split}
\|{\N}(\mathbf{V},\varphi,\varphi_t)\|_{L^2}&\lesssim \|{u}\|_{L^2}\|{u}\|_{H^1}+\|(\varphi_x,\varphi_t)\|_{H^2\times H^1}(\|\mathbf{V}\|_{H^2\times H^1}+\|\varphi_x\|_{L^2}).
\end{split}\end{equation}
\end{lemma}

Using $\varphi(0)=0$ and applying Duhamel's principle to \eqref{hatv},  we obtain the following integral representation
\begin{align}\label{hatv1}
\mathbf{V}(t)+\varphi(t)\bar{\mathbf{U}}'=e^{\ml t}\mathbf{{V}}_0+\int_0^t e^{\ml (t-s)}\mathcal{N}(\mathbf{V}(s),\varphi(s),\varphi_t(s))ds+\varphi_x(t)\mathbf{V}(t),
\end{align}
with $\mathbf{V}(0)=\mathbf{{V}}_0$. Recall that
$$e^{\ml t}=\mathbf{\bar{U}}'(x)s_p(t)+\tilde{S}(t),$$
with
$$\tilde{S}(t)=\tilde{S}_c(t)+S_e(t),$$
and $\mathbf{\bar{U}}'(x)s_p(t)$ and $\tilde{S}(t)$ result in the slowest decay rate $(1+t)^{-1/4}$ and faster decay rate $(1+t)^{-3/4}$, respectively. Define $\varphi(t)$ as
\begin{align}\label{varphi}
\varphi(t) := s_p(t)\mathbf{V}_0+\int_0^t s_p(t-s)\mathcal{N}(\mathbf{V}(s),\varphi(s),\varphi_t(s))ds.
\end{align}
Note that the modulated perturbation $\mathbf{V}$ can be represented by
\begin{align}\label{VV}
\mathbf{V}(t,x)=\mathbf{\widetilde{V}}(t,x-\varphi(t,x))+\mathbf{\bar{U}}(x-\varphi(t,x))-\mathbf{\bar{U}}(x).
\end{align}
Then, for a given $\mathbf{\widetilde{V}}$, \eqref{varphi} is an integral equation for $\varphi$. By using Banach fixed point theorem, we prove the following proposition in Appendix.

\begin{prop}\label{exivar}
For $\mathbf{\widetilde{V}}$ and $T_{\rm{max}}$ given in Proposition \ref{tildeV}, there exists a maximal time $\tau_{\rm{max}}\leq T_{\rm{max}}$ such that the integral equation \eqref{varphi} with $\mathbf{V}$ given by \eqref{VV} has a unique solution
$$\varphi\in C([0,\tau_{\rm{max}}),H^3(\mathbb{R}))\cap C^1([0,\tau_{\rm{max}}),H^1(\mathbb{R})),$$
with $\varphi(t)=0$, for $t\in[0,1]$. Moreover, if $\tau_{\rm{max}}<T_{\rm{max}}$, then
$$\lim_{t\uparrow\tau_{\rm{max}}}\|(\varphi(t),\partial_t\varphi(t))\|_{H^3\times H^1}=\infty.$$
\end{prop}
Combining Lemma \ref{V} and Proposition \ref{Vphi} in the appendix, we obtain the following result.

\begin{cor}
 For $\mathbf{\widetilde{V}}$ as in Proposition \ref{tildeV} and $\varphi$ and $\tau_{\rm{max}}$ given by Proposition \ref{exivar}, the modulated perturbation $\mathbf{V}$ defined by \eqref{VV} satisfies
 $$\mathbf{V}\in C([0,\tau_{\rm{max}}),H^2(\mathbb{R})\times H^2(\mathbb{R})).$$
 Moreover,  $\mathbf{V}$ is a solution of the integral equation \eqref{hatv1}.
\end{cor}

Then, substituting \eqref{varphi} into \eqref{hatv1}, we obtain
\begin{align}\label{hatv2}
\mathbf{V}(t)=\tilde{S}(t)\mathbf{{V}}_0+\int_0^t \tilde{S}(t-s) \mathcal{N}(\mathbf{V}(s),\varphi(s),\varphi_t(s))ds+\varphi_x(t)\mathbf{V}(t),\quad \text{ for } t\in[0,\tau_{\rm{max}}).
\end{align}
Based on Proposition \ref{exivar} and $s_p(t)=0$ for $t\in[0,1]$, we obtain
\begin{align}\label{varphix}
\partial_x^l\partial_t^j\varphi(t)=\partial_x^l\partial_t^js_p(t)\mathbf{{V}}_0+\int_0^t \partial_x^l\partial_t^js_p(t-s)\mathcal{N}(\mathbf{V}(s),\varphi(s),\varphi_t(s))ds,
\end{align}
for $l,j\in\mathbb{N}_0$ with $l+2j\leq 3$ and $t\in[0,\tau_{\rm{max}})$. 

Next, we give some estimates for nonlinear terms $\mathcal{N}(\mathbf{V},\varphi,\varphi_t)=\mathcal{Q}(\mathbf{V},\varphi)+\partial_x\mathcal{R}(\mathbf{V},\varphi,\varphi_t)$ in \eqref{hatv}.

\begin{lemma}\label{Inter}
Suppose that $\mathbf{V}=(u,w)^t\in H^2(\mathbb{R})\times H^1(\mathbb{R})$ and $(\varphi,\varphi_t)\in H^3(\mathbb{R})\times H^1(\mathbb{R})$ with the $H^1$-norm of 
 $u$ having a upper bound and $\|\varphi\|_{H^3}\leq 1/2$. Then the following inequalities
\begin{equation}\begin{split}
\|\mathcal{Q}(\mathbf{V},\varphi)\|_{L^1}&\lesssim \|{u}\|_{L^2}^2,\\
\|\partial_x^l\partial_t^js_p(t)(\partial_x\mathcal{R}(\mathbf{V},\varphi,\varphi_t))\|_{L^2}&\lesssim (1+t)^{-\frac{1}{4}-\frac{l+j}{2}}\|(\varphi_x,\varphi_t)\|_{H^2\times H^1}(\|\mathbf{V}\|_{L^2}+\|\varphi_x\|_{L^2}),\\
\|\tilde{S}_c(t)(\partial_x\mathcal{R}(\mathbf{V},\varphi,\varphi_t))\|_{L^2}&\lesssim (1+t)^{-\frac{3}{4}}\|(\varphi_x,\varphi_t)\|_{H^2\times H^1}(\|\mathbf{V}\|_{L^2}+\|\varphi_x\|_{L^2}),
\end{split}\end{equation}
hold for all integers $l,j$ with $0\leq l+2j\leq3$.
\end{lemma}
\begin{proof}
The first inequality is obtained directly  by using the embedding inequality $\|g\|_{H^1}\lesssim \|g\|_{L^\infty}$ and $\|u\|_{H^1}\leq C$. Note that
$$\partial_x\mathcal{R}(\mathbf{V},\varphi,\varphi_t)=\mathcal{R}_1(\varphi,\varphi_t)\mathbf{V}_{xx}+\mathcal{R}_2(\varphi,\varphi_t)\mathbf{V}_{x}+\mathcal{R}_3(\varphi,\varphi_t)\mathbf{V}+\mathcal{R}_4(\varphi).$$
A simple calculation shows that
$$\|\mathcal{R}_k(\varphi,\varphi_t)\|_{H^{3-k}},\|\mathcal{R}_4(\varphi)\|_{L^2}\lesssim \|(\varphi,\varphi_t)\|_{H^2\times H^1},\text {  for  }k=1,2,3,$$
by using $\|\varphi\|_{H^3}\leq 1/2$.
Then the last two inequalities follow from the above estimates, Lemma \ref{crit} and Proposition \ref{part}.
\end{proof}

\subsubsection{Proof of Theorem \ref{nonsta}}
Recall 
\begin{equation*}\begin{split}
&\mathbf{\widetilde{V}}(t)=e^{\ml t}\mathbf{{V}}_0+\int_0^t e^{\ml (t-s)}\widetilde{\mathcal{N}}(\mathbf{\widetilde{V}}(s))ds,\\
&\varphi(t) := s_p(t)\mathbf{V}_0+\int_0^t s_p(t-s)\mathcal{N}(\mathbf{V}(s),\varphi(s),\varphi_t(s))ds,\\
&\mathbf{V}(t)=\tilde{S}(t)\mathbf{{V}}_0+\int_0^t \tilde{S}(t-s) \mathcal{N}(\mathbf{V}(s),\varphi(s),\varphi_t(s))ds+\varphi_x(t)\mathbf{V}(t),
\end{split}\end{equation*}
where $e^{\ml t}$ and $ \tilde{S}(t)=\tilde{S}_c(t)+S_e(t)$ decay with rate $(1+t)^{-1/4}$ and $ (1+t)^{-3/4}$ respectively. Thus, we define a function
\begin{align*}\eta(t):=\sup_{0\leq s\leq t}&[(1+s)^{\frac{3}{4}}(\|\mathbf{V}(s)\|_{L^2}+\|\varphi_x(s)\|_{H^2}+\|\varphi_t(s)\|_{H^1})\\
&+(1+s)^{\frac{1}{4}}(\|\widetilde{\mathbf{V}}(s)\|_{H^4}+\|\partial_x\mathbf{V}(s)\|_{H^1}+\|\varphi(s)\|_{L^2})],
\end{align*}
which is positive, continuous and monotonically increasing. Furthermore, if $\tau_{\rm{max}}<\infty$, then we have
\begin{align}\label{etain}
\lim_{t\uparrow \tau_{\rm{max}}}\eta(t)=\infty.
\end{align}
Referring to many studies for nonlinear stability \cite{HJPd23,JZ11,de24,deS21}, one often needs to prove the following claim.

{\bf Claim.} There exist constants $\mathbf{K}>0$ and $C>1$ such that for any $t\in[0,\tau_{\rm{max}})$ with $\eta(t)\leq \mathbf{K}$, the inequality
\begin{align}\label{eta2}
\eta(t)\leq C(E_0+\eta(t)^2),
\end{align}
holds, where $E_0>0$ is defined in Theorem \ref{nonsta}.

If this claim is true, then it implies that $\tau_{\rm{max}}=\infty$ and $\eta(t)$ is bounded for all $t\geq0$. In fact, by virtue of \eqref{eta2}, we define a function
$$q(z):=Cz^2-z+CE_0,$$
which has two positive zeroes
$$z_1=\frac{1-\sqrt{1-4C^2E_0}}{2C}< 2CE_0,\quad z_2=\frac{1+\sqrt{1-4C^2E_0}}{2C}> 2CE_0,$$
if $E_0<\frac{1}{4C^2}$. Furthermore, $q(z)>0$ for $0<z<z_1$ and $z>z_2$. By the definition of $\eta$, we have $\eta(0)\leq 2CE_0<\mathbf{K}$, provided by $E_0<\frac{\mathbf{K}}{2C}$.
Then we deduce
$$\eta(t)\leq z_1 < 2CE_0<\mathbf{K},\text{ for all }t\in[0,\tau_{\rm{max}}),$$
by continuity of $\eta$,
which contradicts with \eqref{etain}. Then Theorem \ref{nonsta} follows by taking $\e_0=\min\{\frac{1}{4C^2},\frac{\mathbf{K}}{2C}\}$ and $M=2C$.

Next, we devote to prove this claim. Let $\mathbf{K}=\min\{\frac{\e\gamma}{2C_2},\frac{1}{2}\}$ and assume $t\in[0,\tau_{\rm{max}})$ such that $\eta(t)\leq \mathbf{K}$.
Using the estimate in Lemma \ref{modun} and definition of $\eta$, we obtain
$$\|\mathcal{N}(\mathbf{V}(s),\varphi(s),\varphi_t(s))\|_{L^2}\lesssim \eta(s)^2(1+s)^{-1}.$$
Then, by virtue of Lemma \ref{exp}, we have
\begin{equation}\begin{split}
\left\|\int_0^tS_e(t-s)\mathcal{N}(\mathbf{V}(s),\varphi(s),\varphi_t(s))ds\right\|_{L^2}\lesssim \int_0^t\frac{\eta(s)^2e^{-\mu(t-s)}}{1+s}ds\lesssim \frac{\eta(t)^2}{1+t}.
\end{split}\end{equation}
Additionally, for $l,j\in\mathbb{N}_0$ with $l+2j\leq3$ and $s\in[0,t]$, we have
\begin{align*}
\|\mathcal{Q}(\mathbf{V}(s),\varphi(s))\|_{L^1}&\lesssim \eta(s)^2(1+s)^{-\frac{3}{2}},\\
\|\partial_x^l\partial_t^js_p(t-s)(\partial_x\mathcal{R}(\mathbf{V}(s),\varphi(s),\varphi_t(s)))\|_{L^2}&\lesssim \eta(s)^2(1+t-s)^{-\frac{1}{4}-\frac{l+j}{2}}(1+s)^{-\frac{3}{2}},\\
\|\tilde{S}_c(t-s)(\partial_x\mathcal{R}(\mathbf{V}(s),\varphi(s),\varphi_t(s)))\|_{L^2}&\lesssim \eta(s)^2(1+t-s)^{-\frac{3}{4}}(1+s)^{-\frac{3}{2}},
\end{align*}
by using Lemma \ref{Inter} and $\eta(t)\leq \mathbf{K}$. Furthermore, we obtain
\begin{equation}
\begin{split}
\left\|\int_0^ts_p(t-s)\mathcal{N}(\mathbf{V}(s),\varphi(s),\varphi_t(s))ds\right\|_{L^2} \leq&\left\|\int_0^ts_p(t-s)\mathcal{Q}(\mathbf{V}(s),\varphi(s))ds\right\|_{L^2}\\
&+\left\|\int_0^ts_p(t-s)\partial_x\mathcal{R}(\mathbf{V}(s),\varphi(s),\varphi_t(s))ds\right\|_{L^2}\\
\lesssim&\int_0^t\frac{\eta(s)^2}{(1+t-s)^{\frac{1}{4}}(1+s)^{\frac{3}{2}}}ds\\
\lesssim&\frac{\eta(t)^2}{(1+t)^\frac{1}{4}},\\
\left\|\int_0^t\partial_x^l\partial_t^js_p(t-s)\mathcal{N}(\mathbf{V}(s),\varphi(s),\varphi_t(s))ds\right\|_{L^2} \lesssim&\int_0^t\frac{\eta(s)^2}{(1+t-s)^{\frac{3}{4}}(1+s)^{\frac{3}{2}}}ds\\
\lesssim&\frac{\eta(t)^2}{(1+t)^\frac{3}{4}},\quad \text { for }1\leq l+2j\leq3,\\
\left\|\int_0^t\tilde{S}_c(t-s)\mathcal{N}(\mathbf{V}(s),\varphi(s),\varphi_t(s))ds\right\|_{L^2} \lesssim&\int_0^t\frac{\eta(s)^2}{(1+t-s)^{\frac{3}{4}}(1+s)^{\frac{3}{2}}}ds\\
\lesssim&\frac{\eta(t)^2}{(1+t)^\frac{3}{4}}.
\end{split}
\end{equation}
Then using expressions \eqref{varphi}, \eqref{hatv2} and \eqref{varphix}, we arrive at
\begin{equation}\label{varv}\begin{split}
\|\varphi(t)\|_{L^2}\lesssim & \frac{E_0+\eta(t)^2}{(1+t)^{\frac{1}{4}}},\\
\|\mathbf{V}(t)\|_{L^2},\|\partial_x\varphi(t)\|_{H^2},\|\partial_t\varphi(t)\|_{H^1}\lesssim &\frac{E_0+\eta(t)^2}{(1+t)^{\frac{3}{4}}}.
\end{split}\end{equation}

Next, we turn to control $H^4$-norms of $\mathbf{\widetilde{V}}$. Note that
$$\mathbf{V}(t,x)-\mathbf{\widetilde{V}}(t,x)=\mathbf{U}(t,x-\varphi(t,x))-\mathbf{U}(t,x),$$
which implies
\begin{equation}|\mathbf{V}(t,x)-\mathbf{\widetilde{V}}(t,x)|\leq \|\mathbf{U}_x(t)\|_{L^\infty}|\varphi(t,x)|\leq (\|\mathbf{\bar{U}}'\|_{L^\infty}+\|\mathbf{\widetilde{V}}(t)\|_{H^2})|\varphi(t,x)|,
\end{equation}by using the mean value theorem.
Furthermore, we establish that
\begin{align}\label{diffe}
\|\mathbf{V}(t)-\mathbf{\widetilde{V}}(t)\|_{L^2}&\leq (\|\mathbf{\bar{U}}'\|_{L^\infty}+\|\mathbf{\widetilde{V}}(t)\|_{H^2})\|\varphi(t)\|_{L^2},
\end{align}
Then, combining \eqref{varv} and \eqref{diffe} yields
\begin{align*}
\|\mathbf{\widetilde{V}}(t)\|_{L^2} \lesssim \frac{E_0+\eta(t)^2}{(1+t)^{\frac{1}{4}}},
\end{align*}
where we use $\eta(t)\leq \mathbf{K}$.
Furthermore, by using nonlinear damping estimates, i.e., Proposition \ref{damp},
we have
\begin{equation}\begin{split}
\|\mathbf{\widetilde{V}}(t)\|_{H^4}^2  & \leq e^{-\e\gamma t}\|\mathbf{{V}}_0\|^2_{H^4}+C\int_0^t e^{-\e\gamma(t-s)}\|\mathbf{\widetilde{V}}(s)\|^2_{L^2}ds\\
&\lesssim e^{-\e\gamma t} E_0^2+C\int_0^te^{-\e\gamma(t-s)}\frac{(E_0+\eta(s)^2)^2}{(1+s)^{\frac{1}{2}}}ds\\
&\lesssim \left(e^{-\e\gamma t}+\int_0^te^{-\e\gamma(t-s)}\frac{1}{(1+s)^{\frac{1}{2}}}ds\right)(E_0+\eta(t)^2)^2\\
&\lesssim
\frac{(E_0+\eta(t)^2)^2}{(1+t)^{\frac{1}{2}}}.
\end{split}\end{equation}

Last, we estimate $\|\partial_x \mathbf{V}\|_{H^1}$. Similar to \eqref{diffe}, we have
\begin{equation}\begin{split}
\|\mathbf{V}_x(t)-\mathbf{\widetilde{V}}_x(t)\|_{L^2}=&\|\mathbf{U}_x(\cdot-\varphi(t,\cdot),t)-\mathbf{U}_x(t,\cdot)-\mathbf{U}_x(\cdot-\varphi(t,\cdot),t)\varphi_x(t,\cdot)\|_{L^2}\\
\leq & (\|\mathbf{\bar{U}}'\|_{L^\infty}+\|\mathbf{\widetilde{V}}(t)\|_{H^2})\|\varphi_x(t)\|_{L^2}\\
&+(\|\mathbf{\bar{U}}''\|_{L^\infty}+\|\mathbf{\widetilde{V}}(t)\|_{H^3})\|\varphi(t)\|_{L^2},\\
\|\mathbf{V}_{xx}(t)-\mathbf{\widetilde{V}}_{xx}(t)\|_{L^2}=& \|-\mathbf{U}_x(\cdot-\varphi(t,\cdot),t)\varphi_{xx}(t,\cdot)-2\mathbf{U}_{xx}(\cdot-\varphi(t,\cdot),t)\varphi_x(t,\cdot)\\
 &+\mathbf{U}_{xx}(\cdot-\varphi(t,\cdot),t)\varphi_x^2(t,\cdot)+\mathbf{U}_{xx}(\cdot-\varphi(t,\cdot),t)-\mathbf{U}_{xx}(t,\cdot)\|_{L^2}\\
\leq & (\|\mathbf{\bar{U}}'\|_{L^\infty}+\|\mathbf{\widetilde{V}}(t)\|_{H^2})\|\varphi_{xx}(t)\|_{L^2}\\
&+(\|\mathbf{\bar{U}}'''\|_{L^\infty}+\|\mathbf{\widetilde{V}}(t)\|_{H^4})\|\varphi(t)\|_{L^2}\\
&+(\|\mathbf{\bar{U}}''\|_{L^\infty}+\|\mathbf{\widetilde{V}}(t)\|_{H^3})\|\varphi_x(t)\|_{L^2}(2+\|\varphi(t)\|_{H^2}).
\end{split}\end{equation}
As a consequence, we arrive at
\begin{equation}
\begin{split}
\|\mathbf{V}_x(t)\|_{L^2}&\leq \|\mathbf{\widetilde{V}}_x(t)\|_{L^2}+(\|\mathbf{\bar{U}}\|_{W^{2,\infty}}+\|\mathbf{\widetilde{V}}(t)\|_{H^3})\|\varphi(t)\|_{H^1}\lesssim \frac{E_0+\eta(t)^2}{(1+t)^{\frac{1}{4}}},\\
\|\mathbf{V}_{xx}(t)\|_{L^2}& \lesssim \|\mathbf{\widetilde{V}}_{xx}(t)\|_{L^2}+(\|\mathbf{\bar{U}}\|_{W^{3,\infty}}+\|\mathbf{\widetilde{V}}(t)\|_{H^4})\|\varphi(t)\|_{H^2}\lesssim \frac{E_0+\eta(t)^2}{(1+t)^{\frac{1}{4}}},
\end{split}
\end{equation}
where we use $\eta(t)\leq \mathbf{K}$. Thus, we complete the proof of the claim.


\bigskip

\appendix
\section{A general boundary-value problem}\label{A}
\renewcommand{\theequation}{A.\arabic{equation}}
\setcounter{equation}{0}

In this subsection, we first solve an inhomogeneous boundary value problem
\begin{equation}\label{213}
\begin{aligned}
(\romannumeral 1) \quad & w_{1,\pm}'(x) = F_u(h_1(x),  {\bf 0}) w_{1,\pm}(x) + G_{1,\pm}(x),  \quad x \in (-L_1, L_1),\\
(\romannumeral 2) \quad & w_{2,\pm}'(x) = F_u(h_2(x-L_1-L_2),  {\bf 0}) w_{2,\pm}(x) + G_{2,\pm}(x),  \quad  x \in (L_1, L_1+2L_2),\\
(\romannumeral 3) \quad & Q_1w_{1,\pm}(0) = {\bf 0}, \\
(\romannumeral 4) \quad & Q_2w_{2,\pm}(L_1+L_2) =  {\bf 0}, \\
(\romannumeral 5) \quad & w_{1,+}(0) - w_{1,-}(0) \in Y_1^\bot, \\
(\romannumeral 6) \quad & w_{2,+}(L_1+L_2) - w_{2,-}(L_1+L_2) \in Y_2^\bot, \\
(\romannumeral 7) \quad & w_{1,+}(L_1) - w_{2,-}(L_1)= D_1, \\
(\romannumeral 8) \quad & w_{2,+}(L_1+2L_2) - e^{i\xi T}w_{1,-}(-L_1) = D_2,
\end{aligned}
\end{equation}
for elements \( G = (G_{1,-}, G_{1,+},G_{2,-},G_{2,+}) \in V_w \) and $ D =(D_1, D_2) \in \mathbb{C}^{2n}$.

 
We seek solutions of equation (\ref{213}) for $\xi \in [-\pi/T, \pi/T)$. To solve (\ref{213}), we first introduce a variation of constants formula that can capture all solutions of (\ref{213}) (\romannumeral 1)-(\romannumeral 4). Then, we solve (\romannumeral 5)-(\romannumeral 8) in Lemmas \ref{lem2.2} and \ref{lem2.3}. Finally, we apply this formula in Lemma \ref{lem2.4} to derive a reduced equation involving the period $T := 2(L_1 + L_2)$ and the right-hand side function $G$.

Using the variation of constants formula, the general form of the solution satisfying the inhomogeneous equation (\ref{213})(\romannumeral 1)-(\romannumeral 2) takes the form of:
\begin{equation}\label{214}
\begin{aligned}
(\romannumeral 1) \quad	w_{1,-}(x) =& \Phi^s_{1,-}(x, -L_1) w_{1,-}(-L_1) + \Phi^u_{1,-}(x, 0) w_{1,-}(0) \\
&+ \int_{0}^{x} \Phi^u_{1,-}(x, y) G_{1,-}(y) dy + \int_{-L_1}^{x} \Phi^s_{1,-}(x, y) G_{1,-}(y) dy,\\
(\romannumeral 2) \quad	w_{1,+}(x) =& \Phi^u_{1,+}(x, L_1) w_{1,+}(L_1) + \Phi^s_{1,+}(x, 0) w_{1,+}(0) \\
&+ \int_{0}^{x} \Phi^s_{1,+}(x, y) G_{1,+}(y) dy + \int_{L_1}^{x} \Phi^u_{1,+}(x, y) G_{1,+}(y) dy,\\
(\romannumeral 3) \quad	w_{2,-}(x) =& \Phi^s_{2,-}(x-L_1-L_2, -L_2) w_{2,-}(-L_2) + \Phi^u_{2,-}(x-L_1-L_2, 0) w_{2,-}(0) \\
&+ \int_{L_1+L_2}^{x} \Phi^u_{2,-}(x-L_1-L_2, y-L_1-L_2) G_{2,-}(y) dy  \\
&+ \int_{L_1}^{x} \Phi^s_{2,-}(x-L_1-L_2, y-L_1-L_2) G_{2,-}(y) dy, \\
(\romannumeral 4) \quad	w_{2,+}(x) =& \Phi^u_{2,+}(x-L_1-L_2, L_2) w_{2,+}(L_2) + \Phi^s_{2,+}(x-L_1-L_2, 0) w_{2,+} (0) \\
&+ \int_{L_1+L_2}^{x} \Phi^s_{2,+}(x-L_1-L_2, y-L_1-L_2) G_{2,+}(y) dy  \\
 &+ \int_{L_1+2L_2}^{x} \Phi^u_{2,+}(x-L_1-L_2, y-L_1-L_2) G_{2,+}(y) dy.
\end{aligned}
\end{equation}

Define the spaces over the complex field
\begin{align*}
V_w &:= C^0([-L_1, 0], \mathbb{C}^n) \times C^0([0, L_1], \mathbb{C}^n) \times C^0([L_1, L_1+L_2], \mathbb{C}^n) \times C^0([L_1+L_2, L_1+2L_2], \mathbb{C}^n),\\
V_a &:= R(P_1^s) \times R(P_2^u) \times R(P_2^s) \times R(P_1^u),\\
V_b &:= Y_1^u \times Y_1^s \times Y_2^u \times Y_2^s,
\end{align*}
 and let
\begin{align*}
w &= (w_{1,-}, w_{1,+}, w_{2,-}, w_{2,+}) \in V_w,\\
a &= (a_{1,-}, a_{1,+}, a_{2,-}, a_{2,+}) \in V_a,\\
b &= (b_{1,-}, b_{1,+}, b_{2,-}, b_{2,+}) \in V_b.
\end{align*}

Due  to the estimation
$$|P^s_{1,-}(x) - P^s_1| \leq Ce^{-\alpha_i^u |x|},  $$
we can replace $w_{1,-}(-L_1)$ with $a_{1,-}$. This ensures that $a_{1,-}$ is independent of $L_1$.  Similarly, we replace $w_{1,+}(L_1)$, $w_{2,-}(-L_2)$ and $w_{2,+}(L_2)$ in (\ref{214})  with $a_{1,+}$, $a_{2,-}$ and $a_{2,+}$, respectively.
Based on the definition of the projection \( Q_i \) and the definitions of \( Y_i^u \) and \( Y_i^s \) in (\ref{206}), it is possible to select \( (w_{1,-}(0), w_{1,+}(0), w_{2,-}(0), w_{2,+}(0)) \in V_b \) such that \eqref{214} becomes a solution to \eqref{213}(i)-(iv). Therefore, the general form of the solution to (\ref{213})(\romannumeral 1)-(\romannumeral 4) can be expressed as: 
\begin{equation}\label{215}
	\begin{aligned}
	(\romannumeral 1) \quad	w_{1,-}(x) =& \Phi^s_{1,-}(x, -L_1) a_{1,-} + \Phi^u_{1,-}(x, 0) b_{1,-}  \\
		&+ \int_{0}^{x} \Phi^u_{1,-}(x, y) G_{1,-}(y) dy + \int_{-L_1}^{x} \Phi^s_{1,-}(x, y) G_{1,-}(y) dy,\\
	(\romannumeral 2) \quad	w_{1,+}(x) =& \Phi^u_{1,+}(x, L_1) a_{1,+} + \Phi^s_{1,+}(x, 0) b_{1,+} \\
		&+ \int_{0}^{x} \Phi^s_{1,+}(x, y) G_{1,+}(y) dy + \int_{L_1}^{x} \Phi^u_{1,+}(x, y) G_{1,+}(y) dy,\\
	(\romannumeral 3) \quad	w_{2,-}(x) =& \Phi^s_{2,-}(x-L_1-L_2, -L_2) a_{2,-} + \Phi^u_{2,-}(x-L_1-L_2, 0) b_{2,-} \\
		&+ \int_{L_1+L_2}^{x} \Phi^u_{2,-}(x-L_1-L_2, y-L_1-L_2) G_{2,-}(y) dy  \\
     &+ \int_{L_1}^{x} \Phi^s_{2,-}(x-L_1-L_2, y-L_1-L_2) G_{2,-}(y) dy, \\
	(\romannumeral 4) \quad	w_{2,+}(x) =& \Phi^u_{2,+}(x-L_1-L_2, L_2) a_{2,+} + \Phi^s_{2,+}(x-L_1-L_2, 0) b_{2,+} \\
		&+ \int_{L_1+L_2}^{x} \Phi^s_{2,+}(x-L_1-L_2, y-L_1-L_2) G_{2,+}(y) dy  \\
     &+ \int_{L_1+2L_2}^{x} \Phi^u_{2,+}(x-L_1-L_2, y-L_1-L_2) G_{2,+}(y) dy,
	\end{aligned}
\end{equation}
where the elements \( a \in V_a \) and \( b \in V_b \) are chosen arbitrarily.

\begin{lemma}\label{lem2.2}
There exist constants \( C \) and \( L_* \) such that the following statement holds for each \( L_1 + L_2 > L_*\):\\
(1) The right-hand side of (\ref{215}) defines a linear operator
\[ W_1 : V_a \times V_b \times V_w \to V_w, \quad (a, b, G) \mapsto W_1(a, b, G).\]
(2) There is a linear operator \( B_1 : V_a \times V_w \to V_b \) such that \( w \) satisfies (\ref{213})(\romannumeral 1)-(\romannumeral 6) if and only if,
\begin{equation} \label{216}
w = W_1(a, B_1(a, G), G),
\end{equation}
which implies \( b = B_1(a, G) \). \\
(3) The estimates are given as follows:
\begin{equation}\label{217}
\begin{aligned}
|B_1(a, G)| &\leq C(e^{-\alpha L} |a| + |G|), \\
|W_1(a, b, G)| &\leq C(|a| + |b| + |G|), \\
|W_1(a, B_1(a, G), G)|& \leq C(|a| + |G|), \\
\end{aligned}
\end{equation}
where \(L = \min\{L_1, L_2\}\).
\end{lemma}

\begin{proof}

The results of \( W_1 \) follow from the exponential decay estimates for the evolution operators, as established in Lemma \ref{lem2.1}. It is now necessary to show that (\ref{213})(\romannumeral 5)-(\romannumeral 6) hold with an appropriate choice of \( b \). By substituting \( x = 0 \) into (\ref{215})(\romannumeral 1), we obtain 
\begin{equation}\label{218}
\begin{aligned}
w_{1,+} (0) - w_{1,-} (0) =& b_{1,+} - b_{1,-} + \Phi^u_{1,+}(0, L_1) a_{1,+} - \Phi^s_{1,-}(0, -L_1) a_{1,-} \\
&- \int_0^{L_1} \Phi^u_{1,+}(0, y) G_{1,+}(y) dy - \int_{-L_1}^0 \Phi^s_{1,-}(0, y) G_{1,-}(y) dy,
\end{aligned}
\end{equation}
since \( \Phi^s_{1,+}(0, 0) b_{1,+} = b_{1,+} \) and \( \Phi^u_{1,-}(0, 0) b_{1,-} = b_{1,-} \), which are derived from (\ref{206}) and \( b \in V_b \).

Then, substituting \( x = L_1+L_2 \) into (\ref{215})(\romannumeral 2), we also get
\begin{equation}\label{219}
\begin{aligned}
w_{2,+} (L_1+L_2) - w_{2,-} (L_1+L_2) =& b_{2,+} - b_{2,-} + \Phi^u_{2,+}(0, L_2) a_{2,+} - \Phi^s_{2,-}(0, -L_2) a_{2,-} \\
&- \int_{L_1+L_2}^{L_1+2L_2} \Phi^u_{2,+}(0, y-L_1-L_2) G_{2,+}(y) dy \\
&- \int_{L_1}^{L_1+L_2} \Phi^s_{2,-}(0, y-L_1-L_2) G_{2,-}(y) dy,
\end{aligned}
\end{equation}
since \( \Phi^s_{2,+}(0, 0) b_{2,+} = b_{2,+} \) and \( \Phi^u_{2,-}(0, 0) b_{2,-} = b_{2,-} \), which are derived from (\ref{206}) and \( b \in V_b \).

By (\ref{213})(\romannumeral 3)(\romannumeral 4), we have that
\begin{align*}
    &P(Y_1^c, Y_1^s \times Y_1^u \times Y_1^\bot) w_{1,\pm} (0) = 0,\\
     &P(Y_2^c, Y_2^s \times Y_2^u \times Y_2^\bot) w_{2,\pm} (L_1+L_2) = 0,
\end{align*}
where \( P(X, Y)\) represents a projection on \( X \) along \( Y\).

To solve (\ref{213})(\romannumeral 5)(\romannumeral 6), it is only necessary to require
\begin{align*}
    &P(Y_1^u, Y_1^s \times Y_1^c \times Y_1^\bot) (w_{1,+}-w_{1,-}) (0) = 0,\\
    &P(Y_1^s, Y_1^c \times Y_1^u \times Y_1^\bot) (w_{1,+}-w_{1,-}) (0) = 0,\\
     &P(Y_2^u, Y_2^s \times Y_2^c \times Y_2^\bot) (w_{2,+}-w_{2,-}) (L_1+L_2) = 0,\\
     &P(Y_2^s, Y_2^c \times Y_2^u \times Y_2^\bot) (w_{2,+}-w_{2,-}) (L_1+L_2) = 0.
\end{align*}

Thus, by projecting (\ref{218}) onto \( Y_1^s \times Y_1^u \) and using the definition of \( Y_1^s \) and \( Y_1^u \) from (\ref{206}), we see that \( w_{1,+} (0) - w_{1,-}(0) \in Y_1^\bot \) if and only if,
\begin{align*}
 b_{1,-} &= P(Y_1^u, Y_1^s \times Y_1^c \times Y_1^\bot) \left( \Phi^u_{1,+}(0, L_1) a_{1,+}  - \int_0^{L_1} \Phi^u_{1,+}(0, y) G_{1,+}(y) dy \right), \\
b_{1,+} &= P(Y_1^s, Y_1^c \times Y_1^u \times Y_1^\bot) \left( \Phi^s_{1,-}(0, -L_1) a_{1,-} + \int_{-L_1}^0 \Phi^s_{1,-}(0, y) G_{1,-}(y) dy\right).
\end{align*}

By projecting (\ref{219}) onto \( Y_2^s \times Y_2^u \) and using the definition of \( Y_2^s \) and \( Y_2^u \) from (\ref{206}), we see that
$$ w_{2,+} (L_1 + L_2) - w_{2,-}(L_1 + L_2) \in Y_2^\bot, $$
if and only if,
\begin{align*}
 b_{2,-} &= P(Y_2^u, Y_2^s \times Y_2^c \times Y_2^\bot) \left( \Phi^u_{2,+}(0, L_2) a_{2,+}  - \int_{L_1+L_2}^{L_1+2L_2} \Phi^u_{2,+}(0, y-L_1-L_2) G_{2,+}(y) dy \right), \\
b_{2,+} &= P(Y_2^s, Y_2^c \times Y_2^u \times Y_2^\bot) \left( \Phi^s_{2,-}(0, -L_2) a_{2,-} + \int_{L_1}^{L_1+L_2} \Phi^s_{2,-}(0, y-L_1-L_2) G_{2,-}(y) dy\right).
\end{align*}
The bounded linear operator \( B_1 : V_a \times V_w \to V_b \) is defined according to the right-hand sides of these equations, which satisfies the inequality
\[|B_1(a, G)| \leq C(e^{-\alpha L} |a| + |G|),\]
based on the exponential dichotomies in Lemma \ref{lem2.1}. This completes the proof of the lemma. 
\end{proof}

In the following lemma, we shall solve (\ref{213})(\romannumeral 7)\(w_{1,+}(L_1) - w_{2,-}(L_1)= D_1 \) and (\romannumeral 8) \(w_{2,+}(L_1+2L_2) - e^{i\xi T}w_{1,-}(-L_1) = D_2\).

\begin{lemma}\label{lem2.3}
For fixed $l\geq 0$, there exist constants $C$ and $L_*$ such that the following statement holds for any \( L_1 + L_2 > L_* \):\\
(1) There exist analytic maps
\[A_2: [-\pi/T,\pi/T)
\rightarrow \mathbf{L}(\mathbb{C}^{2n} \times V_w, V_a),
\quad B_2: [-\pi/T,\pi/T) \rightarrow \mathbf{L}(\mathbb{C}^{2n} \times V_w, V_b),\]
\[W_2: [-\pi/T,\pi/T) \rightarrow \mathbf{L}(\mathbb{C}^{2n} \times V_w, V_w),\]
such that \( w \) satisfies (\ref{213}) if and only if \( w \) is given in (\ref{215}) with
\begin{equation}\label{220}
a = A_2(\xi)(D, G), \quad b = B_2(\xi)(D, G),
\end{equation}
that is,
\begin{equation}\label{221}
w = W_2(\xi)(D, G) = W_1(A_2(\xi)(D, G), B_2(\xi)(D, G), G),
\end{equation}
where $D= \left(D_1, D_2\right)$ and \( W_1 \) is defined in Lemma \ref{lem2.2}. \\
(2) The operators \( A_2 \), \( B_2 \), and \( W_2 \) satisfy
\begin{equation}\label{222}
\left| \partial_\xi^{\ell} A_2 \right| + \left| \partial_\xi^{\ell} B_2 \right| + \left| \partial_\xi^{\ell} W_2 \right| \leq C_l.
\end{equation}
(3) Furthermore, we have the expansion
\begin{equation}\label{223}
\begin{aligned}
a = (a_{1,-}, a_{1,+}, a_{2,-}, a_{2,+})  
  = (-e^{-i\xi T} P_{1}^s D_2, P_{2}^u D_1, -P_{2}^s D_1, P_{1}^u D_2) + A_3(\xi)(D,G),
  \end{aligned}
\end{equation}
where \( A_3 \)  satisfies
\begin{equation}\label{224}
\left| A_3(\xi)(D, G) \right| \leq C(e^{-\alpha L} |D| + |G|).
\end{equation}
\end{lemma}
\begin{proof}
Calculating equation (\ref{215})(\romannumeral 1)  at \( x = -L_1 \), (\romannumeral 2)(\romannumeral 3) at \( x = L_1 \), and  (\romannumeral 4) at \( x = L_1+2L_2 \), the results are given as follows
\begin{equation}\label{225}
\begin{aligned}
w_{1,-}(-L_1) &= a_{1,-} + \left( P_{1,-}^s( -L_1) - P^s_1 \right) a_{1,-} + \Phi^u_{1,-}( -L_1, 0) b_{1,-} \\
&\quad+ \int_0^{-L_1} \Phi^u_{1,-}(-L_1, y) G_{1,-}(y) dy, \\
w_{1,+}(L_1) &= a_{1,+} + \left( P_{1,+}^u(L_1) - P^u_2 \right) a_{1,+} + \Phi^s_{1,+}(L_1, 0) b_{1,+} \\
&\quad + \int_0^{L_1} \Phi^s_{1,+}(L_1, y) G_{1,+}(y) dy, \\
w_{2,-}(L_1) &= a_{2,-} + \left( P_{2,-}^s(-L_2) - P^s_2 \right) a_{2,-} + \Phi^u_{2,-}(-L_2, 0) b_{2,-}  \\
&\quad+ \int_{L_1+L_2}^{L_1} \Phi^u_{2,-}(-L_2, y-L_1-L_2) G_{2,-}(y) dy,  \\
w_{2,+}(L_1+2L_2) &= a_{2,+} + \left( P_{2,+}^u(L_2) - P^u_1 \right) a_{2,+} + \Phi^s_+(L_2, 0) b_{2,+} \\
&\quad+ \int_{L_1+L_2}^{L_1+2L_2} \Phi^s_{2,+}(L_2, y-L_1-L_2) G_{2,+}(y) dy,
\end{aligned}
\end{equation}
where the projections \( P^s_{i,-}(x) \) and \( P^u_{i,+}(x) \) 
are defined in Lemma \ref{lem2.1}. It is clear that \( P_1^s a_{1,-} = a_{1,-} \), \( P_2^u a_{1,+} = a_{1,+} \), \( P_2^s a_{2,-} = a_{2,-} \), and \( P_1^u a_{2,+} = a_{2,+} \), due to \( a \in V_a \).

Substituting these expressions into (\ref{213})(\romannumeral 7) and (\romannumeral 8), we obtain
\begin{equation}\label{226}
\begin{aligned}
D_1 =& w_{1,+}(L_1) - w_{2,-}(L_1) \\
=& a_{1,+} - a_{2,-} + \left( P_{1,+}^u(L_1) - P^u_2 \right) a_{1,+} - \left( P_{2,-}^s(-L_2) - P^s_2 \right) a_{2,-} \\
&+ \Phi^s_{1,+}(L_1, 0) b_{1,+}  - \Phi^u_{2,-}( -L_2, 0) b_{2,-} \\
& + \int_0^{L_1} \Phi^s_{1,+}(L_1, y) G_{1,+}(y) dy - \int_{L_1+L_2}^{L_1} \Phi^u_{2,-}(-L_2, y-L_1-L_2) G_{2,-}(y) dy, \\
D_2 =& w_{2,+}(L_1+2L_2) - e^{i\xi T} w_{1,-}(-L_1) \\
=& a_{2,+} - e^{i\xi T} a_{1,-} + \left( P_{2,+}^u(L_2) - P^u_1 \right) a_{2,+} - e^{i\xi T} \left( P_{1,-}^s( -L_1) - P^s_1 \right) a_{1,-} \\
& + \Phi^s_{2,+}(L_2, 0) b_{2,+}  - e^{i\xi T} \Phi^u_{1,-}( -L_1, 0) b_{1,-} \\
& + \int_{L_1+L_2}^{L_1+2L_2} \Phi^s_{2,+}(L_2, y-L_1-L_2) G_{2,+}(y) dy - e^{i\xi T} \int_0^{-L_1} \Phi^u_{1,-}(-L_1, y) G_{1,-}(y) dy,
\end{aligned}
 \end{equation}
where \( b = B_1(a, G) \) is defined in Lemma \ref{lem2.2} with \( |B_1(a, G)| \leq C(e^{-\alpha L} |a| + |G|) \). Equation (\ref{226}) can be rewritten in a more compact form
\begin{equation}\label{227}
\begin{pmatrix}
D_1 \\
D_2
\end{pmatrix}
=
\begin{pmatrix}
a_{1,+} - a_{2,-} \\
a_{2,+} - e^{i\xi T} a_{1,-}
\end{pmatrix}
+ A_1(\xi)(a, G),
\end{equation}
where \( A_1: [-\pi/T,\pi/T) \rightarrow \mathbf{L}(V_a \times V_w, \mathbb{C}^{2n}) \) is analytic. 
This equation will be solved for \( a \). Using the estimates from Lemma \ref{lem2.1} and Lemma \ref{lem2.2}, it can be verified directly that for every fixed \( \ell \geq 0 \), there exists a constant $C>0$ such that 
\begin{equation}\label{228}
 \left|\partial^\ell_{\xi} A_1(\xi)(a, G)\right|\leq C(e^{-\alpha L}|a|+|G|).
\end{equation}
The principal part of (\ref{227}) is given by the map \( J_1: [-\pi/T,\pi/T) \rightarrow \mathbf{L}(V_a, \mathbb{C}^{2n}) \), taking the form of 
\[J_1(\xi): V_a \rightarrow \mathbb{C}^{2n}, \quad a=(a_{1,-}, a_{1,+}, a_{2,-}, a_{2,+}) \mapsto \begin{pmatrix}
a_{1,+} - a_{2,-} \\
a_{2,+} - e^{i\xi T} a_{1,-}
\end{pmatrix}.\]
According to \( {\rm {R}}(P_1^u) \oplus {\rm {R}}(P_1^s) = \mathbb{C}^n \) and \( {\rm {R}}(P_2^u) \oplus {\rm {R}}(P_2^s) = \mathbb{C}^n \), the linear operator \( J_1(\xi) \) is an isomorphism for every \( \xi \).
Based on the fact that $A_1(\xi)(a,G)$ is linear with respect to $a$ and $G$, we can write
$$A_1(\xi)(a,G) = A_1(\xi)(a,0) + A_1(\xi)(0,G),$$
with
$$\left|\partial^\ell_{\xi} A_1(\xi)(a, 0)\right|\leq Ce^{-\alpha L}|a|. $$
Therefore, there exists the constant \( L_* > 0 \) such that, for every \( L_1 + L_2 > L_* \), the operator
\[a \mapsto J_1(\xi) a + A_1(\xi)(a, 0)\]
is invertible. Therefore, system (\ref{227}) can be solved abstractly, and its solution is given by
\begin{equation}\label{229}
a = (J_1(\xi) + A_1(\xi) I_1)^{-1} (D - A_1(\xi)(0, G)) =: A_2(\xi)(D, G), \end{equation}
where \( I_1 a = (a, 0) \). Note that \( A_2(\xi) \) is analytic in \( \xi \), linear in \( (D, G) \) and bounded uniformly in \( L_1 + L_2 > L_* \) and \( \xi \in [-\pi/T,\pi/T) \). 
Then, equations (\ref{213})(\romannumeral 7)(\romannumeral 8) remains to be solved. If we set
\begin{equation*}
    \begin{aligned}
B_2(\xi)(D, G) &:= B_1(A_2(\xi)(D, G), G), \\
W_2(\xi)(D, G) &:= W_1(A_2(\xi)(D, G), B_2(\xi)(D, G), G),
    \end{aligned}
\end{equation*}
according to the definition of \( B_1 \) and \( W_1 \), see Lemma \ref{lem2.2}, and the above construction of \( A_2 \), it can be concluded that for every \( (D, G) \), \( w = W_2(\xi)(D, G) \) meets the integral equation (\ref{215}).

After the entire boundary-value problem (\ref{213}) is solved, it is also necessary to verify the uniform estimate (\ref{222}), the expansion (\ref{223}), and the estimate (\ref{224}). It can be observed that, due to the boundedness of the operator \( (J_1(\xi) + A_1(\xi) I_1)^{-1} \) and the estimate (\ref{228}), we can derive
\begin{equation*}
\begin{aligned}
\left| A_2(D, G) \right| =  |(J_1(\xi) + A_1(\xi) I_1)^{-1}| |(D - A_1(\xi)(0, G))| \leq C(|D|+ |G| ),
\end{aligned}
\end{equation*}
and it can be easily verified that
\begin{equation}\label{230}
\left| \partial^\ell_{\xi} A_2 \right| \leq C.
\end{equation}
Utilizing this estimate along with (\ref{217}) and the definitions of \( B_2 \) and \( W_2 \),  it can be obtained that
\[\left| \partial^\ell_{\xi} B_2 \right| + \left| \partial^\ell_{\xi} W_2 \right| \leq C.\]
This proves the uniform estimate (\ref{222}). Applying projections $P_1^s, P_1^u$ and projections $P_2^s, P_2^u$ respectively to $D_1$ and $D_2$, given by (\ref{226}) can obtain
\begin{equation*}
    \begin{aligned}
a_{1,-} + e^{-i \xi T}P_1^s D_2 = & \left( P_{1,-}^s( -L_1) - P^s_1 \right) a_{1,-} + P_1^s \big( e^{-i \xi T}\Phi^s_{2, +}(L_2, 0) b_{2,+}  -  \Phi^u_{1,-}( -L_1, 0) b_{1,-} \\
& + e^{-i \xi T} \int_{L_1+L_2}^{L_1+2L_2} \Phi^s_{2,+}(L_2, y-L_1-L_2) G_{2,+}(y) dy  \\
& -  \int_0^{-L_1} \Phi^u_{1,-}(-L_1, y) G_{1,-}(y) dy \big) , \\
a_{1,+} - P_2^u D_1 = & - \left( P_{1,+}^u(L_1) - P^u_2 \right) a_{1,+} - P_2^u \big( \Phi^s_{1,+}(L_1, 0) b_{1,+}  - \Phi^u_{2,-}( -L_2, 0) b_{2,-} \\
& + \int_0^{L_1} \Phi^s_{1,+}(L_1, y) G_{1,+}(y) dy \\
& -  \int_{L_1+L_2}^{L_1} \Phi^u_{2,-}(-L_2, y-L_1-L_2) G_{2,-}(y) dy \big) , \\
a_{2,-} + P_2^s D_1 = & - \left( P_{2,-}^s(-L_2) - P^s_2 \right) a_{2,-} + P_2^s \big( \Phi^s_{1,+}(L_1, 0) b_{1,+}  - \Phi^u_{2,-}( -L_2, 0) b_{2,-} \\
& + \int_0^{L_1} \Phi^s_{1,+}(L_1, y) G_{1,+}(y) dy \\
& -  \int_{L_1+L_2}^{L_1} \Phi^u_{2,-}(-L_2, y-L_1-L_2) G_{2,-}(y) dy \big) , \\
a_{2,+} - P_1^u D_2 = & - \left( P_{2,+}^u(L_2) - P^u_1 \right) a_{2,+} - P_1^u \big( \Phi^s_{2,+}(L_2, 0) b_{2,+}  - e^{i\xi T} \Phi^u_{1,-}( -L_1, 0) b_{1,-} \\
& + \int_{L_1+L_2}^{L_1+2L_2} \Phi^s_{2,+}(L_2, y-L_1-L_2) G_{2,+}(y) dy  \\
& - e^{i\xi T} \int_0^{-L_1} \Phi^u_{1,-}(-L_1, y) G_{1,-}(y) dy \big) .
    \end{aligned}
\end{equation*}
The right-hand sides of these equations define a bounded and linear operator \( A_3 :  [-\pi/T,\pi/T) \rightarrow \mathbf{L}(\mathbb{C}^{2n} \times V_w, V_a) \) that satisfies
$$\left| A_3(\xi)(D, G) \right| \leq C(e^{-\alpha L} |D| + |G|).$$
The above estimate can be directly obtained by utilizing (\ref{228}). The proof is complete.
\end{proof}

\begin{prop} \label{prop2.1}
Based on the results of Lemma \ref{lem2.3}, we have
\begin{equation}\label{231}
\begin{aligned}
&(w_{1,-}(-L_1), w_{1,+}(L_1), w_{2,-}(L_1), w_{2,+}(L_1+2L_2)) \\
=& (-e^{-i\xi T} P_{1}^s D_2, P_{2}^u D_1, -P_{2}^s D_1, P_{1}^u D_2) + \widetilde{\mathcal{W}}_2(\xi)(D, G),
\end{aligned}
\end{equation}
where $\widetilde{\mathcal{W}}_2 : [-\pi/T,\pi/T) \rightarrow \mathbf{L}(\mathbb{C}^{2n} \times V_w, \mathbb{C}^{4n})$ is analytic in $\xi$ and, for fixed $\ell \geq 0$, we have
\begin{equation}\label{232}
\left|\partial^{\ell}_\xi \widetilde{\mathcal{W}}_2(\xi)(D, G)\right| \leq C(e^{-\alpha L}|D| + |G|).
\end{equation}

Furthermore, we have \( w = W_{2}(\xi)(D, 0) + W_{2}(\xi)(0, G) \) and
\begin{equation}\label{233}
    \begin{aligned}
\left| \partial^{\ell}_{\xi} W_{2}(\xi)(D, 0)(x) \right| &\leq C e^{-\alpha (L_1+x)} |D|, \quad \text{for } x \in [-L_1, 0]\\
\left| \partial^{\ell}_{\xi} W_{2}(\xi)(D, 0)(x) \right| &\leq C e^{-\alpha (L_1-x)} |D|, \quad \text{for } x \in [0, L_1]\\
\left| \partial^{\ell}_{\xi} W_{2}(\xi)(D, 0)(x) \right| &\leq C e^{-\alpha (x-L_1)} |D|, \quad \text{for } x \in [L_1, L_1+L_2]\\
\left| \partial^{\ell}_{\xi} W_{2}(\xi)(D, 0)(x) \right| &\leq C e^{-\alpha (L_1+2L_2-x)} |D|, \quad \text{for } x \in [L_1+L_2, L_1+2L_2].
    \end{aligned}
\end{equation}

\end{prop}

\begin{remark}
    This proposition is used to prove the second and third inequalities in \eqref{303}.
\end{remark}

\begin{proof}
By substituting \eqref{223} and \eqref{224} into \eqref{225}, it can be directly concluded that expression \eqref{231} and the corresponding estimate \eqref{232} hold.

The operator $W_2(\xi)$ defined by equation \eqref{221} is linear with respect to $D$ and $G$, where
$W_1(a,b,G)$ is defined by the right-hand side of the variation of constants formula \eqref{215}. The estimate (\ref{233}) can be obtained by setting $G = 0$ in equation \eqref{215} and applying the estimates from Lemma \ref{lem2.1} and Lemma \ref{lem2.2}. This concludes the proof.
\end{proof}

In Lemma \ref{lem2.3}, we have demonstrated that \( w \) is a bounded solution of (\ref{213})(\romannumeral 1)-(\romannumeral 8) if and only if, \( w \) is given by the variation-of-constant formula (\ref{215}), where \( (a, b) \) are defined by the bounded linear operators \( A_{2} \) and \( B_{2} \) appearing in (\ref{220}).  Recall that, having solved (\ref{213})(\romannumeral 5)(\romannumeral 6),
then \( w_{1,+}(0) = w_{2,-}(0) \) and $ w_{2,+}(L_1+L_2) = w_{2,-}(L_1+L_2)$ if and only if, the jumps
\begin{equation}\label{234}
    \begin{aligned}
        \Xi_1 &:= \langle \psi_{1}(0), w_{1,+}(0) - w_{1,-}(0) \rangle = 0, \\
        \Xi_2 &:= \langle \psi_{2}(0), w_{2,+}(L_1+L_2) - w_{2,-}(L_1+L_2)  \rangle = 0. \\
    \end{aligned}
\end{equation}
In the following lemma, we derive a concrete expression for the jumps $\Xi_{i}, i=1,2$.

\begin{lemma}\label{lem2.4}
For fixed $l\geq0$, there exist constants \( C \) and \( L_* \) such that the following statement holds for any \( L_1+L_2 > L_* \). Let \( w \) be given by (\ref{221}). The jumps \( \Xi_1 \) and \( \Xi_2 \) defined above are then given by
\begin{equation}\label{235}
    \begin{aligned}
     \Xi_1 =& \langle \psi_1(L_1), P_2^u D_1 \rangle + e^{-i\xi T} \langle \psi_1(-L_1), P_1^s D_2 \rangle \\
     &- \int_0^{L_1} \langle \psi_1(x), G_{1,+}(x) \rangle dx \\
     &- \int_{-L_1}^0 \langle \psi_1(x), G_{1,-}(x) \rangle dx + \hat{R}_{1,1}(\xi)(D, G),
    \end{aligned}
\end{equation}
and
\begin{equation}\label{236}
    \begin{aligned}
     \Xi_2 =& \langle \psi_2(L_2), P_1^u D_2 \rangle + \langle \psi_2(-L_2), P_2^s D_1 \rangle \\
     &- \int_{L_1+L_2}^{L_1+2L_2} \langle \psi_2(x-L_1-L_2), G_{2,+}(x) \rangle dx \\
     &- \int_{L_1}^{L_1+L_2} \langle \psi_2(x-L_1-L_2), G_{2,-}(x) \rangle dx + \hat{R}_{1,2}(\xi)(D, G),
    \end{aligned}
\end{equation}
for a certain analytic remainder term \( \hat{R}_{1} : [-\pi/T,\pi/T) \rightarrow \mathbf{L}(\mathbb{C}^{2n} \times V_w, \mathbb{C}^2) \) that satisfies
\begin{equation}\label{237}
\left| \partial^{\ell}_{\xi} \hat{R}_{1}(\xi)(D, G) \right| \leq C e^{-\alpha L} (e^{-\alpha L} |D| + |G|),
\end{equation}
 where $\hat{R}_{1} = (\hat{R}_{1,1}, \hat{R}_{1,2}).$
\end{lemma}

\begin{proof}
 Throughout the proof, we choose \( (a, b, w) \) according to (\ref{220}) and (\ref{221}) from Lemma \ref{lem2.3} so that (\ref{213}) is solved by the variation-of-constant formula (\ref{215}). Substituting the expressions for \( w_{1,\pm} \) and \( w_{2,\pm} \) from (\ref{215}) into equation (\ref{234}) for the jumps, we see that the jumps \( \Xi_{i} \) for $i=1,2$ are linear in \( (D, G) \) and analytic in \( \xi \).

 More precisely, taking the scalar product of equation (\ref{218}) for \( w_{1,+}(0) - w_{1,-}(0) \) with \( \psi_1(0) \) and using that \( \langle \psi_1(0), b_{1,+} \rangle = 0 \) and \( \langle \psi_1(0), b_{1,-} \rangle = 0 \) by definition of \( V_b \), we obtain
\begin{equation}\label{238}
    \begin{aligned}
&\langle \psi_1(0),  w_{1,+}(0) - w_{1,-}(0) \rangle \\
= &\langle \psi_1(0), \Phi_{1,+}^{u}(0, L_1) a_{1,+} \rangle - \langle \psi_1(0), \Phi_{1,-}^{s}(0, -L_1) a_{1,-} \rangle \\
& - \int_0^{L_1} \langle \psi_1(0), \Phi_{1,+}^{u}(0, x) G_{1,+}(x) \rangle dx - \int_{-L_1}^0 \langle \psi_1(0), \Phi_{1,-}^{s}(0, x) G_{1,-}(x) \rangle dx \\
= &\langle \psi_1(L_1), a_{1,+} \rangle - \langle \psi_1(-L_1), a_{1,-} \rangle \\
& - \int_0^{L_1} \langle \psi_1(x), G_{1,+}(x) \rangle dx - \int_{-L_1}^0 \langle \psi_1(x), G_{1,-}(x) \rangle dx.
    \end{aligned}
\end{equation}
In the second equality, we utilized  the fact that
\begin{equation*}
    \begin{aligned}
(\Phi_{1,+}^u(0, L_1))^* = \widetilde{\Phi}_{1,+}^s(L_1, 0),
\quad (\Phi_{1,-}^s(0, -L_1))^* = \widetilde{\Phi}_{1,-}^u(-L_1, 0).
    \end{aligned}
\end{equation*}
By substituting the expression (\ref{223}) from Lemma \ref{lem2.3} for \( a \), and using estimate (\ref{224}) for the remainder term, we obtain
\begin{equation}\label{239}
    \begin{aligned}
&\langle \psi_1(L_1), a_{1,+} \rangle - \langle \psi_1(-L_1), a_{1,-} \rangle \\
=& \langle \psi_1(L_1), P_2^u D_1 \rangle + e^{-i\xi T} \langle \psi_1(-L_1), P_1^s D_2 \rangle  + O(e^{-\alpha L} (e^{-\alpha L} |D| + |G|)),
    \end{aligned}
\end{equation}
Due to the complexity of the remainder term, we omit its specific form. We can directly obtain the estimate for the higher-order terms above, by using Proposition \ref{prop2.1} and the fact $ |\psi_{i}(x)| \leq Ce^{-\alpha |x|}, \quad  i=1,2.$

Similarly, by taking the scalar product of equation (\ref{219}) for \( w_{2,+}(L_1+L_2) - w_{2,-}(L_1+L_2) \) with \( \psi_2(0) \) and using  the fact that \( \langle \psi_2(0), b_{2,+} \rangle = 0 \) and \( \langle \psi_2(0), b_{2,-} \rangle = 0 \) by definition of \( V_b \), we obtain
\begin{equation}\label{241}
    \begin{aligned}
&\langle \psi_2(0), w_{2,+}(L_1+L_2) - w_{2,-}(L_1+L_2) \rangle \\
=& \langle \psi_2(0), \Phi_{2,+}^{u}(0, L_2) a_{2,+} \rangle - \langle \psi_2(0), \Phi_{2,-}^{s}(0, -L_2) a_{2,-} \rangle \\
& - \int_{L_1+L_2}^{L_1+2L_2} \langle \psi_2(0), \Phi_{2,+}^{u}(0, x-L_1-L_2) G_{2,+}(x) \rangle dx \\
& - \int_{L_1}^{L_1+L_2} \langle \psi_2(0), \Phi_{2,-}^{s}(0, x-L_1-L_2) G_{2,-}(x) \rangle dx \\
= &\langle \psi_2(L_2), a_{2,+} \rangle - \langle \psi_2(-L_2), a_{2,-} \rangle \\
& - \int_{L_1+L_2}^{L_1+2L_2} \langle \psi_2(x-L_1-L_2), G_{2,+}(x) \rangle dx - \int_{L_1}^{L_1+L_2} \langle \psi_2(x-L_1-L_2), G_{2,-}(x) \rangle dx,
    \end{aligned}
\end{equation}
and
\begin{equation}\label{242}
    \begin{aligned}
&\langle \psi_2(L_2), a_{2,+} \rangle - \langle \psi_1(-L_2), a_{2,-} \rangle \\
= &\langle \psi_2(L_2), P_1^u D_2 \rangle + \langle \psi_2(-L_2), P_2^s D_1 \rangle + O(e^{-\alpha L} (e^{-\alpha L} |D| + |G|)) .
    \end{aligned}
\end{equation}
This completes the proof of the lemma.
\end{proof}

\noindent{\bf proof of Lemma \ref{lem2.5}.}
Comparing the boundary-value problem (\ref{243}) with the equation (\ref{213}) reveals that it is necessary to set
\begin{equation} \label{249}
\begin{aligned}
G_{i,\pm} = H_i w_{i,\pm} + g_{i,\pm}, \quad i=1, 2
\end{aligned}
\end{equation}
in the unperturbed problem. Thus, \( w \) satisfies (\ref{243}) if and only if,
\[w = W_{2}(\xi)(D, Hw + g),\]
where \( W_{2}(\xi) \) is the solution operator for equation (\ref{213}), defined in (\ref{221}); see Lemma \ref{lem2.3}. Rewrite the above equation as
\begin{equation}\label{250}
w = W_{3}(\xi, H) w + W_{2}(\xi)(D, g),
\end{equation}
where \( W_{3}(\xi, H)w = W_{2}(\xi)(0, Hw) \) is analytic in \( (\xi, H) \) and
\[|W_{3}(\xi, H) w| \leq C |H| |w| \leq C \delta |w|,\]
by the estimate (\ref{222}).  Thus, for \( \delta > 0 \) sufficiently small but independent of \( \xi \) and \( L_1 +L_2 > L_{*} \), the operator \( \text{id} - W_{3}(\xi, H) \) can be inverted. Then we obtain a unique solution of (\ref{250}), which is  represented by the operator
 \begin{equation}\label{251}
w = W(\xi, H)(D, g) = (\text{id} - W_{3}(\xi, H))^{-1} W_{2}(\xi)(D, g).
\end{equation}
The operator $W(\xi, H)$ is analytic in \((\xi,H)\), and for fixed $\ell \geq 0$, it satisfies
\begin{equation*}
| \partial^{\ell}_{(\xi,H)} W | \leq C.
\end{equation*}
The expansion (\ref{245}) follows from Proposition \ref{prop2.1} and the estimates obtained earlier. This proves the first part of the lemma. Next, we calculate the jumps of the perturbed boundary-value problem. By substituting \( G = Hw + g \) into the expansion (\ref{235}) and (\ref{236}) for the jumps, as described in Lemma \ref{lem2.4}, we obtain
\begin{equation}\label{252}
\begin{aligned}
\Xi_1 =& \langle \psi_1(L_1), P_2^u D_1 \rangle + e^{-i\xi T} \left\langle \psi_1(-L_1), P_1^s D_2 \right\rangle \\
&- \int_0^{L_1} \langle \psi_1(x), H_1(x)w_{1,+}(x) \rangle dx
- \int_{-L_1}^0 \langle \psi_1(x), H_1(x)w_{1,-}(x) \rangle dx \\
&- \int_{-L_1}^{L_1} \langle \psi_1(x), g_1(x) \rangle dx   + \hat{R}_{3,1}(\xi , H)(D, g),
\end{aligned}
\end{equation}
and,
\begin{equation}\label{253}
    \begin{aligned}
     \Xi_2 =& \langle \psi_2(L_2), P_1^u D_2 \rangle + \langle \psi_2(-L_2), P_2^s D_1 \rangle \\
     &- \int_{L_1+L_2}^{L_1+2L_2} \langle \psi_2(x-L_1-L_2), H_2(x)w_{2,+}(x) \rangle dx \\
     &- \int_{L_1}^{L_1+L_2} \langle \psi_2(x-L_1-L_2), H_2(x)w_{2,-}(x) \rangle dx \\
     &- \int_{L_1}^{L_1+2L_2} \langle \psi_2(x-L_1-L_2), g_2(x) \rangle dx + \hat{R}_{3,2}(\xi, H )(D, g),
    \end{aligned}
\end{equation}
where \( w = W(\xi, H)(D, g) \) and
\begin{equation}\label{254}
    \begin{aligned}
&\hat{R}_{3}: [-\pi/T,\pi/T) \times U_{\delta} \rightarrow \mathbf{L}(\mathbb{C}^{2n} \times V_{w}, \mathbb{C}^2),\\
&\hat{R}_{3}(\xi, H)(D, g) := \hat{R}_{1}(\xi)(D, Hw + g),
    \end{aligned}
\end{equation}
 is analytic, where $\hat{R}_{3}= (\hat{R}_{3,1}, \hat{R}_{3,2})$. From (\ref{250}), it can be deduced that
\[w = W_{2}(\xi)(D, g) + W_{3}(\xi, H) W(\xi, H)(D, g).\]

Set
\begin{align*}
\hat{R}_{4, 1}(\xi, H)(D, g) = & \int_0^{L_1} \langle \psi_1(x), H_1(x)w_{1,+}(x) \rangle dx + \int_{-L_1}^0 \langle \psi_1(x), H_1(x)w_{1,-}(x) \rangle dx,
\end{align*}
then, we know that \( R_{4,1} \) is analytic in \( (\xi, H) \) and
\begin{equation} \label{255}
\hat{R}_{4,1}(\xi, H)(D, g) = \hat{T}_{1,1}(\xi)(D, g)[H] + \hat{T}_{2,1}(\xi, H)(D, g)[H, H],
\end{equation}
where
\begin{align*}
\hat{T}_{1,1}(\xi)(D, g)[H] =& \int_{0}^{L_1} \langle \psi_{1}(x), H_1(x)W_{2}(\xi)(D, g)(x) \rangle dx \\
& + \int_{-L_1}^{0} \langle \psi_{1}(x), H_1(x)W_{2}(\xi)(D, g)(x) \rangle dx ,\\
\hat{T}_{2,1}(\xi, H)(D, g)[H, H] = &\int_{0}^{L_1} \langle \psi_{1}(x), H_1(x)W_{3}(\xi, H) W(\xi, H)(D, g)(x) \rangle dx \\
& + \int_{-L_1}^{0} \langle \psi_{1}(x), H_1(x)W_{3}(\xi, H) W(\xi, H)(D, g)(x) \rangle dx.
\end{align*}
Based on the estimates for \( W_{3} \) and \( W \), it follows that 
 \[\left| \partial^{\ell}_{(\xi, H)} \hat{T}_{2,1}(\xi, H) \right| \leq C.\]
 Furthermore, using the estimate (\ref{233}) in Proposition \ref{prop2.1},
we obtain
\begin{align*}
 \left| \int_{0}^{L_1} \langle \psi_{1}(x), H_1(x) W_{2}(\xi)(D, 0)(x) \rangle dx \right|
\leq &\int_{0}^{L_1} C \left| \psi_{1}(x) \right| \left| H_1 \right| e^{-\alpha(L_1+x)} |D| dx \\
\leq &Ce^{-\alpha L} |H| |D|,
\end{align*}
and similarly
\begin{align*}
\left| \int_{-L_1}^{0} \langle \psi_{1}(x), H_1(x) W_{2}(\xi)(D, 0)(x) \rangle dx \right| \leq Ce^{-\alpha L} |H| |D|.
\end{align*}

We conclude that
\[\left|\partial^{\ell}_{\xi} T_{1,1}(\xi)(D, g)\right| \leq C(e^{-\alpha L} |D| + |g|).\]
 In summary, we have
 \begin{align*}
\Xi_1 = &\langle \psi_{1}(L_1), P_{2}^u D_1 \rangle + e^{-i \xi T} \langle \psi_{1}(-L_1), P_{1}^s D_2 \rangle \\
&- \int_{-L_1}^{L_1} \langle \psi_{1}(x), g_{1}(x) \rangle dx + \hat{R}_{2, 1}(\xi, H )(D, g),
  \end{align*}
where
\[
\hat{R}_{2,1}(\xi, H) = \hat{R}_{3,1}(\xi, H) + \hat{R}_{4,1}(\xi, H).
\]
This last expression can be written as
\[
\hat{R}_{2,1}(\xi, H)(D, g) = T_{0,1}(\xi)(D, g) + T_{1,1}(\xi)(D, g)[H] + T_{2,1}(\xi, H)(D, g)[H, H],
\]
where
\begin{align*}
T_{0,1}(\xi)(D, g) =& \hat{R}_{1,1}(\xi)(D, g),\\
T_{1,1}(\xi)(D, g)[H] =& \hat{R}_{1,1}(\xi)(0, HW(\xi, 0)(D, g)) + \hat{T}_{1,1}(\xi)(D, g)[H],\\
T_{2,1}(\xi, H)(D, g)[H, H] =& \hat{R}_{1,1}(\xi)(0, HW(\xi, H)(D, g) - HW(\xi, 0)(D, g)) \\
&+ \hat{T}_{2,1}(\xi)(D, g)[H, H],
\end{align*}
are analytic; see (\ref{254}) and (\ref{255}). It follows from the estimates for \(\hat{T}_{1,1}\) and \(\hat{T}_{2,1}\) as well as (\ref{237}) that
\begin{equation*}
\begin{aligned}
\left| \partial^{\ell}_{\xi} T_{0,1}(\xi)(D, g) \right| &\leq C e^{-\alpha L}(\left|D\right| + \left|g\right|),\\
\left| \partial^{\ell}_{\xi} T_{1,1}(\xi)(D, g) \right| &\leq C (e^{-\alpha L}\left|D\right| + \left|g\right|),\\
\left| \partial^{\ell}_{(\xi, H)} T_{2,1}(\xi, H) \right| &\leq C.
\end{aligned}
\end{equation*}

By repeating the above argument, we obtain
  \begin{align*}
\Xi_2 =& \langle \psi_{2}(L_2), P_{1}^u D_2 \rangle + \langle \psi_{2}(-L_2), P_{2}^s D_1 \rangle \\
&- \int_{L_1}^{L_1+2L_2} \langle \psi_{2}(x-L_1-L_2), g_{2}(x) \rangle dx + \hat{R}_{2, 2}(\xi, H )(D, g),
  \end{align*}
where the remainder term
\[
\hat{R}_{2,2}(\xi, H) = T_{0,2}(\xi)(D, g) + T_{1,2}(\xi)(D, g)[H] + T_{2,2}(\xi, H)(D, g)[H, H]
\]
also satisfies the estimate (\ref{248}).
This completes the proof of the lemma.

\section{Proof of the existence of $\varphi$}\label{B}\renewcommand{\theequation}{B.\arabic{equation}}
\setcounter{equation}{0}
\begin{lemma}\label{V}
  For $\mathbf{\widetilde{V}}$ given by Proposition \ref{tildeV}, the mapping $V:H^2(\mathbb{R})\times [0, T_{\rm{max}})\to H^2(\mathbb{R})\times H^2(\mathbb{R})$ given by
  $$V(\varphi,t)[x]=\mathbf{\widetilde{V}}(x-\varphi(x),t)+\mathbf{\bar{U}}(x-\varphi(x))-\mathbf{\bar{U}}(x)$$
  is well-defined, continuous in $t$, and locally Lipschitz continuous in $\varphi$ (uniformly in $t$ on compact subintervals of $[0,T_{\rm{max}})$).
\end{lemma}
\begin{proof}
Since $\mathbf{\widetilde{V}}\in C([0,T_{\rm{max}}), H^4(\mathbb{R}))$, we have
\begin{align}\label{CB3}\mathbf{\widetilde{V}}\in C([0,T_{\rm{max}}), C_b^3(\mathbb{R})),\end{align}
by virtue of the embedding $H^4(\mathbb{R})\hookrightarrow C_b^3(\mathbb{R})$.
Then applying the mean value theorem, we obtain the inequality
\begin{equation}\label{V2}\begin{split}
\|V(\varphi_1,t)-V(\varphi_2,t)\|_{H^2}&\lesssim \|\mathbf{\widetilde{V}}(t)+\mathbf{\bar{U}}\|_{W^{3,\infty}}\|\varphi_1-\varphi_2\|_{H^2}\\
&\lesssim \|\mathbf{\widetilde{V}}(t)+\mathbf{\bar{U}}\|_{H^4}\|\varphi_1-\varphi_2\|_{H^2}
\end{split}
\end{equation}for $\varphi_{1,2}\in H^2(\mathbb{R})$ and $t\in[0,T_{\rm{max}})$. Thus, $V$ is locally Lipschitz continuous in $\varphi$ (uniformly in $t$ on compact subintervals of $[0,T_{\rm{max}})$). Moreover,   the inequality \eqref{V2} also shows that $V$ is well-defined by taking $\varphi_2=0$ and noting $V(0,t)=\mathbf{\widetilde{V}}(t)\in H^2(\mathbb{R})$.
By virtue of \eqref{CB3}, we have the equality $$\mathbf{\widetilde{V}}(x_1,t)=\mathbf{\widetilde{V}}(x_2,t)+\int_{x_2}^{x_1}\mathbf{\widetilde{V}}_x(y,t)dy, \text{ for fixed }t\in[0,T_{\rm{max}}).$$
Then, we arrive at
\begin{equation}\begin{split}
\|V(\varphi,t)-V(\varphi,s)\|_{L^2}\leq\ &
\|(V(\varphi,t)-V(0,t))-(V(\varphi,s)-V(0,s))\|_{L^2}+\|(V(0,t)-V(0,s))\|_{L^2}\\
=\ &\|\mathbf{\widetilde{V}}(x-\varphi(x),t)-\mathbf{\widetilde{V}}(t,x)-(\mathbf{\widetilde{V}}(x-\varphi(x),s)-\mathbf{\widetilde{V}}(x,s))\|_{L^2}\\
&+\|(V(0,t)-V(0,s))\|_{L^2}\\
\leq\ & \|\mathbf{\widetilde{V}}(t)-\mathbf{\widetilde{V}}(s)\|_{L^\infty}\|\varphi\|_{L^2}+\|\mathbf{\widetilde{V}}(t)-\mathbf{\widetilde{V}}(s)\|_{L^2}.
\end{split}
\end{equation}
A similar calculation shows that
\begin{align*}
\|V(\varphi,t)-V(\varphi,s)\|_{H^2}\leq\|\mathbf{\widetilde{V}}(t)-\mathbf{\widetilde{V}}(s)\|_{W^{3,\infty}}\|\varphi\|_{L^2}+\|\mathbf{\widetilde{V}}(t)-\mathbf{\widetilde{V}}(s)\|_{H^2},
\end{align*}
which implies that $V$ is continuous in $t$.
\end{proof}

\begin{prop}\label{Vphi}
For $\mathbf{\widetilde{V}}$ given by Proposition \ref{tildeV}, let $V:H^2(\mathbb{R})\times [0, T_{\rm{max}})\to H^2(\mathbb{R})\times H^2(\mathbb{R})$ be the mapping in Lemma \ref{V}. Then there exists a maximal time $\tau_{\rm{max}}\in(0,T_{\rm{max}}]$ such that system
\begin{equation}\label{phiphit}\begin{split}
\varphi(t)&=s_p(t)\mathbf{V}_0+\int_0^ts_p(t-s)\mathcal{N}(V(\varphi(s),s),\varphi(s),\varphi_t(s))ds,\\
\varphi_t(t)&=\partial_ts_p(t)\mathbf{V}_0+\int_0^t\partial_ts_p(t-s)\mathcal{N}(V(\varphi(s),s),\varphi(s),\varphi_t(s))ds,
\end{split}
\end{equation}has a unique solution
$$(\varphi,\varphi_t)\in C([0,\tau_{\rm{max}}),H^3({\mathbb{R}})\times H^1(\mathbb{R})).$$
Furthermore, if $\tau_{\rm{max}}<T_{\rm{max}}$, then
\begin{align}\label{blowup}
\lim_{t\uparrow \tau_{\rm{max}}}\|(\varphi,\varphi_t)\|_{H^3\times H^1}=\infty.
\end{align}
\end{prop}
\begin{proof}
By virtue of Lemma \ref{modun} and Lemma \ref{V}, the nonlinear map $N: H^3(\mathbb{R})\times H^1(\mathbb{R})\times [0,T_{\rm{max}})\to L^2(\mathbb{R})$ given by
$$N(\varphi,\varphi_t,t)=\mathcal{N}(V(\varphi,t),\varphi,\varphi_t)$$
is well-defined, continuous in $t$, and locally Lipschitz continuous in $(\varphi,\varphi_t)$ (uniformly in $t$ on compact subintervals of $[0,T_{\rm{max}})$), where we use the embedding inequality $\|g\|_{L^\infty}\lesssim \|g\|_{H^1}$.
Define
$$\Psi(\varphi,\varphi_t)[t]=\left(\begin{array}{cc}s_p(t)\mathbf{V}_0\\\partial_ts_p(t)\mathbf{V}_0\end{array}\right)+\int_0^t\left(\begin{array}{cc}s_p(t-s)\mathcal{N}(V(\varphi(s),s),\varphi(s),\varphi_t(s))\\\partial_ts_p(t-s)\mathcal{N}(V(\varphi(s),s),\varphi(s),\varphi_t(s))\end{array}\right)ds.$$
Note that the operators $s_p(t): L^2(\mathbb{R})\times L^2(\mathbb{R})\to H^3(\mathbb{R})$ amd $\partial_ts_p(t): L^2(\mathbb{R})\times L^2(\mathbb{R})\to H^1(\mathbb{R})$ are uniformly bounded with respect to $t$ and strong continuous on $[0,\infty)$. Then, by standard arguments (see Proposition 4.3.9 of \cite{CH98} or Theorem 6.1.3 of \cite{pazy}), for $\|\mathbf{V}_0\|_{L^2}$ sufficiently small, there exist constants $R>0$ and $\tau\in(0,T_{\rm{max}}]$ such that
$$\Psi: C([0,\tau],B(R))\to C([0,\tau],B(R)),$$
is a well-defined contraction mapping, where $B(R)=\{(g,h)\in H^3({\mathbb{R}})\times H^1({\mathbb{R}}): \|(g,h)\|_{H^3\times H^1}\leq R\}$. Thus, by Banach fixed point theorem, $\Psi$ has a unique fixed point, which yields a unique solution $(\varphi,\varphi_t)\in C([0,\tau],H^3({\mathbb{R}})\times H^1(\mathbb{R}))$ of \eqref{phiphit}. Suppose that $\tau_{\rm{max}}\in(0,T_{\rm{max}}]$ is the supremum of all such $\tau$. Then we obtain a maximally defined solution of \eqref{phiphit} $(\varphi,\varphi_t)\in C([0,\tau_{\rm{max}}),H^3({\mathbb{R}})\times H^1(\mathbb{R}))$.

Next, we assume that $\tau_{\rm{max}}<T_{\rm{max}}$ and \eqref{blowup} do not hold. For any fixed $t_0$, we can prove that there exist constants $R_1,\delta>0$ (independent of $t_0$) such that $\Psi_{t_0}: C([t_0,t_0+\tilde{\delta}],B(R_1))\to C([t_0,t_0+\tilde{\delta}],B(R_1))$ defined as
\begin{align*}
\Psi_{t_0}(\tilde{\varphi},\tilde{\varphi}_t)[t]=&\left(\begin{array}{cc}s_p(t)\mathbf{V}_0\\\partial_ts_p(t)\mathbf{V}_0\end{array}\right)+\int_0^t\left(\begin{array}{cc}s_p(t-s)\mathcal{N}(V(\varphi(s),s),\varphi(s),\varphi_t(s))\\\partial_ts_p(t-s)\mathcal{N}(V(\varphi(s),s),\varphi(s),\varphi_t(s))\end{array}\right)ds\\
&+\int_{t_0}^{t_0+\tilde{\delta}}\left(\begin{array}{cc}s_p(t-s)\mathcal{N}(V(\tilde{\varphi}(s),s),\tilde{\varphi}(s),\tilde{\varphi}_t(s))\\\partial_ts_p(t-s)\mathcal{N}(V(\tilde{\varphi}(s),s),\tilde{\varphi}(s),\tilde{\varphi}_t(s))\end{array}\right)ds
\end{align*}
is a well-defined contraction mapping. Then $\Psi_{t_0}$ has a unique solution $(\tilde{\varphi},\tilde{\varphi}_t)\in C([t_0,t_0+\tilde{\delta}],H^3({\mathbb{R}})\times H^1(\mathbb{R}))$. Choose $t_0=\tau_{\rm{max}}-\tilde{\delta}/2$ and define
$$(\check{\varphi}(t),\check{\varphi}(t))=\left\{\begin{array}{ll}(\varphi(t),\varphi_t(t)),\quad t\in[0,t_0],\\(\tilde{\varphi}(t),\tilde{\varphi}_t(t)),\quad t\in[t_0,t_0+\tilde{\delta}].\end{array}\right.$$
Then $(\check{\varphi}(t),\check{\varphi}(t))\in C([0,\tau_{\rm{max}}],H^3({\mathbb{R}})\times H^1(\mathbb{R}))$ is a solution of \eqref{phiphit}, which contradicts the maximality of $\tau_{\rm{max}}$. Then we complete the proof.
\end{proof}

\medskip

\end{document}